\documentclass[12pt]{amsbook}

\usepackage{amsmath,amsthm, amssymb}
\usepackage[T1]{fontenc}
\usepackage{latexsym}
\usepackage{xspace}
\usepackage{wasysym}
\usepackage{setspace}
\usepackage{bbm}
\usepackage{stmaryrd}
\usepackage{todonotes}

\newcommand{\NN}{\mathbb{N}}  \newcommand{\ZZ}{\mathbb{Z}}  \newcommand{\QQ}{\mathbb{Q}}  
\newcommand{\RR}{\mathbb{R}}    \newcommand{\LL}{\mathbb{L}}

\newcommand{\bB}{\mathbb{B}}  \newcommand{\DD}{\mathbb{D}}  \newcommand{\XX}{\mathbb{X}}

 \newcommand{\mcC}{\mathcal{C}}  \newcommand{\mcF}{\mathcal{F}}

         \newcommand{\mcA}{\mathcal{A}}

  \newcommand{\be}{{\bf e}} \newcommand{\bfa}{{\bf a}}  \newcommand{\bfb}{{\bf b}}

    \newcommand{\bfX}{{\bf X}}  \newcommand{\bfY}{\bf Y}   \newcommand{\bfB}{\bf B} 
\newcommand{\bfZ}{\bf Z}     \newcommand{\bfP}{\bf P} 
 
\newcommand{\bbX}{\mathbb{X}}    \newcommand{\bbA}{\mathbb{A}}   \newcommand{\bbZ}{\mathbb{Z}}     \newcommand{\bbB}{\mathbb{B}}

\newcommand{\EE}{\mathbb{E}} \newcommand{\PP}{\mathbb{P}}

\newcommand{\st}{such that }

\newcommand{\varep}{\varepsilon}
\newcommand{\ep}{\epsilon}
\renewcommand{\leq}{\leqslant}
\renewcommand{\geq}{\geqslant}

\newcommand{\al}{\alpha}

\newcommand{\wrt}{with respect to }
\renewcommand{\st}{such that }

\newcommand{\ssk}{\smallskip}

\newtheorem{thm}{\hspace{-0.15cm}  {\sc Theorem} }
\newtheorem{cor}[thm]{\hspace{-0.15cm}  {\sc Corollary} }
\newtheorem{lem}[thm]{\hspace{-0.15cm}  {\sc Lemma} }
\newtheorem{prop}[thm]{\hspace{-0.15cm} {\sc Proposition}}
\newtheorem{defn}[thm]{ \hspace{-0.3cm} {\sc Definition}}

\newtheorem{Rems}[thm]{{\sc Remarks}}

\numberwithin{equation}{section} 

\newenvironment{Dem}{%
    \begin{list}{\hspace{0.5cm}{\sc Proof --}}{%
        \setlength{\topsep}{0pt}%
        \setlength{\leftmargin}{0pt}%
        \setlength{\rightmargin}{0pt}%
        \setlength{\listparindent}{0pt}%
        \setlength{\itemindent}{0pt}%
        \setlength{\parsep}{0pt}%
        \addtolength{\leftmargin}{20pt}%
        \addtolength{\rightmargin}{0pt}%
    } \item }{\hfill{\space $\rhd$}\end{list}\smallskip}

    \newenvironment{DemPropEstimatesRDEWarmUp}{%
    \begin{list}{\hspace{0.5cm}{\sc Proof of proposition \ref{PropEstimatesRDEWarmUp} --}}{%
        \setlength{\topsep}{0pt}%
        \setlength{\leftmargin}{0pt}%
        \setlength{\rightmargin}{0pt}%
        \setlength{\listparindent}{0pt}%
        \setlength{\itemindent}{0pt}%
        \setlength{\parsep}{0pt}%
        \addtolength{\leftmargin}{20pt}%
        \addtolength{\rightmargin}{0pt}%
    } \item }{\hfill{\space $\rhd$}\end{list}\smallskip}

    \newenvironment{DemPropFundamentalEstimate}{%
    \begin{list}{\hspace{0.5cm}{\sc Proof of proposition \ref{PropFundamentalEstimate} --}}{%
        \setlength{\topsep}{0pt}%
        \setlength{\leftmargin}{0pt}%
        \setlength{\rightmargin}{0pt}%
        \setlength{\listparindent}{0pt}%
        \setlength{\itemindent}{0pt}%
        \setlength{\parsep}{0pt}%
        \addtolength{\leftmargin}{20pt}%
        \addtolength{\rightmargin}{0pt}%
    } \item }{\hfill{\space $\rhd$}\end{list}\smallskip}

      \newenvironment{DemThmSchilderBrownianRP}{%
    \begin{list}{\hspace{0.5cm}{\sc Proof of theorem \ref{ThmSchilderBrownianRP} --}}{%
        \setlength{\topsep}{0pt}%
        \setlength{\leftmargin}{0pt}%
        \setlength{\rightmargin}{0pt}%
        \setlength{\listparindent}{0pt}%
        \setlength{\itemindent}{0pt}%
        \setlength{\parsep}{0pt}%
        \addtolength{\leftmargin}{20pt}%
        \addtolength{\rightmargin}{0pt}%
    } \item }{\hfill{\space $\rhd$}\end{list}\smallskip}

\title{A flow-based approach to rough differential equations}
\date{\today}
\author{I. Bailleul\footnote{IRMAR, 263 Avenue du General Leclerc, 35042 RENNES, France, {\sf ismael.bailleul@univ-rennes1.fr}}}

\begin{document}

\maketitle

\tableofcontents

\chapter{Introduction}
\label{SectionIntroduction}

This course is dedicated to the study of some class of dynamics in a Banach space, index by time $\RR_+$. Although there exists many recipes to cook up such dynamics, those generated by differential equations or vector fields on some configuration space are the most important from a historical point of view. Classical mechanics reached for example its top with the description by Hamilton of the evolution of any classical system as the solution of a first order differential equation with a universal form.  The outcome, in the second half of the twentieth centary, of the study of random phenomena did not really change that state of affair, with the introduction by It\^o of stochastic integration and stochastic differential equations.

\ssk

Classically, one understands a differential equation as the description of a point motion, the set of all these motions being gathered into a single object called a \textit{flow}.\index{Flow} It is a familly $\varphi = \big(\varphi_{ts}\big)_{0\leq s\leq t\leq T}$ of maps from the state space to itself, such that $\varphi_{tt} = \textrm{Id}$, for all $0\leq t\leq T$, and $\varphi_{ts} = \varphi_{tu}\circ\varphi_{us}$, for all $0\leq s\leq u\leq t\leq T$. The first aim of the approach to some class of dynamics that is proposed is this course is the construction of flows, as opposed to the construction of trajectories started from some given point.

\ssk

I will explain in the first part of the course a simple method for constructing a flow $\varphi$ from a family $\mu = \big(\mu_{ts}\big)_{0\leq s \leq t\leq T}$ of maps that almost forms a flow. The two essential points of this construction are that
\begin{itemize}
   \item[{\bf i)}] $\varphi_{ts}$ is loosely speaking the composition of infinitely many $\mu_{t_{i+1}t_i}$ along an infinite partition $s<t_1<\cdots<t$ of the interval $[s,t]$, with infinitesimal mesh,
   \item[{\bf ii)}] $\varphi$ depends continuously on $\mu$ in some sense.
\end{itemize}
Our main application of this general machinery will be to study some general class of controlled ordinary differential equations, that is differential equations of the form
$$
dx_t = \sum_{i=1}^\ell V_i(x_t)dh^i_t,
$$
where the $V_i$ are vector fields on $\RR^d$, say, and the controls $h^i$ are real-valued. Giving some meaning and solving such an equation in some general framework is highly non-trivial outside the framework of absolutely continuous controls, without any extra input like probability, under the form of stochastic calculus for instance. It requires Young integration theory for controls with finite $p$-variation, for $1\leq p<2$, and Terry Lyons' theory of rough paths for "rougher" controls! Probabilists are well-acquainted with this kind of situation as stochastic differential equations driven by some Brownian motion are nothing but an example of the kind of problem we intend to tackle. (With no probability!) It is the aim of this course to give you all the necessary tools to understand what is going on here, in the most elementary way as possible, while aiming at some generality. 

\bigskip

The general machinery of approximate flows is best illustrated by looking at the classical Cauchy-Lipschitz theory.Fix some Lipschitz continuous vector fields $V_i$ and some real-valued controls $h^i$ of class $\mcC^1$. It will appear in our setting that a good way of understanding what it means to be a solution to the ordinary differential equation on $\RR^n$ 
\begin{equation}
\label{EqModelODE}
\dot x_t = \sum_{i=1}^\ell V_i(x_t)\dot h^i_t =:  V_i(x_t)\dot h^i_t,
\end{equation}
is to say that the path $x_\bullet$ satisfies at any time $s$ the Taylor-type expansion formula
$$
x_t = x_s + \big(h^i_t-h^i_s\big)V_i(x_s) + o(t-s),
$$
and even 
$$
f\big(x_t\big) = f\big(x_s\big) + \big(h^i_t-h^i_s\big)\big(V_if\big)(x_s) + O\big(|t-s|^2\big),
$$
for any function $f$ of class $\mcC^2_b$, with $V_if$ standing for the derivative of $f$ in the direction of $V_i$. Setting $\mu_{ts}(x) := x+\big(h^i_t-h^i_s\big)V_i(x)$, the preceeding identity rewrites
$$
f\big(x_t\big) = f\big(\mu_{ts}(x_s)\big) + O\big(|t-s|^2\big),
$$
so the elementary map $\mu_{ts}$ provides a very accurate description of the dynamics. It almost forms a flow under mild regularity assumptions on the driving vector fields $V_i$, and its flow associated by the above "almost-flow to flow" machinery happens to be flow classically generated by  equation \eqref{EqModelODE}.

\ssk

Going back to a probabilistic setting, what insight does this machinery provide on Stratonovich stochastic differential equations 
\begin{equation}
\label{EqStratoSDE}
{\circ d}x_t = V_i(x_t) {\circ d}w_t
\end{equation}
driven by some Brownian motion $w$? The use of this notion of differential enables to write the following kind of Taylor-type expansion of order 2 for any function $f$ of class $\mcC^3$.
\begin{equation}
\begin{split}
f\big(x_t\big) &= f\big(x_s\big) + \int_s^t \big(V_if\big)(x_r)\,{\circ d}w_r \\
                      &= f\big(x_s\big) + \big(w^i_t-w^i_s\big)\big(V_if\big)(x_s) + \int_s^t\int_s^r \big(V_j(V_if)\big)(x_u)\,{\circ d}w_u\,{\circ d}w_r \\
                      &= f\big(x_s\big) + \big(w^i_t-w^i_s\big)\big(V_if\big)(x_s) + \left(\int_s^t\int_s^r {\circ d}w_u\,{\circ d}w_r\right)\big(V_j(V_if)\big)(x_s) + \int_s^t\int_s^r\int_s^u (\cdots)
\end{split}
\end{equation}
For any choice of $2<p<3$, the Brownian increments $w^i_{ts} := w^i_t-w^i_s$ have almost-surely a size of order $(t-s)^\frac{1}{p}$, the iterated integrals $\int_s^t\int_s^r {\circ d}w_u\,{\circ d}w_r$ have size $(t-s)^\frac{2}{p}$, and the triple integral size $(t-s)^\frac{3}{p}$, with $\frac{3}{p}>1$. What will come later out of this formula is that a solution to equation \eqref{EqStratoSDE} is precisely a path $x_\bullet$ for which one can write for any function $f$ of class $\mcC^3$ a Taylor-type expansion of order 2 of the form
$$
f\big(x_t\big) =  f\big(x_s\big) + \big(w^i_t-w^i_s\big)\big(V_if\big)(x_s) + \left(\int_s^t\int_s^r {\circ d}w_u\,{\circ d}w_r\right)\big(V_j(V_if)\big)(x_s) + O\big(|t-s|^a\big)
$$
at any time $s$, for some exponent $a>1$ independent of $s$. This conclusion puts forward the fact that what the dynamics really see of the Brownian control $w$ is not only its increments  $w_{ts}$ but also its iterated integrals $\int_s^t\int_s^r {\circ d}w_u\,{\circ d}w_r$. The notion of $p$-rough path ${\bfX} = \big(X_{ts},\XX_{ts}\big)_{0\leq s\leq t\leq T}$ is an abstraction of this family of pairs of quantities, for $2<p<3$ here. This \textit{multi-level object} satisfies some constraints of analytic type (size of its increments) and algebraic type, coming from the higher level parts of the object. As they play the role of some iterated integrals, they need to satisfy some identities consequences of the Chasles relation for elementary integrals: $\int_s^t = \int_s^u + \int_u^t$. These constraints are all what these rough paths ${\bfX} = (X,\XX)$ need to satisfy to give a sense to the equation
\begin{equation}
\label{EqRDE}
dx_t = \textrm{F}^\otimes (x_t){\bfX}(dt)
\end{equation}
for a collection $\textrm{F} = \big(V_1,\dots,V_\ell\big)$ of vector fields on $\RR^n$, by defining a solution as a path $x_\bullet$ for which one can write some uniform Taylor-type expansion of order 2
\begin{equation}
\label{EqSolRDE}
f\big(x_t\big) =  f\big(x_s\big) + X^i_{ts}\big(V_if\big)(x_s) + \bbX^{jk}_{ts}\big(V_j(V_kf)\big)(x_s) + O\big(|t-s|^a\big),
\end{equation}
for any function $f$ of class $\mcC^3_b$. The notation $\textrm{F}^\otimes$ is used here to insist on the fact that it is not only the collection F of vector fields  that is used in this definition, but also the differential operators $V_jV_k$ constructed from F. The introduction and the study of $p$-rough paths and their collection is done in the second part of the course.

\ssk

Guided by the results on flows of the first part, we shall reinterpret equation \eqref{EqRDE} to construct directly a flow $\varphi$ solution to the equation
\begin{equation}
\label{EqFlowRDE}
d\varphi = \textrm{F}^\otimes{\bfX}(dt),
\end{equation}
in a sense to be made precise in the third part of the course. The recipe of construction of $\varphi$  will consist in associating to F and $\bfX$ a $\mcC^1$-approximate flow $\mu = \big(\mu_{ts}\big)_{0\leq s\leq t\leq T}$ having everywhere a behaviour similar to that described by equation \eqref{EqSolRDE}, and then to apply the theory described in the first part of the course. The maps $\mu_{ts}$ will be constructed as the time 1 maps associated with some ordinary differential equation constructed from F and ${\bfX}_{ts}$ in a simple way. As they will depend continuously on $\bfX$, the continuous dependence of $\varphi$ on $\bfX$ will come as a consequence of point {\bf ii)} above.

\medskip

All that will be done in a deterministic setting. We shall see in the fourth part of the course how this approach to dynamics is useful in giving a fresh viewpoint on stochastic differential equations and their associated dynamics. The key point will be the fundamental fact that Brownian motion has a natural lift to a Brownian $p$-rough path, for any $2<p<3$. Once this object will be constructed by probabilistic means, the deterministic machinery for solving rough differential equations, described in the third part of the course, will enable us to associate to any realization of the Brownian rough path a solution to the rough differential equation \eqref{EqRDE}. This solution coincides almost-surely with the solution to the Stratonovich differential equation \eqref{EqStratoSDE}! One shows in that way that this solution is a continuous function of the Brownian rough path, in striking contrast with the fact that it is only a measurable function of the Brownian path itself, with no hope for a more regular dependence in a generic setting. This fact will provide a natural and easy road to the deep results of Wong-Zakai, Stroock \& Varadhan or Freidlin \& Wentzell.

\bigskip

Several other approaches to rough differential equations are available, each with their own pros and cons. We refer the reader to the books \cite{LyonsStFlour} and \cite{LyonsQianBook} for an account of Lyons' original approach; she/he is refered to the book \cite{FVBook} for a thourough account of the Friz-Victoir approach, and to the lecture note \cite{BaudoinLectureNotes} by Baudoin for an easier account of their main ideas and results, and to the forthcoming excellent lecture notes \cite{FH14} by Friz and Hairer on Gubinelli's point of view. The present approach building on \cite{BailleulFlows} does not overlap with the above ones.\footnote{Comments on these lecture notes are most welcome. Please email them at the address \textsf{ismael.bailleul@univ-rennes1.fr}}

\vfill
\pagebreak

\chapter{Flows and approximate flows}
\label{SectionFlows}

\ssk

\todo[inline, backgroundcolor=white, bordercolor=black]{Guide for this chapter}

\medskip

This first part of the course will present the backbone of our approach to rough dynamics under the form of a simple recipe for constructing flows of maps on some Banach space. Although naive, it happens to be robust enough to provide a unified treatment of ordinary, rough and stochastic differential equations. We fix throughout a Banach space V.

\medskip

The main technical difficulty is to deal with the non-commutative character of the space of maps from V to itself, endowed with the composition operation. To understand the part of the problem that does not come from non-commutativity, let us consider the following model problem. Suppose we are given a family $\mu=\big(\mu_{ts}\big)_{0\leq s\leq t\leq 1}$ of elements of some Banach space depending continuously on $s$ and $t$, and such that $\big|\mu_{ts}\big| = o_{t-s}(1)$. Is it possible to construct from $\mu$ a family $\varphi=\big(\varphi_{ts}\big)_{0\leq s\leq t\leq 1}$ of elements of that Banach space, depending continuously on $s$ and $t$, and such that we have 
\begin{equation}
\label{EqAdditiveFunctional}
\varphi_{tu}+\varphi_{us} = \varphi_{ts}
\end{equation}
for all $0\leq s\leq u\leq t\leq 1$? This additivity property plays the role of the flow property. Would the time interval $[0,1]$ be a finite discrete set $t_1<\cdots<t_n$, the additivity property \eqref{EqAdditiveFunctional} would mean that $\varphi_{ts}$ is the sum of the $\varphi_{t_{i+1}t_i}$, whose definition should be $\mu_{t_{i+1}t_i}$, as these are the only quantities we are given if no arbitrary choice is to be done. Of course, this will not turn $\varphi$ into an additive map, in the sense that property \eqref{EqAdditiveFunctional} holds true, in this discrete setting, but it suggest the following attempt in the continuous setting of the time interval $[0,1]$.
 
\ssk 
 
 Given a partition $\pi = \big\{0<t_1<\cdots<1\big\}$ of $[0,1]$ and $0\leq s\leq t\leq 1$, set 
 $$
 \varphi^\pi_{ts} = \sum_{s\leq t_i<t_{i+1}\leq t} \mu_{t_{i+1}t_i}. 
 $$
This map almost satisfies relation \eqref{EqAdditiveFunctional} as we have
$$
\varphi^\pi_{tu} + \varphi^\pi_{us} =  \varphi^\pi_{ts} -\mu_{u^+u^-} =  \varphi^\pi_{ts} + o_{|\pi|}(1),
$$
for all $0\leq s\leq u\leq t\leq 1$, where $u^-,u^+$ are the elements of $\pi$ such that $u^-\leq u<u^+$, and $|\pi| = \max\,\{t_{i+1}-t_i\}$ stands for the mesh of the partition. So we expect to find a solution $\varphi$ to our problem under the form $\varphi^\pi$, for a partition of $[0,1]$ of infinitesimal mesh, that is as a limit of $\varphi^\pi$'s, say along a sequence of refined partitions $\pi_n$ where $\pi_{n+1}$ has only one more point than $\pi_n$, say $u_n$. However, the sequence $\varphi^{\pi_n}$  has no reason to converge without assuming further conditions on $\mu$. To fix further the setting, let us consider partitions $\pi_n$ of $[0,1]$ by dyadic times, where we exhaust first all the dyadic times multiples of $2^{-k}$, in any order, before taking in the partition points multiples of $2^{-(k+1)}$. Two dyadic times $s$ and $t$ being given, both multiples of $2^{-k_0}$, take $n$ big enough for them to be points of $\pi_n$. Then, denoting by $u_n^-, u_n^+$ the two points of $\pi_n$ such that $u_n^-<u_n<u_n^+$, the quantity $\varphi^{\pi_{n+1}}_{ts}-\varphi^{\pi_n}_{ts}$ will either be null if $u_n\notin [s,t]$, or 
\begin{equation}
\label{EqDiffVarphi}
\varphi^{\pi_{n+1}}_{ts}-\varphi^{\pi_n}_{ts} = \big(\mu_{u_n^+u_n}+\mu_{u_n u_n^-}\big) - \mu_{u_n^+u_n^-},
\end{equation}
otherwise. A way to control this quantity is to assume that the map $\mu$ is \textit{approximately additive}, in the sense that  we have some positive constants $c_0$ and $a>1$ such that the inequality
\begin{equation}
\label{EqDefnAlmostAdditivity}
\big|\big(\mu_{tu}+\mu_{us}\big)-\mu_{ts}\big| \leq c_0\,|t-s|^a
\end{equation}
holds for all $0\leq s\leq u\leq t\leq 1$. Under this condition, we have
$$
\Big|\varphi^{\pi_{n+1}}_{ts}-\varphi^{\pi_n}_{ts}\Big| \leq c_0 2^{-am},
$$
where $\big|\pi_{n+1}\big| = 2^{-m}$. There will be $2^m$ such terms in the formal series $\sum_ {n\geq 0} \big( \varphi^{\pi_{n+1}}_{ts}-\varphi^{\pi_n}_{ts}\big)$, giving a total contribution for these terms of size $2^{-(a-1)m}$, summable in $m$. So this sum converges to some quantity $\varphi_{ts}$ which satisfies \eqref{EqAdditiveFunctional} by construction (on dyadic times only, as defined as above). Note that commutativity of the addition operation was used implicitly to write down equation \eqref{EqDiffVarphi}.

\medskip

Somewhat surprisingly, the above approach also works in the non-commutative setting of maps from V to itself under a condition which essentially amounts to replacing the addition operation and the norm $|\cdot|$ in condition \eqref{EqDefnAlmostAdditivity} by the composition operation and the $\mcC^1$ norm. This will be the essential content of theorem \ref{ThmConstructingFlows} below, taken from the work \cite{BailleulFlows}.

\section{$C^1$-approximate flows and their associated flows}
\label{SubsectionApproximateFlows}

We start by defining what will play the role of an approximate flow, in the same way as $\mu$ above was understood as an approximately additive map under condition \eqref{EqDefnAlmostAdditivity}. \index{$\mcC^1$-approximate flow}

\begin{defn}
A \textbf{$\mcC^1$-approximate flow on \emph{V}} is a family $\mu = \big(\mu_{ts}\big)_{0\leq s\leq t\leq T}$ of $\mcC^2$ maps from V to itself, depending continuously on $s,t$ in the topology of uniform convergence, such that 
\begin{equation}
\label{EqContinuityMu}
\big\|\mu_{ts}-\textrm{\emph{Id}}\big\|_{\mcC^2} = o_{t-s}(1)
\end{equation}
and there exists some positive constants $c_1$ and $a>1$, such that the inequality
\begin{equation}
\label{EqMuMu}
\big\|\mu_{tu}\circ\mu_{us}-\mu_{ts}\big\|_{\mcC^1} \leq c_1 |t-s|^a
\end{equation}
holds for all $0\leq s\leq u\leq t\leq T$.
\end{defn}

Note that $\mu_{ts}$ is required to be $\mcC^2$ close to the identity while we ask it to be an approximate flow in a $\mcC^1$ sense. Given a partition $\pi_{ts} = \{s=s_0<s_1<\cdots<s_{n-1}<s_n=t\}$ of an interval $[s,t]\subset [0,T]$, set
$$
\mu_{\pi_{ts}} = \mu_{t_nt_{n-1}}\circ\cdots\circ\mu_{t_1t_0}.
$$

\begin{thm}[Constructing flows on a Banach space]
\label{ThmConstructingFlows}
A $\mcC^1$-approximate flow defines a unique flow $\varphi=\big(\varphi_{ts}\big)_{0\leq s\leq t\leq T}$ on \emph{V} such that the inequality
\begin{equation}
\label{EqMuApproxVarphi}
\big\|\varphi_{ts}-\mu_{ts}\big\|_\infty \leq c |t-s|^a
\end{equation}
holds for some positive constant $c$, for all $0\leq s\leq t\leq T$ sufficiently close, say $t-s\leq \delta$. This flow satisfies the inequality
\begin{equation}
\label{EqApproxVarphiMu}
\big\|\varphi_{ts}-\mu_{\pi_{ts}}\big\|_\infty \leq \frac{2}{1-2^{1-a}}\,c_1^2T\big|\pi_{ts}\big|^{a-1}
\end{equation}
for any partition $\pi_{ts}$ of any interval $(s,t)$ of mesh $\big|\pi_{ts}\big|\leq \delta$.
\end{thm}

Note that the conclusion of theorem \ref{ThmConstructingFlows} holds in $\mcC^0$-norm. This loss of regularity \wrt the controls on $\mu$ given by equations \eqref{EqContinuityMu} and \eqref{EqMuMu} roughly comes from the use of uniform $\mcC^1$-estimates on some functions $f_{ts}$ to control some increments of the form $f_{ts}\circ g_{ts} - f_{ts}\circ g'_{ts}$, for some $\mcC^0$-close maps $g_{ts},g'_{ts}$. Note that if $\mu$ depends continuously on some parameter, then $\varphi$ also depends continuously on that parameter, as a uniform limit of continuous functions, equation \eqref{EqApproxVarphiMu}. 

The remainder of this section will be dedicated to the proof of theorem \ref{ThmConstructingFlows}. We shall proceed in two steps, by proving first that one can construct $\varphi$ as the \textit{uniform limit} of the $\mu_\pi$'s provided one can control uniformly their Lipschitz norm. This control will be proved in a second step.

\subsection{First step}
\label{SubsubSection1stStep}

Let us introduce the following inductive definition to prepare the first step. \index{Partition $\ep$-special}

\begin{defn}
Let $\ep\in (0,1)$ be given. A partition $\pi=\{s=s_0<s_1<\cdots<s_{n-1}<s_n=t\}$ of $(s,t)$ is said to be \emph{\textbf{$\ep$-special}} if it is either trivial or 
\begin{itemize}
   \item one can find an $s_i\in\pi$ sucht that $\ep\leq\frac{s_i-s}{t-s}\leq 1-\ep$,
   \item and for any choice $u$ of such an $s_i$, the partitions of $[s,u]$ and $[u,t]$ induced by $\pi$ are both $\ep$-special.
\end{itemize}
\end{defn}

A partition of any interval into sub-intervals of equal length has special type $\frac{1}{2}$. Given a partition $\pi=\{s=s_0<s_1<\cdots<s_{n-1}<s_n=t\}$ of $(s,t)$ of special type $\ep$ and $u\in\{s_1,\dots,s_{n-1}\}$ with $\ep\leq\frac{u-s}{t-s}\leq 1-\ep$, the induced partitions of the intervals $[s,u]$ and $[u,t]$ are also $\ep$-special. Set $m_\ep = \underset{\ep\leq \beta\leq 1-\ep}{\sup}\,\beta^a+(1-\beta)^a<1$, and pick a constant 
$$
L > \frac{2c_1}{1-m_\ep},
$$ 
where $c_1$ is the constant that appears in the dfinition of a $\mcC^1$-approximate flow, in equation \eqref{EqMuMu}.

\begin{lem}
\label{LemImprovedEqMumu}
Let $\mu = \big(\mu_{ts}\big)_{0\leq s\leq t\leq T}$ be a $\mcC^1$-approximate flow on \emph{V}. Given $\ep>0$, there exists a positive constant $\delta$ such that for any $0\leq s\leq t\leq T$ with $t-s\leq \delta$, and any special partition of type $\ep$ of an interval $(s,t)\subset [0,T]$, we have 
\begin{equation}
\label{EqImprovedEqMumu}
\big\|\mu_{\pi_{ts}}-\mu_{ts}\big\|_\infty \leq L|t-s|^a.
\end{equation}
 \end{lem}

\ssk

\begin{Dem}
We proceed by induction on the number $n$ of sub-intervals of the partition. The case $n=2$ is the $\mcC^ 0$ version of identity \eqref{EqMuMu}. Suppose the statement has been proved for $n\geq 2$. Fix $0\leq s<t\leq T$ with $t-s\leq \delta$, and let $\pi_{ts}=\{s_0=s<s_1<\cdots<s_n<s_{n+1}=t\}$ be an $\ep$-special partition of $[s,t]$, splitting the interval $[s,t]$ into $(n+1)$ sub-intervals. Let $u$ be one of the points of the partition sucht that $\ep \leq \frac{t-u}{t-s} \leq 1-\ep$, so the two partitions $\pi_{tu}$ and $\pi_{us}$ are both $-\ep$-special, with respective cardinals no greater than $n$. Then
\begin{equation*}
\label{FirstEstimate}
\begin{split}
\big\|\mu_{\pi_{ts}} - \mu_{ts}\big\|_\infty &\leq \big\|\mu_{\pi_{tu}}\circ\mu_{\pi_{us}} - \mu_{tu}\circ\mu_{\pi_{us}}\big\|_\infty + \big\|\mu_{tu}\circ\mu_{\pi_{us}} - \mu_{ts}\big\|_\infty \\
&\leq \big\|\mu_{\pi_{tu}} - \mu_{tu}\big\|_\infty + \big\|\mu_{tu}\circ\mu_{\pi_{us}} - \mu_{tu}\circ\mu_{us}\big\|_\infty + \big\|\mu_{tu}\circ\mu_{us} - \mu_{ts}\big\|_\infty \\
&\leq L|t-u|^a+ \big(1+o_\delta(1)\big)L\,|u-s|^a  + c_1|t-s|^a,
\end{split}
\end{equation*}
by the induction hypothesis and \eqref{EqContinuityMu} and \eqref{EqMuMu}. Set $u-s=\beta(t-s)$, with $\ep\leq \beta\leq 1-\ep$.  The above inequality rewrites 
\begin{equation*}
\begin{split}
\big\|\mu_{\pi_{ts}}-\mu_{ts}\big\|_\infty &\leq \Big\{\big(1+o_{\delta}(1)\big)\big( (1-\beta)^a+ \beta^a\big) L + c_1 \Big\} \,|t-s|^a. 
\end{split}
\end{equation*}
In order to close the induction, we need to choose $\delta$ small enough for the condition
\begin{equation}
\label{ConditionsLDelta}
c_1 +\big(1+o_\delta(1)\big)\,m_\ep L  \leq L
\end{equation}
to hold; this can be done since $m_\ep<1$.
\end{Dem}

\medskip

As a shorthand, we shall write $\mu^n_{ts}$ for $\bigcirc_{i=0}^{n-1} \mu_{t_{i+1}t_i}$, where $s_i = s+\frac{i}{n}(t-s)$. The next proposition is to be understood as the core of our approach.

\begin{prop}[Step 1]
\label{PropStep1}
Let $\mu = \big(\mu_{ts}\big)_{0\leq s\leq t\leq T}$ be a $\mcC^1$-approximate flow on \emph{V}. Assume the existence of a positive constant $\delta$ such that the maps $\mu^n_{ts}$, for $n\geq 2$ and $t-s\leq \delta$, are all Lipschitz continuous, with a Lipschitz constant uniformly bounded above by some constant $c_2$, then there exists a unique flow $\varphi=\big(\varphi_{ts}\big)_{0\leq s\leq t\leq T}$ on \emph{V} such that the inequality
\begin{equation}
\label{Eq}
\big\|\varphi_{ts}-\mu_{ts}\big\|_\infty \leq c |t-s|^a
\end{equation}
holds for some positive constant $c$, for all $0\leq s\leq t\leq T$ with $t-s\leq \delta$. This flow satisfies the inequality
\begin{equation}
\label{EqApproxVarphiMu}
\big\|\varphi_{ts}-\mu_{\pi_{ts}}\big\|_\infty \leq c_1c_2T\big|\pi_{ts}\big|^{a-1}
\end{equation}
for \emph{any} partition $\pi_{ts}$ of $(s,t)$, of mesh $\big|\pi_{ts}\big|\leq \delta$.
\end{prop}

\begin{Dem}
The existence and uniqueness proofs both rely on the elementary identity
\begin{equation}
\label{ElementaryIdentity}
f_N\circ\cdots\circ f_1 \,- \,g_N\circ\cdots\circ g_1 = \sum_{i=1}^N \Big(g_N\circ\cdots \circ g_{N-i+1}\circ f_{N-i}\, -\, g_N\circ\cdots \circ \,g_{N-i+1}\circ\, g_{N-i}\Big)\circ f_{N-i-1} \circ\cdots\circ f_1,
\end{equation}
where the $g_i$ and $f_i$ are maps from V to itself, and where we use the obvious convention concerning the summand for the first and last term of the sum. In particular, if all the maps $g_N\circ\cdots\circ g_k$ are Lipschitz continuous, with a common upper bound $c'$ for their Lipschitz constants, then 
\begin{equation}
\label{CsqElementaryIdentity}
\big\|f_N\circ\cdots\circ f_1 - g_N\circ\cdots\circ g_1\big\|_\infty \leq c'\sum_{i=1}^N \|f_i-g_i\|_\infty.
\end{equation}

\ssk

\textbf{a) Existence.} Set $\textrm{D}_\delta := \big\{0\leq s\leq t\leq T\,;\,t-s\leq \delta\big\}$ and write $\DD_\delta$ for the intersection of $\textrm{D}_\delta $ with the set of dyadic real numbers. Given $s=a2^{-k_0}$ and $t=b2^{-k_0}$ in $\DD_\delta$, define for $n\geq k_0$
$$
\mu^{(n)}_{ts} := \mu_{ts}^{2^n} = \mu_{s_{N(n)}s_{N(n)-1}}\circ\cdots\circ\mu_{s_1s_0} ,
$$
where $s_i=s+i2^{-n}$ and $s_{N(n)}=t$. Given $n\geq k_0$, write 
$$
\mu^{(n+1)}_{ts} = \overset{N(n)-1}{\underset{i=0}{\bigcirc}}\big(\mu_{s_{i+1}s_i+2^{-n-1}}\circ\mu_{s_i+2^{-n-1}s_i}\big)
$$
and use \eqref{ElementaryIdentity} with $f_i = \mu_{s_{i+1}s_i+2^{-n-1}}\circ\mu_{s_i+2^{-n-1}s_i}$ and $g_i=\mu_{s_{i+1}s_i}$ and the fact that all the maps $\mu_{s_{N(n)}s_{N(n)-1}} \circ \cdots \circ \mu_{s_{N(n)-i+1} s_{N(n)-i}} = \mu_{s_{N(n)}s_{N(n)-i}}^i$ are Lipschitz continuous with a common Lipschitz constant $c_2$, by assumption, to get by \eqref{CsqElementaryIdentity} and \eqref{EqMuMu}
$$
\Big\|\mu^{(n+1)}_{ts} - \mu^{(n)}_{ts}\Big\|_\infty \leq c_2\sum_{i=0}^{N(n)-1}\big\|\mu_{s_{i+1}s_i+2^{-n-1}}\circ\mu_{s_i+2^{-n-1}s_i} - \mu_{s_{i+1}s_i}\big\|_\infty \leq  c_1c_2 T\,2^{-(a-1)n};
$$
so $\mu^{(n)}$ converges uniformly on $\DD_\delta$ to some continuous function $\varphi$. We see that $\varphi$ satisfies inequality \eqref{EqMuApproxVarphi} on $\DD_\delta$ as a consequence of \eqref{EqImprovedEqMumu}. As $\varphi$ is a uniformly continuous function of $(s,t)\in\DD_\delta$, by \eqref{EqMuApproxVarphi}, it has a unique continuous extension to $\textrm{D}_\delta$, still denoted by $\varphi$. To see that it defines a flow on $\textrm{D}_\delta$, notice that for dyadic times $s\leq u\leq t$, we have $\mu^{(n)}_{ts} = \mu^{(n)}_{tu}\circ\mu^{(n)}_{us}$, for $n$ big enough; so, since the maps $\varphi^{(n)}_{tu}$ are uniformly Lipschitz continuous, we have $\varphi_{ts}=\varphi_{tu}\circ\varphi_{us}$ for such triples of times in $\DD_\delta$, hence for all times since $\varphi$ is continuous. The map $\varphi$ is easily extended as a flow to the whole of $\{0\leq s\leq t\leq T\}$. Note that $\varphi$ inherits from the $\mu^n$'s their Lipschitz character, for a Lipschitz constant bounded above by $c_2$.

\medskip

\noindent \textbf{b) Uniqueness.} Let $\psi$ be any flow satisfying condition \eqref{EqMuApproxVarphi}. With formulas \eqref{ElementaryIdentity} and \eqref{CsqElementaryIdentity} in mind, rewrite \eqref{EqMuApproxVarphi} under the form $\psi_{ts}=\mu_{ts}+O_c\bigl(|t-s|^a\bigr)$, with obvious notations. Then
\begin{equation*}
\begin{split}
\psi_{ts} &= \psi_{s_{2^n}s_{2^n-1}}\circ\cdots\circ\psi_{s_1s_0} = \Bigl(\mu_{s_{2^n}s_{2^n-1}} + O_c\bigl(2^{-an}\bigr)\Bigr) \circ\cdots\circ \Bigl(\mu_{s_1s_0} + O_c\bigl(2^{-an}\bigr)\Bigr)  \\
             &= \mu_{s_{2^n}s_{2^n-1}}\circ\cdots\circ\mu_{s_1s_0} + \Delta_n = \mu^{(n)}_{ts} + \Delta_n,
\end{split}
\end{equation*}
where $\Delta_n$ is of the form of the right hand side of \eqref{ElementaryIdentity}, so is bounded above by a constant multiple of $2^{-(a-1)n}$, since all the  maps $\mu_{s_{2^n}s_{2^n-1}} \circ \cdots \circ \mu_{s_{2^n-\ell+1} s_{2^n-\ell}}$ are Lipschitz continuous with a common upper bound for their Lipschitz constants, by assumption.  Sending $n$ to infinity shows that $\psi_{ts}=\varphi_{ts}$.

\medskip

\noindent {\bf c) Speed of convergence.} Given any partition $\pi = \{s_0=s<\cdots<s_n=t\}$ of $(s,t)$, writing $\varphi_{ts} = \bigcirc_{i=0}^{n-1} \varphi_{s_{i+1}s_i}$, and using their uniformly Lipschitz character, we see as a consequence of \eqref{CsqElementaryIdentity} that we have for $\big|\pi_{ts}\big|\leq \delta$
\begin{equation*}
\big\|\varphi_{ts}-\mu_{\pi_{ts}}\big\|_\infty \leq c_2\sum_{i=0}^{n-1} \big\|\varphi_{s_{i+1}s_i}-\mu_{s_{i+1}s_i}\big\|_\infty \leq c_1c_2\sum_{i=0}^{n-1} |s_{i+1}-s_i|^a \leq c_1c_2T\,\big|\pi_{ts}\big|^{a-1}.
\end{equation*}
\end{Dem}

\medskip

Compare what is done in the above proof with what was done in the introduction to this part of the course in a commutative setting.

\medskip

\subsection{Second step}
\label{SubsubsectionSecondStep}

The uniform Lipschitz control assumed in proposition \ref{PropStep1} actually holds under the assumption that $\mu$ is a $\mcC^1$-approximate flow. The results of this paragraph could have been proved just after lemma \ref{LemImprovedEqMumu} and do not use the result proved in the fundamental proposition \ref{PropStep1}. Recall $L$ stands for a constant strictly greater than $\frac{2c_1}{1-m_\ep}$.

\begin{prop}[Uniform Lipschitz controls]
\label{PropSecondStep}
Let $\mu = \big(\mu_{ts}\big)_{0\leq s\leq t\leq T}$ be a $\mcC^1$-approximate flow on \emph{V}. Then, given $\ep>0$, there exists a positive constant $\delta$ such that the inequality
$$
\big\|\mu_{\pi_{ts}} - \mu_{ts}\big\|_{\mcC^1} \leq L|t-s|^ a
$$
holds for any partition $\pi_{ts}$ of $[s,t]$ of special type $\ep$, whenever $t-s\leq \delta$.
\end{prop}

\begin{Dem}
We proceed by induction on the number $n$ of sub-intervals of the partition as in the proof of lemma \ref{LemImprovedEqMumu}. The case $n=2$ is identity \eqref{EqMuMu}. Suppose the statement has been proved for $n\geq 2$. Fix $0\leq s<t\leq T$ with $t-s\leq \delta$, and let $\pi_{ts}=\{s_0=s<s_1<\cdots<s_n<s_{n+1}=t\}$ be an $\ep$-special partition of $[s,t]$ of special, splitting the interval $[s,t]$ into $(n+1)$ sub-intervals. Let $u$ be a point of the partition with $\ep\leq \frac{u-s}{t-s}\leq 1-\ep$, so that the two partitions $\pi_{tu}$ and $\pi_{us}$ are both $\ep$-special, with respective cardinals no greater than $n$. Then, for any $x\in \textrm{V}$, one can write $D_x\mu_{\pi_{ts}}-D_x\mu_{ts}$ as a telescopic sum which involve only some controlled quantities.
\begin{equation*}
\label{FirstEstimate}
\begin{split}
D_x &\mu_{\pi_{ts}} - D_x\mu_{ts} = D_x\big(\mu_{\pi_{tu}}\circ\mu_{\pi_{us}}\big) - D_x\mu_{ts} \\
        													   &= \Big(D_{\mu_{\pi_{us}}(x)}\mu_{\pi_{tu}} - D_{\mu_{\pi_{us}}(x)}\mu_{tu} \Big)\big(D_x\mu_{\pi_{us}}\big) + \Big( \big(D_{\mu_{\pi_{us}}(x)}\mu_{tu} - D_{\mu_{us}(x)}\mu_{tu}\big)\big(D_x\mu_{\pi_{us}}\big)\Big) \\
        													   &\quad +\big(D_{\mu_{us}(x)}\mu_{tu}\big)\Big(D_x\mu_{\pi_{us}}-D_x\mu_{us}  \Big) + \Big(\big(D_{\mu_{us}(x)}\mu_{tu}\big)\big(D_x\mu_{us}\big) - D_x\mu_{ts}\Big) \\
        													   &=: (1) + (2) + (3) + (4)
\end{split}
\end{equation*}
We treat each term separately using repeatedly the induction hypothesis, continuity assumption \eqref{EqContinuityMu} for $\mu_{ts}$ in $\mcC^2$ topology, and lemma \ref{LemImprovedEqMumu} when needed. We first have 
$$
\big|(1)\big| \leq L|t-u|^a\,\big(1+o_\delta(1)\big).
$$ 
Also, 
$$
\Big|D_{\mu_{\pi_{us}}(x)}\mu_{tu} - D_{\mu_{us}(x)}\mu_{tu}\Big| \leq o_{t-u}(1)\,\big|\mu_{\pi_{us}}(x) - \mu_{us}(x)\big|\leq o_{t-u}(1)\,L|u-s|^a,
$$
As the term $D_x\mu_{\pi_{us}}$ has size no greater than $\big(1+o_\delta(1)\big) + L|u-s|^a$, we have
$$
\big|(2)\big| \leq o_\delta(1)\,|u-s|^a.
$$
Last, we have the upper bound
$$
\big|(3)\big| \leq \big(1+o_\delta(1)\big) L|u-s|^a,
$$
while  $\big|(4)\big| \leq \big\|\mu_{tu}\circ\mu_{us}-\mu_{ts}\big\|_{\mcC^1} \leq c_1|t-s|^a$ by \eqref{EqMuMu}. All together, and writing $t-u=\beta (t-s)$, for some $\beta\in [\ep,1-\ep]$, this gives 
\begin{equation*}
\begin{split}
\big|D_x \mu_{\pi_{ts}} - D_x\mu_{ts}\big| &\leq \Big(\big(1+o_\delta(1)\big)\big(\beta^a + (1-\beta)^a\big)L + c_1+o_\delta(1)\Big)|t-s|^a \\
&\leq L\,|t-s|^a
\end{split}
\end{equation*}
for $\delta$ small enough, as $m_\ep<1$.
\end{Dem}

\medskip

Propositions \ref{PropStep1} and \ref{PropSecondStep} together prove theorem \ref{ThmConstructingFlows}. Note that an explicit choice of $\delta$ is possible as soon as one has a quantitative version of the estimate $\big\|\mu_{ts}-\textrm{Id}\big\|_{\mcC^2} = o_{t-s}(1)$. Note also that proposition \ref{PropSecondStep} provides an explicit control on the Lipschitz norm of the $\varphi_{ts}$, in terms of the Lipschitz norm of $\mu_{ts}$ and  $L$.

\section{Exercices on flows}
\label{SubsectionExercicesFlows}

{\small To get a hand on the machinery of $\mcC^1$-approximate flows, we shall first see how theorem \ref{ThmConstructingFlows} gives back the classical Cauchy-Lipschitz theory of ordinary differential equations for \textit{bounded} Lipschitz vector fields on $\RR^d$. Working with \textit{unbounded} Lipschitz vector fields requires a slightly different notion of \textit{local} $\mcC^1$-approximate flow -- see \cite{BailleulFlows}. 

\ssk

Theorem \ref{ThmConstructingFlows} can be understood as a non-commutative analogue of Feyel-de la Pradelle's sewing lemma \cite{FdlP}, first introduced by Gubinelli \cite{Gubinelli04} as an abstraction of a fundamental mechanism invented by Lyons \cite{Lyons97}. Exercices 2-4 are variations on this commutative version of theorem \ref{ThmConstructingFlows}, as already sketched in the introduction to this part.

\bigskip

{\bf 1. Ordinary differential equations.} Let $V_1,\dots,V_\ell$ be $\mcC^2_b$ vector fields on $\RR^d$ (or a Banach space), and $h_1,\dots,h_\ell$ be real-valued $\mcC^1$ controls. Let $\varphi$ stand for the flow associated with the ordinary differential equation
$$
dx_t = V_i(x_t)dh^i_t. 
$$

{\bf $\quad$ a)} Show that one defines a $\mcC^1$-approximate flow setting for all $x\in\RR^d$
$$
\mu_{ts}(x) = x+\big(h_t-h_s\big)^iV_i(x).
$$

{\bf $\quad$ b)} Prove that $\varphi$ is equal to the flow associated to $\mu$ by theorem \eqref{ThmConstructingFlows}. In that sense, a path $x$ is a solution to the above ordinary differential equation if and only if it satisfies at any time $s$ the Taylor-type expansion formula
$$
f\big(x_t\big) = f\big(x_s\big) + \big(h^i_t-h^i_s\big)\big(V_if\big)(x_s) + o(t-s),
$$
for any function $f$ of class $\mcC^2_b$. Show that the above reasoning holds true if we only assume that the $\RR^\ell$-valued control $h$ is globally Lipschitz continuous. (It is actually sufficient to suppose $h$ is $\alpha$-H\"older, for some $\alpha>\frac{1}{2}$.)\vspace{0.1cm}

{\bf $\quad$ c)} Does anything go wrong with the above reasoning if the Lipschitz continuous vector fields $V_i$ are not bounded?  \vspace{0.1cm}

{\bf $\quad$ d)} Show that $\varphi$ depends continuously on $h$ in the uniform topology for $\varphi$ and the Lipschitz topology for $h$, defined by the distance
$$
d(h,h') = \big|h_0-h'_0\big| + \textrm{Lip}(h-h'),
$$
where $\textrm{Lip}(h-h')$ stands for the Lipschitz norm of $h-h'$. (A similar result holds if $h$ is $\alpha$-H\"older, for some $\alpha>\frac{1}{2}$, with the Lipschitz norm replaced by the $\alpha$-H\"older norm.)

\bigskip

{\bf 2. Feyel-de la Pradelle' commutative sewing lemma.} \index{Sewing lemma} Let V be a Banach space and $\mu=\big(\mu_{ts}\big)_{0\leq s\leq t\leq 1}$ be a V-valued continuous function. The following commutative version of theorem \ref{ThmConstructingFlows} was first proved under this form by Feyel and de la Pradelle in \cite{FdlP}; see also \cite{FdlPM}. Suppose there exists some positive constants $c_0$ and $a>1$ such that we have
\begin{equation}
\label{EqAlmostIncrement}
\big|(\mu_{tu}+\mu_{us})-\mu_{ts}\big| \leq c_0|t-s|^a
\end{equation}
for all $0\leq s\leq u\leq t\leq 1$. We say the $\mu$ is an \textit{almost-additive functional} (or map). 

Simplify the proof of theorem \ref{ThmConstructingFlows} to show that there exists a unique map $\varphi = \big(\varphi_t\big)_{0\leq t\leq 1}$, whose increments $\varphi_{ts} := \varphi_t-\varphi_s$, satisfy
$$
\big|\varphi_{ts}-\mu_{ts}\big| \leq c|t-s|^a
$$
for some positive constant $c$ and all $0\leq s\leq t\leq 1$. \vspace{0.3cm}

{\bf 3. Integral products.} Let $\alpha>\frac{1}{2}$ be given, and $\big(A_t\big)_{0\leq t\leq 1}$ be an $\alpha$-H\"older path with values in the space $\textrm{L}_c(\textrm{V})$ of continuous linear maps from V to itself. Set $A_{ts} = A_t-A_s$, for $0\leq s\leq t\leq 1$, and use the notation $|\cdot|$ for the operator norm on $\textrm{L}_c(\textrm{V})$.\vspace{0.1cm}

{\bf $\quad$ a)} Use theorem \eqref{ThmConstructingFlows} to show that setting $\mu_{ts} = \textrm{Id}+A_{ts}$ determines a unique flow $\varphi$ on V. A good notation for $\varphi_{ts}$ is $\prod_{s\leq r\leq t} \big(\textrm{Id}+dA_r\big) = \prod_{s\leq r\leq t} e^{dA_r}$. \vspace{0.1cm}

{\bf $\quad$ b)} Let $\big(B_t\big)_{0\leq t\leq 1}$ be another $\textrm{L}_c(\textrm{V})$-valued $\alpha$-H\"older path. Show that we define a $\mcC^1$-approximate flow $\mu$ setting $\mu_{ts}=\big(\textrm{Id}+A_{ts}\big)\big(\textrm{Id}+B_{ts}\big)$. Prove as a consequence the well-known formula 
$$
e^{A_0+B_0} = \lim_n\,\Big(e^{\frac{A_0}{2^n}}e^{\frac{B_0}{2^n}}\Big)^{2^n}.
$$

\bigskip

{\bf 4. Controls and finite-variation paths.} \index{Finite variation paths} A \textit{control} is a non-negative map $\omega = \big(\omega_{ts}\big)_{0\leq s\leq t\leq 1}$, null on the diagonal, and such that we have
$$
\omega_{tu} + \omega_{us} \leq \omega_{ts}
$$
for all $0\leq s\leq u\leq t\leq 1$. \vspace{0.1cm}

{\bf $\quad$ a)} Show that Feyel and de la Pradelle' sewing lemma holds true if we replace $t-s$ by $\omega_{ts}$, if we suppose that $\omega^a$ is a control. \vspace{0.1cm}

{\bf $\quad$ b)} Recall that a V-valued path $x= \big(x_t\big)_{0\leq t\leq 1}$ is said to have finite $p$-variation, for some $p\geq 1$, if the following quantity is finite for all $0\leq s\leq t\leq 1$:
$$
|x|^p_{p-\textrm{var}; [s,t]} := \sup\,\sum \Big|x_{t_{i+1}}-x_{t_i}\Big|^p,
$$
with a sum over the partition points $t_i$ of a given partition of the interval $[s,t]$, and a supremum over the set of all partitions of $[s,t]$. Such a definition is invariant by any reparametrization of the time interval $[s,t]$. Given such a path, show that setting $\omega_{ts} = |x|^p_{p-\textrm{var}; [s,t]}$, defines a control $\omega$. \vspace{0.1cm}

{\bf $\quad$ c)} Show that a path with finite $p$-variation can be reparametrized into a $\frac{1}{p}$-H\"older path, with $1$-H\"older paths being understood as Lipschitz continuous paths. Given an $\RR^\ell$-valued path $h$ with finite $1$-variation, set 
$$
\zeta_s = \inf\{t\geq 0\,;\,|h|_{1-\textrm{var};[0,t]} \geq s\}.
$$
We define a solution $x_\bullet$ to the ordinary differential equation 
$$
dx_t = V_i(x_t)dh^i_t
$$
driven by $h$ as a path $x_\bullet$ such that the reparametrized path $y:=x\circ\zeta$ is a solution to the ordinary differential equation 
$$
dy_s = V_i(y_s)d(h\circ\zeta)^i_t,
$$
driven by the globally Lipschitz path $h\circ\zeta$. \vspace{0.1cm}

{\bf $\quad$ d)} Prove that the flow $\varphi$ constructed in this case from theorem \ref{ThmConstructingFlows} depends continuously on $h$ in the uniform norm for $\varphi$ and the $1$-variation topology associated with the norm $|\cdot|_{1-\textrm{var}}$ for $h$. (Following the remarks of exercice 1, one can actually prove the results of questions c) and d) for paths with finite $p$-variation, for $1\leq p<2$.)

\vspace{0.3cm}

{\bf 5. Young integral.} \index{Young integral} Given another Banach space E, denote by $\textrm{L}_c(\textrm{V,E})$ the space of continuous linear maps from V to E equipped with the operator norm. Let $\alpha$ and $\beta$ be positive real numbers such that $\alpha+\beta>1$. Given any $0<\alpha<1$, we denote by $\textrm{Lip}_\alpha(\textrm{E})$ the set of $\alpha$-H\"older maps. This unusual notation will be justified in the third path of the course.

\vspace{0.1cm}

{\bf $\quad$ a)} Given an $\textrm{L}_c(\textrm{V,E})$-valued $\alpha$-Lipschitz map $x= \big(x_{s}\big)_{0\leq s\leq 1}$ and a V-valued $\beta$-Lipschitz map $y= \big(y_{s}\big)_{0\leq s\leq 1}$, show that setting
$$
\mu_{ts} = x_s\big(y_t-y_s\big)
$$
for all $0\leq s\leq t\leq 1$, defines an E-valued function $\mu$ that satisfies equation \eqref{EqAlmostIncrement}, with a constant $c_0$ to be made explicit. \vspace{0.1cm}

{\bf $\quad$ b)} The associated function $\varphi$ is denoted by $\varphi_t = \int_0^t x_sdy_s$, for all $0\leq t\leq 1$. Show that it is a continuous function of $x\in\textrm{Lip}_\alpha\big(\textrm{L}_c(\textrm{V,E})\big)$ and $y\in\textrm{Lip}_\beta(\textrm{V})$. \vspace{0.3cm}

\medskip

 {\bf 6. Lipschitz dependence of $\varphi$ on $\mu$.} As emphasized in the remark following theorem \ref{ThmConstructingFlows}, inequality \eqref{EqApproxVarphiMu} implies that $\varphi$, understood as a function of $\mu$, is continuous in the $\mcC^0$-norm on the sets of $\mu$'s of the form 
$\big\{\mu\,;\,\textrm{\eqref{EqMuMu} holds uniformly}\big\}$, equipped with the $\mcC^0$-norm. One can actually prove that it depends Lipschitz continuously on $\mu$ in the following sense. \vspace{0.1cm}

Let $\mu=\big(\mu_{ts}\big)_{0\leq s\leq t\leq 1}$ and $\mu'=\big(\mu'_{ts}\big)_{0\leq s\leq t\leq 1}$ be $\mcC^1$-approximate flows on V, with associated flows $\varphi$ and $\varphi'$. Suppose that we have 
$$
\big\|\mu_{ts}-\mu'_{ts}\big\|_{\mcC^1} \leq \ep\,|t-s|^\frac{1}{p}
$$
and
$$
\big\|\big(\mu_{tu}\circ\mu_{us}-\mu_{ts}\big) - \big(\mu'_{tu}\circ\mu'_{us}-\mu'_{ts}\big)\big\|_{\mcC^1} \leq \ep\,|t-s|^a
$$
for a positive constant $\ep$, with $a>1$ as in the definition of the $\mcC^1$-approximate flows $\mu,\mu'$, for all $0\leq s\leq u\leq t\leq 1$. Prove that one has
$$
\big\|\big(\varphi_{ts}-\mu_{ts}\big) - \big(\varphi'_{ts}-\mu'_{ts}\big)\big\|_\infty \leq c\,\ep\,|t-s|^a,
$$
for all for $0\leq s\leq t\leq 1$, for some explicit positive constant $c$.
}

\vfill
\pagebreak


\chapter{Rough paths}
\label{SectionRoughPaths}

\ssk

\todo[inline, backgroundcolor=white, bordercolor=black]{Guide for this chapter}

\medskip

H\"older $p$-rough paths, which control the rough differential equations 
$$
dx_t = \textrm{F}(x_t){\bfX}(dt),  \quad d\varphi_t = \textrm{F}^\otimes{\bfX}(dt),
$$
and play the role of the control $h$ in the model classical ordinary differential equation 
$$
dx_t = V_i(x_t)\,dh^i_t = \textrm{F}(x_t)\,dh_t
$$ 
are defined in section \ref{SubsubsectionDefnRoughPaths}. As $\RR^\ell$-valued paths, they are not regular enough for the formula 
$$
\mu_{ts}(x) = x+ X_{ts}^iV_i(x)
$$
to define an approximate flow, as in the classical Euler scheme studied in exercice 1. The missing bit of information needed to stabilize the situation is a substitute of the non-existing iterated integrals $\int_s^t X^j_rdX^k_r$, and higher order iterated integrals, which provide a partial description of what happens to $X$ during any time interval $(s,t)$. A (H\"older) $p$-rough path is a \textit{multi-level object} whose higher order parts provide precisely that information.  We saw in the introduction that iterated integrals appear naturally in Taylor-Euler expansions of solutions to ordinary differential equations; they provide higher order numerical schemes like Milstein' second order scheme. It is an important fact that $p$-rough paths take values in a very special kind of algebraic structure, whose basic features are explained in section \ref{SubsectionAlgebraicInterlude}. A H\"older $p$-rough path will then appear as a kind of $\frac{1}{p}$-H\"older path in that space. We shall then study in section \ref{SubsectionRoughPathSpace} the space of $p$-rough path for itself.

\section{Definition of a H\"older $p$-rough path}
\label{SubsectionDefnRoughPath}

Iterated integrals, as they appear for instance under the form $\int_s^t\int_s^y dh^j_r\,dh^k_u$ or $\int_s^t\int_s^y\int_s^r (\cdots)$, are multi-indexed quantities. A useful formalism to work with such objects is provided by the notion of tensor product. We first start our investigations by recalling some elementary facts about that notion. Eventually, all what will be used for practical computations on rough differential equations will be a product operation very similar to the product operation on polynomials. This abstract setting however greatly clarifies the meaning of these computations.

\subsection{An algebraic prelude: tensor algebra over $\RR^\ell$ and free nilpotent Lie group}
\label{SubsectionAlgebraicInterlude}

Let first recall what the algebraic tensor product\index{Tensor product} ${\textrm U}\otimes {\textrm V}$ of any two Banach spaces U and V is. Denote by V' the set of all continuous linear forms on V. Given $u\in{\textrm U}$ and $v\in{\textrm V}$, we define a continuous linear map on V' setting 
$$
(u\otimes v)(v') = (v',v)\,u,
$$
for any $v'\in{\textrm V}'$. The algebraic tensor product ${\textrm U}\otimes {\textrm V}$ is the set of all finite linear combinations of such maps. Its elementary elements $u\otimes v$ are $1$-dimensional rank maps. Note that an element of ${\textrm U}\otimes {\textrm V}$ can have several different decompositions as a sum of elementary elements; this has no consequences as they all define the same map from V' to U.

\ssk

As an example, $\RR^\ell\otimes (\RR^\ell)'$ is the set of all linear maps from $\RR^\ell$ to itself, that is $\textrm{L}(\RR^\ell)$. We keep that interpretation for $\RR^\ell\otimes \RR^\ell$, as $\RR^\ell$ and $(\RR^\ell)'$ are canonically identified. To see which element of $\textrm{L}(\RR^\ell)$ corresponds to $u\otimes v$, it suffices to look at the image of the $j^{\textrm{th}}$ vector $\ep_j$ of the canonical basis by the linear map $u\otimes v$; it gives the $j^{\textrm{th}}$ column of the matrix of $u\otimes v$ in the canonical basis. We have
$$
(u\otimes v)\big(\ep_j\big) = (v,\ep_j)\,u.
$$
The family $\big(\ep_{i_1}\otimes \cdots\otimes \ep_{i_k}\big)_{1\leq i_1,\dots, i_k\leq \ell}$ defines the canonical basis of $(\RR^\ell)^{\otimes k}$.

\medskip

$\bullet$ For $N\in\NN\cup\{\infty\}$, write $T^{(N)}_\ell$ for the direct sum $\underset{r=0}{\overset{N}{\bigoplus}}\big(\RR^\ell\big)^{\otimes r}$, with the convention that $\big(\RR^\ell\big)^{\otimes 0}$ stands for $\RR$. Denote by $\bfa =  \underset{r=0}{\overset{N}{\oplus}} a^r$ and $\bfb =  \underset{r=0}{\overset{N}{\oplus}} b^r$ two generic elements of $T^{(N)}_\ell$. The vector space $T^{(N)}_\ell$ is an algebra for the operations
\begin{equation}
\label{AlgebraOperations}
\begin{split}
&\bfa + \bfb = \underset{r=0}{\overset{N}{\oplus}} (a^r+b^r), \\
&\bfa\bfb = \underset{r=0}{\overset{N}{\oplus}} c^r, \quad\textrm{with }\; c^r=\sum_{k=0}^r a^k\otimes b^{r-k}\,\in(\RR^\ell)^{\otimes r}
\end{split}
\end{equation}
It is called the \textbf{(truncated) tensor algebra of $\RR^\ell$} (if $N$ is finite). \index{Tensor algebra} Note the similarity between these rules and the analogue rules for addition and product of polynomials.

\ssk

The exponential map $\exp : T^{(\infty)}_\ell\rightarrow T^{(\infty)}_\ell$ and the logarithm map $\log : T^{(\infty)}_\ell\rightarrow T^{(\infty)}_\ell$ are defined by the usual series \index{Exponential and logarithm}
\begin{equation}
\label{EqExpLog}
\exp(\bfa) = \sum_{n\geq 0}\frac{\bfa^n}{n !},\quad \log(\bfb) = \sum_{n\geq 1}\frac{(-1)^n}{n}(1-\bfb)^n,
\end{equation}
with the convention ${\bfa}^0=1\in \RR\subset T^{(\infty)}_\ell$. Denote by $\pi_N : T^{(\infty)}_\ell\rightarrow T^{(N)}_\ell$ the natural projection. We also denote by $\exp$ and $\log$ the restrictions to $T^{(N)}_\ell$ of the maps $\pi_N\circ\exp$ and $\pi_N\circ\log$ respectively. Denote by $T^{(N),1}_\ell$, resp. $T^{(N),0}_\ell$, the elements $a_0\oplus\cdots\oplus c_N$ of $T^{(N)}_\ell$ \st $a_0=0$, resp. $a_0=1$. All the elements of $T^{(N),1}_\ell$ are invertible, and $\exp : T^{(N),0}_\ell\rightarrow T^{(N),1}_\ell$ and $\log : T^{(N),1}_\ell\rightarrow T^{(N),0}_\ell$ are smooth reciprocal bijections. As an example, with ${\bfa} = 1\oplus a^1\oplus a^2\in T^{2,1}_\ell$, we have
$$
\log {\bfa} = 0\oplus (-a^1)\oplus \Big(\frac{1}{2}a^1\otimes a^1-a^2\Big).
$$

The set $T^{(N),1}_\ell$ is naturally equipped with a norm defined by the formula \index{Homogeneous norm on a tensor product}
\begin{equation}
\label{EqDefnHomogeneousNormAmbiantSpace}
\|{\bfa}\| := \sum_{i=1}^\ell \big\|a^i\big\|_{\textrm{Eucl}}^{\frac{1}{i}},
\end{equation}
where $\big\|a^i\big\|_{\textrm{Eucl}}$ stands for the Euclidean norm of $a^i\in(\RR^\ell)^{\otimes i}$, identified with an element of $\RR^{\ell^i}$ by  looking at its coordinates in the canonical basis of $(\RR^\ell)^{\otimes i}$. The choice of power $\frac{1}{i}$ comes from the fact that $T^{(N),1}_\ell$ is naturally equipped with a \textit{dilation} operation \index{Dilation}
\begin{equation}
\label{EqDefnDeltaLambda}
\delta_\lambda({\bfa}) = \big(1,\lambda a^1,\dots,\lambda^Na^N\big),
\end{equation}
so the norm $\|\cdot\|$ is homogeneous with respect to this dilation, in the sense that one has
$$
\big\|\delta_\lambda({\bfa})\big\| = |\lambda|\|{\bfa}\|
$$
for all $\lambda\in\RR$, and all ${\bfa}\in T^{(N),1}_\ell$.

\ssk

The formula $[\bfa,\bfb]=\bfa\bfb-\bfb\bfa$, defines a Lie bracket on $T^{(N)}_\ell$. Define inductively $\frak{f}:=\frak{f}^1 := \RR^\ell$, considered as a subset of $T^{(\infty)}_\ell$, and $\frak{f}^{n+1}=[\frak{f},\frak{f}^n]\subset T^{(\infty)}_\ell$. \index{Lie Bracket}

\begin{defn}
\begin{itemize}
   \item The Lie algebra $\frak{g}^{N}_\ell$ generated by the $\frak{f}^1,\dots,\frak{f}^N$ in $T^{(N)}_\ell$ is called the \emph{\textbf{$N$-step free nilpotent Lie algebra}}. \index{Nilpotent Lie algebra}
   \item As a consequence of Baker-Campbell-Hausdorf-Dynkin formula, the subset $\exp\big(\frak{g}^{N}_\ell\big)$ of $T^{(N),1}_\ell$ is a group for the multiplication operation. It is called the \emph{\textbf{$N$-step nilpotent Lie group}} on $\RR^\ell$ and denoted by $\frak{G}^{N}_\ell$. \index{Nilpotent Lie group}
\end{itemize}
\end{defn}

As all finite dimensional Lie groups, the $N$-step nilpotent Lie group is equipped with a natural (sub-Riemannian) distance inherited from its manifold structure. \index{Sub-Riemannian distance} Its definition rests on the fact that the element ${\bf a}u$ of $T^{(N)}_\ell$ is for any ${\bf a}\in \frak{G}^{(N)}_\ell$ and $u\in\RR^\ell\subset T^{(N)}_\ell$ a tangent vector to $\frak{G}^{(N)}_\ell$ at point ${\bf a}$ (as $u$ is tangent to $\frak{G}^{(N)}_\ell$ at the identity and tangent vectors are transported by left translation in the group). So the ordinary differential equation
$$
d{\bfa}_t = {\bfa}_t\,\dot h_t
$$
makes sense for any $\RR^\ell$-valued smooth control $h$, and defines a path in $\frak{G}^{(N)}_\ell$ started from the identity. We define the size $|{\bfa}|$ of $\bfa$ by the formula
$$
|{\bfa}| = \inf\,\int_0^1 \big|\dot h_t\big|\,dt,
$$
where the infimum is over the set of all piecewise smooth controls $h$ such that ${\bfa}_1=\bfa$. This set is non-empty as ${\bfa}\in\exp\big(\frak{g}^{N}\big)$ can be written as ${\bfa}_1$ for some piecewise $\mcC^1$ control, as a consequence of a theorem of sub-Riemannian geometry due to Chow; see for instance the textbook \cite{Montgomery} for a nice account of that theorem. The distance between any two points $\bfa$ and $\bfb$ of $\frak{G}^{(N)}_\ell$ is then defined as $\big|{\bfa}^{-1}\bfb\big|$. It is homogeneous in the sense that if ${\bfa} = \exp(u)$, with $u\in\RR^\ell\subset T^{(N)}_\ell$, then $\big|\exp(\lambda u)\big| = |\lambda||{\bfa}|$, for all $\lambda\in\RR$ and all $u\in\RR^\ell\subset T^{(N)}_\ell$. 

This way of defining a distance is intrinsic to $\frak{G}^{(N)}_\ell$ and classical in geometry. From an extrinsic point of view, one can also consider $\frak{G}^{(N)}_\ell$ as a subset of $T^{(N)}_\ell$ and use the ambiant metric to define the distance between any two points $\bfa$ and $\bfb$ of $\frak{G}^{(N)}_\ell$ as $\big\|{\bfa}^{-1}{\bfb}\big\|$. It can be proved (this is elementary, see e.g. proposition 10 in Appendix A of \cite{LejayIntro2}, pp. 76-77) that the two norms $|\cdot|$ and $\|\cdot\|$ on  $\frak{G}^{(N)}_\ell$ are equivalent, so one can equivalently work with one or the other, depending on the context. This will be useful in defining the Brownian rough path for example.

\subsection{Definition of a H\"older $p$-rough path}
\label{SubsubsectionDefnRoughPaths}

The relevance of the algebraic framework provided by the $N$-step nilpotent Lie group for the study of smooth paths was first noted by Chen in his seminal work \cite{Chen77}. Indeed, for any $\RR^\ell$-valued smooth path $(x_s)_{s\geq 0}$, the family of iterated integrals 
$$
{\frak X}^N_{ts} := \left(1,x_t-x_s,\int_s^t\int_s^{s_1}dx_{s_2}\otimes dx_{s_1},\dots,\int_{s\leq s_1\leq \cdots\leq s_N\leq t}dx_{s_1}\otimes\cdots\otimes dx_{s_N}\right)
$$
defines for all $0\leq s\leq t$ an element of $T^{(N),1}_\ell$ with the property that if $x_\bullet$ is scaled into $\lambda x_\bullet$ then ${\frak X}^N$ becomes $\delta_\lambda {\frak X}^N$. We actually have ${\frak X}^N_{ts}\in \frak{G}^{(N)}_\ell$. To see that, notice that, as a function of $t$, the function ${\frak X}^N_{ts}$ satisfies the differential equation 
$$
d{\frak X}^N_{ts} = {\frak X}^N_{ts}\,dx_t,
$$
in $T^{(N)}_\ell$ driven by the $\RR^\ell$-valued smooth control $x$, so it defines a $\frak{G}^{(N)}_\ell$-valued path as an integral curve of a field of tangent vectors. The above differential equation also makes it clear the we have the following \textit{Chen relations} \index{Chen's identity}
$$
{\frak X}^N_{ts}  = {\frak X}^N_{us}\,{\frak X}^N_{tu},
$$
for all $0\leq s\leq u\leq t$, which is nothing but the "flow" property for ordinary differential equation solutions; they imply in particular the identity
$$
{\frak X}^N_{ts} = \Big({\frak X}^N_{s0}\Big)^{-1}{\frak X}^N_{t0},
$$ 
Rough paths and weak geometric rough paths are somehow an abstract version of this family of iterated integrals. \index{H\"older $p$-rough path}\index{Weak geometric H\"older $p$-rough path}

\begin{defn}
\label{DefnRoughPath}
Let $2\leq p$. A \emph{\textbf{H\"older $p$-rough path on $[0,T]$}} is a $T_\ell^{([p]),1}$-valued path ${\bfX} : t\in [0,T]\mapsto 1\oplus X^1_t\oplus X^2_t\oplus\cdots\oplus X^{[p]}_t$ \st 
\begin{equation}
\label{ConditionsHolder}
\big\|X^i\big\|_{\frac{i}{p}} := \underset{0\leq s<t\leq T}{\sup}\,\frac{|X^i_{ts}|}{|t-s|^{\frac{i}{p}}}<\infty,
 \end{equation}
for all $i=1\dots [p]$, where we set ${\bf X}_{ts} := {\bf X}_s^{-1}{\bf X}_t$. We define the norm of $\bfX$ to be 
\begin{equation}
\|{\bfX}\| := \underset{i=1\dots [p]}{\max}\;\big\|X^i\big\|_{\frac{i}{p}},
\end{equation}
and a distance $d(\bf X,\bf Y) = \|\bfX-\bfY\|$ on the set of H\"older $p$-rough path. A \emph{\textbf{H\"older weak geometric $p$-rough path on $[0,T]$}} is a $\frak{G}_\ell^{[p]}$-valued $p$-rough path.
\end{defn}

\ssk

\noindent So a (weak geometric) H\"older $p$-rough path is in a way nothing but a ($\frak{G}^{(N)}_\ell$ or) $T^{(N),1}_\ell$-valued $\frac{1}{p}$-H\"older continuous path, for the $\|\cdot\|$-norm introduced above and the use of ${\bfX}_s^{-1}{\bfX}_t$ in place of the usual ${\bfX}_t-{\bfX}_s$. Note that the Chen relation
$$
{\bfX}_{ts} = {\bfX}_{us}{\bfX}_{tu}
$$
is granted by the definition of ${\bfX}_{ts}={\bfX}_s^{-1}{\bfX}_t$. 

\medskip

For $2\leq p<3$, Chen's relation is equivalent to \index{Chen's identity!with $2<p<3$} \vspace{0.1cm}
\begin{itemize}
   \item[{\bf (i)}] $X^1_{ts} = X^1_{tu}+X^1_{us}$, \vspace{0.1cm}
   \item[{\bf (ii)}]  $X^2_{ts} = X^2_{tu} + X^1_{us}\otimes X^1_{tu} + X^2_{us}$. \vspace{0.1cm}
\end{itemize}
Condition \textbf{(i)} means that $X^1_{ts}=X^1_{t0}-X^1_{s0}$ represents the increment of the $\RR^\ell$-valued path $\big(X^1_{r0}\big)_{0\leq r\leq T}$. Condition \textbf{(ii)} is nothing but the analogue of the elementary property $\int_s^t\int_s^r =\int_s^u\int_s^r + \int_u^t\int_s^u + \int_u^t\int_u^r$, satisfied by any reasonable notion of integral on $\RR$ that satisfies the Chasles relation
$$
\int_s^t = \int_s^u+\int_u^t.
$$
This remark justifies thinking of the $\big(\RR^\ell\otimes\RR^\ell\big)$-part of a rough path as a kind of iterated integral of $X^1$ against itself, although this hypothetical iterated integral does not make sense in itself for lack of an integration operation for a general H\"older path in $\RR^\ell$. In that setting, a $p$-rough path $\bfX$ is a weak geometric $p$-rough path iff the symmetric part of $X^2_{ts}$ is $\frac{1}{2}X^1_{ts}\otimes X^1_{ts}$, for all $0\leq s\leq t\leq T$.

\ssk

Note that the space of H\"older $p$-rough paths is not a vector space; this prevents the use of the classical Banach space calculus.

\medskip

It is clear that considering the iterated integrals of any given smooth path defines a H\"older $p$-rough path above it, for any $p\geq 2$. This \textit{lift is not unique}, as if we are given a H\"older $p$-rough path ${\bfX} = \big(X^1,X^2\big)$, with $2\leq p<3$ say, and any $\frac{2}{p}$-H\"older continuous $(\RR^\ell)^{\otimes 2}$-valued path $\big(M_t\big)_{0\leq t\leq 1}$, we define a new rough path setting $M^2_{ts} = M_t-M_s$, and 
$$
{\bfX}'_{ts} = \big(X^1_{ts},X^2_{ts}+M^2_{ts}\big)
$$
for all $0\leq s\leq t\leq 1$. Relations \textbf{(i)} and \textbf{(ii)} above are indeed easily checked.

\ssk

Last, note that a H\"older $p$-rough path is also a H\"older $q$-rough path for any $p<q<[p]+1$. \vspace{0.2cm}

\section{The metric space of H\"older $p$-rough paths}
\label{SubsectionRoughPathSpace}

The distance $d$ defined in definition \ref{DefnRoughPath} is actually not a distance since only the increments $\bfX_{ts}-\bfY_{ts}$ are taken into account. We define a proper metric on the set of all H\"older $p$-rough paths setting
$$
\overline{d}({\bfX},{\bfY}) = \big|X^1_0-Y^1_0\big| + d({\bfX},{\bfY}).
$$

\begin{prop}
\label{PropMetricSpaceRP}
The metric $\overline{d}$ turns the set of all H\"older $p$-rough paths into a (non-separable) complete metric space.
\end{prop}

\begin{Dem}
Given a Cauchy sequence of H\"older $p$-rough paths ${^{(n)}\bfX}$, there is no loss of generality in supposing that their first level starts from the same point in $\RR^\ell$. It follows from the uniform H\"older bounds for $\Big\|{^{(n)}X}^i_{ts}-{^{(m)}X}^i_{ts}\Big\|_{\frac{i}{p}}$, and (an easily proved version of) Ascoli-Arzela theorem (for 2-parameter maps) that ${^{(n)}X}$ converges uniformly to some H\"older $p$-rough path $\bfX$. To prove the convergence of ${^{(n)}X}$ to $\bfX$ in $d$-distance, it suffices to send $m$ to infinity in the inequality
$$
\Big|{^{(n)}X}^i_{ts}-{^{(m)}X}^i_{ts}\Big|\leq \ep\,|t-s|^\frac{i}{p},
$$
which holds for all $n,m$ bigger than some $N_\ep$, uniformly with respect to $0\leq s\leq t\leq 1$.

An uncountable family of $\RR^\ell$-valued $\frac{1}{p}$-H\"older continuous functions at pairwise $\frac{1}{p}$-H\"older distance bounded below by a positive constant is constructed in example 5.28 of \cite{FVBook}. As the set of all first levels of the set of H\"older $p$-rough paths is a subset of the set of $\RR^\ell$-valued $\frac{1}{p}$-H\"older paths, this examples implies the non-separability of set of all H\"older $p$-rough paths.
\end{Dem}

The following \textit{interpolation result} will be useful in several places to prove rough paths convergence results at a cheap price. It roughly says that uniform bounds in a strong sense together with a convergence property in a weak sense are sufficient to prove a convergence result in a mild sense. \index{Interpolation}

\begin{prop}
\label{PropInterpolation}
Assume ${^{(n)}\bfX}$ is a sequence of H\"older $p$-rough paths with uniform bounds
\begin{equation}
\label{EqUniformBoundsInterpolation}
\sup_n\,\big\|{^{(n)}X}\big\| \leq C<\infty,
\end{equation}
which converge pointwise, in the sense that ${^{(n)}{\bfX}}_{ts}$ converges to some ${\bfX}_{ts}$ for each $0\leq s\leq t\leq 1$. Then the limit object $\bfX$ is a H\"older $p$-rough path, and ${^{(n)}\bfX}$ converges to $\bfX$ as a H\"older $q$-rough path, for any $p<q< [p]+1$.
\end{prop}

\begin{Dem}
(Following the solution of exercice 2.9 in \cite{FH14}) The fact that $\bfX$ is a H\"older $p$-rough path is a direct consequence of the uniform bounds \eqref{EqUniformBoundsInterpolation} and pointwise convergence:
$$
\big|X^i_{ts}\big| = \lim_n\Big|{^{(n)}X}^i_{ts}\Big| \leq C|t-s|^\frac{i}{p}.
$$
Would the convergence of ${^{(n)}\bfX}$ to $\bfX$ be uniform, we could find a sequence $\ep_n$ decreasing to $0$, such that, uniformly in $s,t$,
$$
\Big|X^i_{ts}-{^{(n)}X}^i_{ts}\Big| \leq \ep_n, \quad \Big|X^i_{ts}-{^{(n)}X}^i_{ts}\Big|\leq 2C|t-s|^\frac{i}{p}.
$$
Using the geometric interpolation $a\wedge b\leq a^{1-\theta}b^\theta$, with $\theta = \frac{p}{q}<1$, we would have 
$$
\Big|X^i_{ts}-{^{(n)}X}^i_{ts}\Big| \leq \ep_n^{1-\frac{p}{q}}|t-s|^\frac{i}{q},
$$
which entails the convergence result as a H\"older $q$-rough path.

We proceed as follow to see that pointwise convergence suffices to get the result. Given a partition $\pi$ of $[0,1]$ and any $0\leq s\leq t\leq 1$, denote by  $\overline{s}, \overline{t}$ the nearest points in $\pi$ to $s$ and $t$ respectively. Writing 
\begin{equation}
\label{EqEstimateInterpolationProp}
d\Big({\bfX}_{ts},{^{(n)}\bfX}_{ts}\Big) \leq d\big({\bfX}_{ts},{\bfX}_{\overline{t}\overline{s}}\big)  + d\Big({\bfX}_{\overline{t}\overline{s}},{^{(n)}\bfX}_{\overline{t}\overline{s}}\Big) + d\Big({^{(n)}\bfX}_{\overline{t}\overline{s}},{^{(n)}\bfX}_{ts}\Big)
\end{equation}
and the fact that 
$$
{\bfX}_{\overline{t}\overline{s}} = {\bfX}_{s\overline{s}} {\bfX}_{ts} {\bfX}_{\overline{t} t}, \quad {^{(n)}\bfX}_{\overline{t}\overline{s}} = {^{(n)}\bfX}_{s\overline{s}} {^{(n)}\bfX}_{ts} {^{(n)}\bfX}_{\overline{t} t}
$$
and the uniform estimate \eqref{EqUniformBoundsInterpolation} to see that the first and third terms in the above upper bound can be made arbitrarily small by choosing a partition with a small enough mesh, uniformly in $s,t$ and $n$. The second term is dealt with the pointwise convergence assumption as it involves only finitely many points once the partition $\pi$ has been chosen as above.
\end{Dem}

\bigskip

\section{Controlled paths and rough integral}
\label{SubsectionControlledPath}

It will be the set of H\"older weak geometric $p$-rough paths that will play the main role in the sequel, as a set of driving signals in rough differential equations. Unlike the space of H\"older $p$-rough paths, this set is not a linear space, nor even a(n infinite dimensional) manifold, simply a metric space for which none of the classical tools of Banach space calculus can be applied in a straightforward way. It is fortunate, however, that Gubinelli developped in \cite{Gubinelli04} some intermediate spaces of rough paths which have some Banach space structure built in. Their definition requires the introduction of the notion of controlled path, which will appear as the good notion of integrand in the definition of a rough integral. The formalism of this section will be used in part IV of the course on stochastic analysis, where we shall recast the theory of rough differential equations developped in part III of the course in therms of Taylor-Euler expansion properties, using the rough integral introduced in this section. A reference H\"older $p$-rough path $\bfX$ is fixed throughout this section, with $2<p<3$. \index{Controlled path}

\ssk

\begin{defn}
An $\RR^d$-valued path $z_\bullet$ is said to be a {\textit{\textbf{path controlled by $\bfX$}}} if its increments $Z_{ts} = z_t-z_s$, satisfy 
$$
Z_{ts} := Z'_sX_{ts} + R_{ts},
$$
for all $0\leq s\leq t\leq 1$, for some $\textrm{\emph{L}}\big(\RR^\ell,\RR^d\big)$-valued $\frac{1}{p}$-Lipschitz map $Z'_\bullet$, and some $\RR^d$-valued $\frac{2}{p}$-Lipschitz map $R$. 
\end{defn}

The following example shows that a controlled path $z_\bullet$ may have several derivatives $Z'_\bullet$. Choose a H\"older $p$-rough path $\bfX$ with $X_{ts}=(t-s)v$, for some fixed vector $v\in\RR^\ell$. For any path $z_\bullet$  controlled by $\bfX$, one can write
$$
Z_{ts} = Z'_s (t-s)v + R^{Z'}_{ts}
$$
for \textit{any} choice of $\frac{1}{p}$-H\"older function $Z'$ as the term $Z'_s(t-s)v$ can always be inserted in the remainder. So, strictly speaking, a controlled path is a \textit{pair} $(z,Z')$ with the above properties. We sometimes abuse notations and talk of the controlled path $z_\bullet$.

\ssk

Using the notation $\|\cdot\|_\alpha$ to denote the $\alpha$-H\"older norm of a 1 or 2-indices map, it is straightforward to see that one defines a complete metric on the set of $\RR^d$-valued paths $(z,Z')$ controlled by $\bfX$, together with their derivative, setting 
$$
\big\|(z,Z')\big\| := \|Z'\|_\frac{1}{p} + \|R\|_\frac{2}{p} + \big|z_0\big|,
$$
where $\|\cdot\|_\alpha$ stands for the $\alpha$-Lipschitz norm. It is elementary to see that the image of a controlled path $z_\bullet$ by an $\RR^n$-valued $\mcC^1$ map $F$ on $\RR^d$ is a controlled path $F(z_\bullet)$ with derivative $D_{z_t}F\circ Z'_t$.

\ssk

The definition of a controlled path involves only the first level of the rough path $\bfX$. The reference to $\bfX$ itself comes from the following crucial property of controlled paths: they admit an natural lift into a H\"older $p$-rough path, whose definition involves all of $\bfX$. Given two linear maps $A,B\in \textrm{L}(\RR^\ell,\RR^d)$, and any $a,b\in\RR^\ell$, we set
$$
\big(A\otimes B\big)(a\otimes b) := (Aa)\otimes (Bb).
$$

\begin{prop}
\label{PropDefnRoughIntegral}
Let $(z,Z')$ be an $\RR^d$-valued path controlled by $\bfX$. We define an almost-additive map setting
$$
\mu_{ts} := z_s\otimes Z_{ts} + Z'_s\otimes Z'_s \bbX_{ts},
$$
for all $0\leq s\leq t\leq 1$. Its associated additive map $\varphi_{ts}$ is denoted by
$$
{\bbZ}_{ts} =: \int_s^t z_{us}\otimes dz_u.
$$
The pair $(z,{\bbZ})$ is a H\"older $p$-rough path. 
\end{prop}

\begin{Dem}
An elementary computation using Chen's relation $\bbX_{ts} = \bbX_{tu} + \bbX_{us} + X_{us}\otimes X_{tu}$, for any $0\leq s\leq u\leq t\leq 1$, gives
\begin{equation*}
\begin{split}
\big(\mu_{tu}+\mu_{us}\big) - \mu_{ts} &= Z_{us}\otimes Z_{tu} + \big(Z'_u\otimes Z'_u - Z'_s\otimes Z'_s\big)\bbX_{tu} - \big(Z'_s\otimes Z'_s\big)X_{us}\otimes X_{tu} \\
&= Z_{us}\otimes \big(Z_{tu}-Z'_sX_{tu}\big) + O\big(|t-s|^\frac{3}{p}\big) \\
&= Z_{us}\otimes \big(\big(Z'_u-Z'_s\big)\otimes X_{tu}\big) + O\big(|t-s|^\frac{3}{p}\big) = O\big(|t-s|^\frac{3}{p}\big).
\end{split}
\end{equation*}
The $\frac{2}{p}$-H\"older character of $\ZZ_{ts}$ is immediate from the identity
$$
\ZZ_{ts} = \mu_{ts} + O\big(|t-s|^\frac{3}{p}\big),
$$
while Chen's relations are straightforward to check.
\end{Dem}

\ssk

Let $\big(\textrm{F}_t\big)_{0\leq s\leq t\leq 1}$ be an $\textrm{L}(\RR^\ell,\RR^n)$-valued path controlled by $\bfX$. The same computation as above shows that we define an almost-additive map by the formula
$$
\textrm{F}_sX_{ts} + \textrm{F}'_s\bbX_{ts};
$$
it associated additive map is denoted by 
$$
\int_s^t \textrm{F}d{\bfX},
$$
and called the \textbf{rough integral of} $\textrm{F}$ \textbf{with respect to $\bfX$}. \index{Rough integral}

\medskip

There exists, for any $p\geq 3$, a notion of path controlled by a H\"older $p$-rough path. However, the good algebraic setting to work with these objects is not the tensor algebra introduced in this part of the course, but a Hopf algebra of labelled trees. Rough paths are replaced in that setting by \textit{branched rough paths}. This somewhat heavier algebraic setting makes the use of branched rough paths not so convenient. Fortunately, we shall only need the results contained in this section to investigate stochastic differential equations driven by Brownian motion in part IV of the course. See Gubinelli's original work \cite{GubinelliBranched} on the subject, or the nice account \cite{HairerKelly} given by Hairer and Kelly to get some more insights on this question.

\bigskip

\section{Exercices on rough paths}
\label{SubsectionExercicesRoughPaths}

Exercice 7 presents a fundamental result of Lyons of primary importance in the original formulation of the theory. It essentially means that a $p$-rough path has a unique extension into a $q$-rough path, for any $q\geq [p]+1$. The extension of a rough path to all higher degrees defines an object called the \textit{signature} of the rough path, whose importance for real life data analysis is actively investigated presently. Exercices 8 and 9 emphasize the fact that rough paths naturally appear in highly oscillating systems as a class of  controllers (this fact will appear clearly after reading Part III of the course). Exercises 10 an 11 deal with the question of lifting a path or a pair of rough paths into a single rough path.

\bigskip

{\small 

{\bf 7. Lyons' extension theorem \cite{Lyons97}.} \index{Lyons' extension theorem} Let $n$ be a positive integer. A $T^{n,1}_\ell$-valued map ${\bfX} = \big({\bfX}_{ts}\big)_{0\leq s\leq t\leq 1}$ is said to be \textit{multiplicative} if we have
$$
{\bfX}_{ts} = {\bfX}_{us}{\bfX}_{tu}
$$
for all $0\leq s\leq u\leq t\leq 1$. It is said to be \textit{almost-multiplicative} if we have
$$
\left| X^k_{ts} - \big(X_{us} X_{tu}\big)^k\right|\leq c\,|t-s|^{ka}
$$
for all $0\leq s\leq u\leq t\leq 1$ and $0\leq k\leq n$, for some positive constants $c$ and $a>1$; the notation $X^k$ stands here for the $(\RR^\ell)^{\otimes k}$-component of an element $X$ of $T^{n,1}_\ell$. Prove that if $\bfX$ is a $T^{(n),1}_\ell$-valued multiplicative map and $Y^{n+1}_{ts}$ is a continuous $(\RR^d)^{\otimes (n+1)}$-valued map such that the $T^{n+1,1}_\ell$-valued map
$$
Y := \big(1,X^1,\dots,X^n,Y^{n+1}\big)
$$
is almost-multiplicative, then there exists a unique $(\RR^d)^{\otimes (n+1)}$-valued map $X^{n+1}_{ts}$ with 
$$
\Big|X^{n+1}_{ts} - Y^{n+1}_{ts}\Big| \leq c_1\,|t-s|^{(n+1)a}
$$
for some positive constant $c_1$, such that
$$
Z := \Big(1,X^1,\dots,X^n, X^{n+1}\Big) 
$$
is a $T^{n+1,1}_\ell$-valued multiplicative map. 

Starting from a H\"older $p$-rough path $\bfX$ and $Y^{[p]+1}=0$, one can apply iteratively the above procedure to extend uniquely $\bfX$ into a H\"older $q$-rough path, for any $q\geq [p]+1$, in a consistent way. This provides a $T^{\infty,1}_\ell$-valued extension of $\bfX$ called its \textbf{signature}. \vspace{0.3cm}

{\bf 8. Pure area rough path.} \index{Pure area rough path}  Let $x^n$ be the $\RR^2$-valued path defined in complex notations by the formula
$$
x^n_t = \frac{1}{n}\exp\big(2i\pi n^2 t\big),
$$
for $0\leq t\leq 1$. Let $2<p<3$ be given. \vspace{0.1cm}

{\bf $\quad$ a)} Show that the natural lift ${\bfX}^n = \big(x^n,\bbX^n\big)$ of $x^n$ to a H\"older $p$-rough path converges pointwise to the H\"older $p$-rough path ${\bfX} = (X,\bbX)$ with $X=0$ and 
$$
\bbX_{ts} = \pi\,(t-s)\begin{pmatrix} 0 & 1 \\ -1 & 0\end{pmatrix}. \vspace{0.1cm}
$$

{\bf $\quad$ b)} Prove the uniform bounds $\sup_n\,\big\|x^n\big\|_{\frac{1}{2}}<\infty$ and $\sup_n\,\big\|\bbX^n\big\|_1<\infty$. \vspace{0.1cm}

{\bf $\quad$ c)} Conclude by interpolation that the convergence of ${\bfX}^n$ to $\bfX$ takes place in the space of H\"older $p$-rough paths. \vspace{0.3cm}

{\bf 9. Wild oscillations.} Find a widely oscillating piecewise smooth path converging to $(0,0,t\textrm{I})$ in the space of H\"older $p$-rough paths, for $3<p<4$. The letter I stands here for the element of $(\RR^\ell)^{\otimes 3}$ given in the canonical basis by $\textrm{I}_{ijk}=\delta_{ij}\delta_{jk}$. \vspace{0.3cm}

{\bf 10. Lifting $\alpha$-H\"older paths to rough paths, for $\alpha>\frac{1}{2}$.} Show that using the Young integral defined in exercice 5 one can lift any $\alpha$-H\"older paths, with $\alpha>\frac{1}{2}$, into a H\"older $p$-rough path, for any $p\geq 2$. \vspace{0.3cm}


{\bf 11. Pairing two rough paths.} The problem we address in this exercice is the following. \textit{Given two rough paths defined on some (different) spaces, are these two rough paths "pieces" of a higher dimensional rough path?} This is a non-trivial question, when formulated in this generality, due to the fact that there is no canonical way of constructing the cross-iterated integrals between the two rough paths. However, this question has a simple answer when one of the two rough paths is actually a sufficiently regular H\"older path (this exercice), or when one can use probabilistic arguments to construct the missing iterated integrals (exercice 19).

\ssk

Let $\bfX$ be a H\"older $p$-rough path over $\RR^\ell$, with $2< p<3$, and $h$  be a $\frac{1}{q}$-H\"older $\RR^d$-valued path, with $\frac{1}{p}+\frac{1}{q}>1$; so in particular $\frac{1}{q}>\frac{1}{2}$. We describe an element of $(\RR^\ell\times\RR^d)^{\otimes 2}$ as a $2\times 2$ matrix $\begin{pmatrix} A & C \\ B & D \end{pmatrix}$, with $A$ of size $\ell\times \ell$, $B$ of size $\ell\times d$, $C$ of size $d\times \ell$ and $D$ of size $d\times d$. Show that one defines a rough path ${\bfZ} = (Z,\ZZ)$ over $\RR^\ell\times\RR^d$ setting
$$
Z_{ts} = \big(X_{ts},h_t-h_s\big)
$$
and 
$$
\ZZ_{ts} = \begin{pmatrix} A_{ts} & C_{ts} \\ B_{ts} & D_{ts} \end{pmatrix},
$$
with $A_{ts}=\bbX_{ts}$, $D_{ts}=\int_s^th_{us}\otimes dh_u$ and $B_{ts}=\int_s^t h_{us}\otimes dX_u$, $C_{ts}=\int_s^t X_{us}\otimes dh_u$, where $D_{ts}$ and $C_{ts}$ are Young integrals and the integral $B_{ts}$ is \textit{defined} by the integration by parts formula
$$
\int_s^t h_{us}\otimes dX_u := h_{ts}\otimes X_{ts} - \int_s^t dh_u\otimes X_{us}.
$$
We say that $\bfZ$ is a \textbf{pairing of $\bfX$ and $h$}. Would could you possibly pair them if $\bfX$ were a H\"older $p$-rough path over $\RR^\ell$, without any rstriction on $p>2$, and the condition $\frac{1}{p}+\frac{1}{q}>1$ still holds?
}

\vfill


\chapter{Flows driven by rough paths}
\label{SectionFlowsRoughPaths}

\medskip

\todo[inline, backgroundcolor=white, bordercolor=black]{Guide for this chapter}

\medskip

We have seen in part I of the course that a $\mcC^1$-approximate flow on a Banach space E defines a unique flow $\varphi = \big(\varphi_{ts}\big)_{0\leq s\leq t\leq 1}$ on E such that the inequality
\begin{equation}
\label{EqMuApproxVarphiBis}
\big\|\varphi_{ts}-\mu_{ts}\big\|_\infty \leq c |t-s|^a
\end{equation}
holds for some positive constants $c$ and $a>1$, for all $0\leq s\leq t\leq T$ sufficiently close. The construction of $\varphi$ is actually quite explicit, for if we denote by $\mu_{\pi_{ts}}$ the composition of the maps $\mu_{t_{i+1}t_i}$ along the times $t_i$ of a partition $\pi_{ts}$ of an interval $[s,t]$, the map $\mu_{ts}$ satisfies the estimate
\begin{equation}
\label{EqApproxVarphiMuBis}
\big\|\varphi_{ts}-\mu_{\pi_{ts}}\big\|_\infty \leq \frac{2}{1-2^{1-a}}\,c_1^2\,T\,\big|\pi_{ts}\big|^{a-1},
\end{equation}
where $c_1$ is the constant that appears in the definition of a $\mcC^1$-approximate flow
\begin{equation}
\label{EqMuMuBis}
\big\|\mu_{tu}\circ\mu_{us}-\mu_{ts}\big\|_{\mcC^1} \leq c_1 |t-s|^a.
\end{equation}
It follows in particular from equation \eqref{EqMuApproxVarphiBis} that if $\mu$ depends continuously on some metric space-valued parameter $\lambda$, with respect to the $\mcC^0$-topology, and if identity \eqref{EqMuMuBis} holds uniformly for $\lambda$ moving in a bounded set say, then $\varphi$ depends continuously on $\lambda$, as a uniform limit of continuous functions.

\ssk

The point about the machinery of $\mcC^1$-approximate flows is that they actually pop up naturally in a number of situations, under the form of a local in time description of the dynamics under study; nothing else than a kind of Taylor expansion. This was quite clear in exercice 1 on the ordinary controlled differential equation 
\begin{equation}
\label{EqODEBis}
dx_t = V_i(x_t)\,dh^i_t,
\end{equation}
with $\mcC^1$ real-valued controls $h^1,\dots,h^\ell$ and $\mcC^2_b$ vector fields $V_1,\dots,V_\ell$ in $\RR^d$. The 1-step Euler scheme 
$$
\mu_{ts}(x) = x + \big(h^i_t-h^i_s\big)V_i(x)
$$
defines in that case a $\mcC^1$-approximate flow which has the awaited Taylor-type expansion, in the sense that one has
\begin{equation}
\label{EqTaylorExpansionBis}
f\big(\mu_{ts}(x)\big) = f(x) + \big(h^i_t-h^i_s\big)\big(V_if\big)(x) + O\big(|t-s|^{>1}\big)
\end{equation}
for any function $f$ of class $\mcC^2_b$; but $\mu$ fails to be a flow. Its associated flow is not only a flow, it also satisfies equation \eqref{EqTaylorExpansionBis} as a consequence of identity \eqref{EqMuApproxVarphiBis}.

\medskip

We shall proceed in a very similar way to give some meaning and solve the rough differential equation on flows
\begin{equation}
\label{EqRDEFlows}
d\varphi = Vdt + \textrm{F}^\otimes {\bfX}(dt),
\end{equation}
where $V$ is a Lipschitz continuous vector field on E and $\textrm{F}= \big(V_1,\dots,V_\ell\big)$ is a collection of sufficiently regular vector fields on E, and $\bfX$ is a H\"older $p$-rough path over $\RR^\ell$. A \textbf{solution flow }to equation \eqref{EqRDEFlows} will be defined as a flow on E with a \textbf{uniform Taylor-Euler expansion} of the form \index{Taylor-Euler expansion}
\begin{equation}
\label{EqRDETaylorExpansion}
f\big(\varphi_{ts}(x)\big) = f(x) + \sum_{|I|\leq [p]} X^I_{ts}\big(V_If\big)(x) + O\big(|t-s|^{>1}\big),
\end{equation}
where $I = \big(i_1,\dots,i_k\big)\in\llbracket 1,\ell\rrbracket^k$ is a multi-index with size $k\leq [p]$, and $X^I_{ts}$ stands for the coordinates of ${\bfX}_{ts}$ in the canonical basis of $T^{[p],1}_\ell$. The vector field $V_i$ is seen here as a $1^{\textrm{st}}$-order differential operator, and $V_I = V_{i_1}\cdots V_{i_k}$ as the $k^{\textrm{th}}$-order differential operator obtained by applying successively the operators $V_{i_n}$. 

\ssk

For $V=0$ and $\bfX$ the (weak geometric) $p$-rough path canonically associated with an $\RR^\ell$-valued $\mcC^1$ control $h$, with $2\leq p<3$, equation \eqref{EqRDETaylorExpansion} becomes
\begin{equation}
\label{EqRODETaylorExpansion}
f\big(\varphi_{ts}(x)\big) = f(x) + \big(h^i_t-h^i_s\big)\big(V_if\big)(x)  + \left(\int_s^t\int_s^r dh^j_u\,dh^k_r\right)\,\big(V_jV_kf\big)(x) + O\big(|t-s|^{>1}\big),
\end{equation}
which is nothing else than Taylor formula at order 2 for the solution to the ordinary differential equation \eqref{EqODEBis} started at $x$ at time $s$. Condition \eqref{EqRDETaylorExpansion} is a natural analogue of \eqref{EqRODETaylorExpansion} and its higher order analogues.

\medskip

There is actually a simple way of constructing a map $\mu_{ts}$ which satisfies the Euler expansion \eqref{EqRDETaylorExpansion}. It can be defined as the time 1 map associated with an ordinary differential equation constructed form the $V_i$ and their brackets, and where ${\bfX}_{ts}$ appears as a parameter under the form of its logarithm. That these maps $\mu_{ts}$ form a $\mcC^1$-approximate flow will eventually appear as a consequence of the fact that the time 1 map of a differential equation formally behaves as an exponential map, in some algebraic sense.

\medskip

The notationally simpler case of flows driven by weak geometric H\"older $p$-rough paths, with $2\leq p<3$, is first studied in section \ref{SubsectionRDESimpleCase} before studying the general case in section \ref{SubsectionRDEGeneralCase}. The latter case does not present any additional conceptual difficulty, so a reader which who would like to get the core ideas can read section \ref{SubsectionRDESimpleCase} only, or directly go to section \ref{SubsectionRDEGeneralCase}. The two sections have been written with almost similar words on purpose.

\section{Warm up: working with weak geometric H\"older $p$-rough paths, with $2\leq p< 3$}
\label{SubsectionRDESimpleCase}

Let $V$ be a $\mcC^2_b$ vector field on E and $V_1,\dots,V_\ell$ be $\mcC^3_b$ vector fields on E. Let ${\bfX}=(X,\bbX)$ be a H\"older weak geometric $p$-rough path over $\RR^\ell$, with $2\leq p<3$. Let $\mu_{ts}$ be the well-defined time 1 map associated with the ordinary differential equation \index{Bracket of vector fields}
\begin{equation}
\label{EqODEApproachRDE}
\dot y_u = (t-s)V(y_u) + \left(X^i_{ts}V_i+\frac{1}{2}\,\bbX^{jk}_{ts}\big[V_j,V_k\big]\right)(y_u),\quad 0\leq u\leq 1;
\end{equation}
it associates to any $x\in E$ the value at time 1 of the solution of the above equation started from $x$; it is well-defined since $V$ and the $V_i$ are in particular globally Lipschitz. It is a direct consequence of classical results on ordinary differential equations, and of the definition of the topology on the space of H\"older weak geometric $p$-rough paths, that the maps $\mu_{ts}$ depend continuously on $\big((s,t),\bfX\big)$ in the uniform topology, and that 
\begin{equation}
\label{EqC2EstimateMuTs}
\big\|\mu_{ts}-\textrm{Id}\big\|_{\mcC^2} = o_{t-s}(1).
\end{equation}
Also, considering $y_u$ as a function of $x$, it is elementary to see that one has the estimate
\begin{equation}
\label{EqElementaryEstimateYu}
\big\|y_u-\textrm{Id}\big\|_{\mcC^1} \leq c\big(1+\|{\bfX}\|^3\big)|t-s|^{1/p}, \quad 0\leq u\leq 1,
\end{equation}
for some constant depending only on $V$ and the $V_i$.

\subsection{From Taylor expansions to flows driven by rough paths}
\label{SubsubsectionFromTaylorToRDEs}

The next proposition shows that $\mu_{ts}$ has precisely the kind of Taylor-Euler expansion property that we expect from a solution to a rough differential equation, as described in the introduction to that part of the course.

\begin{prop}
\label{PropEstimatesRDEWarmUp}
There exists a positive constant $c$, depending only on $V$ and the $V_i$, such that the inequality
\begin{equation}
\label{EqFundamentalEstimate}
 \Big\|f\circ\mu_{ts} - \Bigl\{f + (t-s)Vf + X^i_{ts}\big(V_if\big) + \bbX^{jk}_{ts}\big(V_jV_kf\big)\Bigr\}\Big\|_\infty \leq c\Big(1+\|\bfX\|^3\Big)\,\|f\|_{\mcC^3}\,|t-s|^{\frac{3}{p}}
\end{equation}
holds for any $f\in\mcC^3_b$. 
\end{prop}

The proof of this proposition and the following one are based on the following elementary identity, obtained by applying twice the identity
$$
f(y_r) = f(x) + (t-s)\int_0^r(Vf)(y_u)\,du + X^i_{ts}\int_0^r\big(V_if\big)(y_u)\,du + \frac{1}{2}\,\bbX^{jk}_{ts}\int_0^r\Big(\big[V_j,V_k\big]f\Big)(y_u)\,du,
$$ 
first to $f$, then to $Vf, V_if$ and $\big[V_j,V_k\big]f$ inside the integrals. One has
\begin{equation*}
\label{EqExactFormulaMuTs}
\begin{split}
f\big(\mu_{ts}(x)\big) &= f(x) + (t-s)\int_0^1(Vf)\big(y_u\big)du + X^i_{ts}\int_0^1 \big(V_if\big)\big(y_{s_1}\big)ds_1 + \frac{1}{2}\,\bbX^{jk}_{ts} \int_0^1 \Big(\big[V_j,V_k\big]f\Big)\big(y_u\big)du \\
&=f(x) + (t-s)\big(Vf\big)(x) + (t-s)\int_0^1 \big\{(Vf)\big(y_u\big)-(Vf)(x)\big\}du \\
&+ X^i_{ts}\big(V_if\big)(x) + (t-s)X^i_{ts}\int_0^1\int_0^{s_1}\big(VV_if\big)\big(y_{s_2}\big)\,ds_2ds_1 \\ 
&\quad\quad + \frac{1}{2}\,X^{i'}_{ts}X^i_{ts}\big(V_{i'}V_if\big)(x) + X^i_{ts}X^{i'}_{ts}\int_0^1\int_0^{s_1}\Big\{\big(V_{i'}V_if\big)\big(y_{s_2}\big)-\big(V_{i'}V_if\big)(x)\Big\}\,ds_2ds_1 \\
&\quad\quad +\frac{1}{2}\,X^i_{ts}\bbX^{jk}_{ts}\int_0^1\int_0^{s_1}\Big(\big[V_j,V_k\big]V_if\Big)\big(y_{s_2}\big)\,ds_2ds_1 \\
&+ \frac{1}{2}\,\bbX^{jk}_{ts} \Big(\big[V_j,V_k\big]f\Big)(x) + \frac{1}{2}\,\bbX^{jk}_{ts}\int_0^1 \Big\{\Big(\big[V_j,V_k\big]f\Big)\big(y_u\big) - \Big(\big[V_j,V_k\big]f\Big)(x)\}\,du.
\end{split}
\end{equation*}
Note that since the H\"older $p$-rough path $\bfX$ is assumed to be weak geometric, the symmetric part of $\bbX_{ts}$ is equal to $\frac{1}{2}X_{ts}\otimes X_{ts}$, so one has
\begin{equation}
\label{EqFMu}
f\big(\mu_{ts}(x)\big)  = f(x) + (t-s)(Vf)(x) + X^i_{ts}\big(V_if\big)(x) + {\bbX}^{jk}_{ts}\big(V_jV_kf\big)(x) + \ep^f_{ts}(x),
\end{equation}
where the remainder $\ep^f_{ts}$ is defined by the formula
\begin{equation*}
\begin{split}
\ep^f_{ts}(x) &:= (t-s)\int_0^1 \big\{(Vf)\big(y_u\big)-(Vf)(x)\big\}du +  (t-s)X^i_{ts}\int_0^1\int_0^{s_1}\big(VV_if\big)\big(y_{s_2}\big)\,ds_2ds_1  \\
&+  X^i_{ts}X^{i'}_{ts}\int_0^1\int_0^{s_1}\Big\{\big(V_{i'}V_if\big)\big(y_{s_2}\big)-\big(V_{i'}V_if\big)(x)\Big\}\,ds_2ds_1  \\
&+ \frac{1}{2}\,X^i_{ts}\bbX^{jk}_{ts}\int_0^1\int_0^{s_1}\Big(\big[V_j,V_k\big]V_if\Big)\big(y_{s_2}\big)\,ds_2ds_1 \\
&+ \frac{1}{2}\,\bbX^{jk}_{ts}\int_0^1 \Big\{\Big(\big[V_j,V_k\big]f\Big)\big(y_u\big) - \Big(\big[V_j,V_k\big]f\Big)(x)\Big\}\,du.
\end{split}
\end{equation*}

\bigskip

\begin{DemPropEstimatesRDEWarmUp}
It is elementary to use estimate \eqref{EqElementaryEstimateYu} and the regularity assumptions on the vector fields $V,V_i$ to see that the remainder $\ep^f_{ts}$ is bounded above by a quantity of the form $c\big(1+\|\bfX\|^3\big)\,\|f\|_{\mcC^3}\,|t-s|^{\frac{3}{p}}$, for some constant depending only on $V$ and the $V_i$.
\end{DemPropEstimatesRDEWarmUp}

\bigskip

A further look at formula \eqref{EqExactFormulaMuTs} and estimate \eqref{EqElementaryEstimateYu} also make it clear that 
\begin{equation}
\label{EqEstimateRemainder}
\Big\|\ep^f_{ts}\Big\|_{\mcC^1} \leq  c\big(1+\|{\bfX}\|^3\big)|t-s|^{\frac{3}{p}},
\end{equation}
for a constant $c$ depending only on $V$ and the $V_i$. This is the key remark for proving the next proposition.

\begin{prop}
\label{PropMuApproximateFlow}
The family $\big(\mu_{ts}\big)_{0\leq s\leq t\leq T}$ forms a $\mcC^1$-approximate flow. 
\end{prop}

\medskip

It will be convenient in the following proof to slightly abuse notations and write $V_I(x)$ for $\big(V_I\textrm{Id}\big)(x)$, for any multi-index $I$ and point $x$.

\medskip

\begin{Dem}
We first use formula \eqref{EqFMu} to write
\begin{equation*}
\mu_{tu}\big(\mu_{us}(x)\big) = \mu_{us}(x) + (t-u)V\big(\mu_{us}(x)\big) + X^i_{tu}V_i\big(\mu_{us}(x)\big) + \bbX^{jk}_{tu}\big(V_jV_k\big)\big(\mu_{us}(x)\big) + \ep^{\textrm{Id}\,;\,[p]}_{tu}\big(\mu_{us}(x)\big).
\end{equation*}
We deal with the term $(t-u)V\big(\mu_{us}(x)\big)$ using estimate \eqref{EqElementaryEstimateYu} and the Lipschitz character of $V$:
$$
\big|(t-u)V\big(\mu_{us}(x)\big)-(t-u)V(x)\big| \leq c\big(1+\|{\bfX}\|^3\big)\,|u-s|^\frac{3}{p}.
$$
The remainder $\ep^{\textrm{Id}}_{tu}\big(\mu_{us}(x)\big)$ has a $\mcC^1$-norm bounded above by $c\big(1+\|{\bfX}\|^3\big)^2|t-u|^\frac{3}{p}$, by the remark preceeding proposition \ref{PropMuApproximateFlow} and the $\mcC^1$-estimate \eqref{EqElementaryEstimateYu} on $\mu_{us}$. We develop $V_i\big(\mu_{us}(x)\big)$ to deal with the term $X^i_{tu} V_i\big(\mu_{us}(x)\big)$. As 
\begin{equation*}
V_i\big(\mu_{us}(x)\big) = V_i(x) + (u-s)\big(VV_i\big)(x) + X^{i'}_{us}\big(V_{i'}V_i\big)(x) + \bbX^{jk}_{us}\big(V_jV_kV_i\big)(x) + \ep^{V_i}_{us}(x)
\end{equation*}
we have
\begin{equation}
\label{EqEstimateViMuUs}
X^i_{tu}V_i\big(\mu_{us}(x)\big) = X^i_{tu}V_i(x) + X^{i'}_{us}X^i_{tu}\big(V_{i'}V_i\big)(x) + \varep^{V_i}_{tu,us}(x),
\end{equation}
where the remainder $\varep^{V_i}_{tu,us}$ has $\mcC^1$-norm bounded above by 
\begin{equation}
\label{EqEstimateVarep}
\Big\|\varep^{V_i}_{tu,us}\Big\|_{\mcC^1} \leq c\big(1+\|{\bfX}\|^3\big)\,|u-s|^\frac{3}{p},
\end{equation}
for a constant $c$ depending only on $V$ and the $V_n$. Set
$$
\varep_{tu,us}(x) = \sum_{i=1}^\ell \varep^{V_i}_{tu,us}(x).
$$
The term $\bbX^{jk}_{tu}\big(V_jV_k\big)\big(\mu_{us}(x)\big)$ is simply dealt with writing
\begin{equation}
\label{EqEstimateVjkMuUs}
\bbX^{jk}_{tu}\big(V_jV_k\big)\big(\mu_{us}(x)\big) = \bbX^{jk}_{tu}\big(V_jV_k\big)(x) + \bbX^{jk}_{tu}\Big\{\big(V_jV_k\big)\big(\mu_{us}(x)\big) - \bbX^{jk}_{tu}\big(V_jV_k\big)(x)\Big\},
\end{equation}
and using estimate \eqref{EqElementaryEstimateYu} and the $\mcC^1_b$ character of $V_jV_k$ to see that the last term on the right hand side has a $\mcC^1$-norm bounded above by $c\big(1+\|{\bfX}\|^3\big)\,|u-s|^\frac{3}{p}$. All together, this gives
\begin{equation*}
\begin{split}
\mu_{tu}\big(\mu_{us}(x)\big) &= \mu_{us}(x) + (t-u)V(x) + X^i_{tu}V_i(x) + X^{i'}_{us}X^i_{tu}\big(V_{i'}V_i\big)(x) + \bbX^{jk}_{tu}\big(V_jV_k\big)(x) + \varep_{tu,us}(x) \\
&= x + (u-s)V(x) + X^i_{us}V_i(x) +\bbX^{jk}_{us}\big(V_jV_k\big)(x) + \ep^{\textrm{Id}}_{us}(x) + (\cdots) \\
&= x + (t-s)V(x) + X^i_{ts}V_i(x) +\bbX^{jk}_{ts}\big(V_jV_k\big)(x) + \ep^{\textrm{Id}}_{us}(x) + \varep_{tu,us}(x) \\
&= \mu_{ts}(x) + \ep^{\textrm{Id}}_{us}(x) + \varep_{tu,us}(x),
\end{split}
\end{equation*}
so it follows from estimates \eqref{EqEstimateRemainder} and \eqref{EqEstimateVarep} that $\mu$ is indeed a $\mcC^1$-approximate flow.
\end{Dem}

\medskip

The above proof makes it clear that one can take for constant $c_1$ in the $\mcC^1$-approximate flow property \eqref{EqMuMu} for $\mu$ the constant $c\big(1+\|\bfX\|^3\big)$, for a constant $c$ depending only on $V$ and the $V_i$.

\medskip

Recalling proposition \ref{PropEstimatesRDEWarmUp} describing the maps $\mu_{ts}$ in terms of Euler expansion, the following definition of a solution flow to a rough differential equation is to be thought of as defining a notion of solution in terms of uniform Euler expansion
$$
\Big\|f\circ\varphi_{ts}-\Big\{f+X^i_{ts}V_if+\bbX^{jk}_{ts}V_jV_kf\Big\}\Big\|_\infty \leq c\,|t-s|^{>1}.
$$

\medskip

\begin{defn}
\label{DefnGeneralRDESolution}
A \textbf{\emph{flow}} $\big(\varphi_{ts}\big)_{\,0\leq s\leq t\leq T}$ is said to \textbf{\emph{solve the rough differential equation}} \index{Solution flow to a rough differential equation} \index{Rough differential equation}
\begin{equation}
\label{RDEGeneral}
d\varphi = Vdt + \textrm{\emph{F}}^\otimes\,\bfX(dt)
\end{equation}
if there exists a constant $a>1$ independent of $\bfX$ and two possibly $\bfX$-dependent positive constants $\delta$ and $c$ such that
\begin{equation}
\label{DefnSolRDEGeneral}
\big\|\varphi_{ts}-\mu_{ts}\big\|_\infty \leq c\,|t-s|^a
\end{equation}
holds for all $0\leq s\leq t\leq T$ with $t-s\leq\delta$.
\end{defn}

If for instance $\bfX$ is the weak geometric H\"older $p$-rough path canonically associated with an $\RR^\ell$-valued piecewise smooth path $h$, it follows from exercice 1, and the fact that the iterated integral $\int_s^t\int_s^r dh_u\otimes dh_r$ has size $|t-s|^2$, that the solution flow to the rough differential equation 
$$
d\varphi = Vdt + \textrm{F}^\otimes {\bfX}(dt)
$$
is the flow associated with the ordinary differential equation 
$$
\dot y_t = V(y_t)dt+V_i(y_t)\,dh^i_t. 
$$
The following well-posedness result follows directly from theorem \ref{ThmConstructingFlows} on $\mcC^1$-approximate flows and proposition \ref{PropMuApproximateFlow}.

\begin{thm}
\label{ThmMainResultGeneral}
The rough differential equation on flows 
$$
d\varphi = Vdt + \textrm{\emph{F}}^\otimes\,\bfX(dt)
$$ 
has a unique solution flow; it takes values in the space of uniformly Lipschitz continuous homeomorphisms of \emph{E} with uniformly Lipschitz continuous inverses, and depends continuously on $\bfX$.
\end{thm}

\ssk

\begin{Dem}
Applying theorem \ref{ThmConstructingFlows} on $\mcC^1$-approximate flows to $\mu$ we obtain the existence of a unique flow $\varphi$ satisfying condition \eqref{DefnSolRDEGeneral}, for $\delta$ small enough; it further satisfies the inequality
\begin{equation}
\label{EqApproxPhiMu}
\|\varphi_{ts}-\mu_{\pi_{ts}}\|_\infty \leq c\big(1+\|\bfX\|^3\big)^2T\,\big|\pi_{ts}\big|^{a-1},
\end{equation}
for any partition $\pi_{ts}$ of $[s,t]\subset [0,T]$ of mesh $\big|\pi_{ts}\big|\leq \delta$, as a consequence of inequality \eqref{EqApproxVarphiMu}. As this bound is uniform in $(s,t)$, and for $\bfX$ in a bounded set of the space of weak geometric H\"older $p$-rough paths, and since each map $\mu_{\pi_{ts}}$ is a continuous function of $\big((s,t),\bfX\big)$, the flow $\varphi$ depends continuously on $\big((s,t),\bfX\big)$.

To prove that $\varphi$ is a homeomorphism, note that, with the notations of part I of the course,
$$
\Big(\mu^{(n)}_{ts}\Big)^{-1} = \mu_{s_1s_0}^{-1}\circ\cdots\circ\mu_{s_{2^n}s_{2^n-1}}^{-1}, \quad s_i = s+i2^{-n}(t-s),
$$
can actually be written $\big(\mu^{(n)}_{ts}\big)^{-1} = \overline\mu_{s_{2^n}s_{2^n-1}}\circ\cdots\circ\overline\mu_{s_1s_0}$, for the time 1 map $\overline\mu$ associated with the rough path $\bfX_{t-\bullet}$. As $\overline\mu$ enjoys the same properties as $\mu$, the maps $\big(\mu^{(n)}_{ts}\big)^{-1}$ converge uniformly to some continuous map $\varphi_{ts}^{-1}$ which satisfies by construction $\varphi_{ts}\circ\varphi_{ts}^{-1} = \textrm{Id}$.

\ssk

Recall that proposition \ref{PropSecondStep} provides a uniform control of the Lipschitz norm of the maps $\varphi_{ts}$; the same holds for their inverses in view of the preceeding paragraph. We propagate this property from the set $\big\{(s,t)\in [0,T]^2\,;\,s\leq t, \; t-s\leq \delta\big\}$ to the whole of the $\big\{(s,t)\in [0,T]^2\,;\,s\leq t\big\}$ using the flow property of $\varphi$.
\end{Dem}

\bigskip

\begin{Rems}
\label{RemsThm2p3}
\begin{enumerate}
   \item {\bf Friz-Victoir approach to rough differential equations.} \index{Friz-Victoir approach} The continuity of the solution flow \wrt the driving rough path ${\bf X}$ has the following consequence, which justifies the point of view adopted by Friz and Victoir in their works. Suppose the H{\"o}lder weak geometric $p$-rough path ${\bf X}$ is the limit in the rough path metric of the canonical H{\"o}lder weak geometric $p$-rough paths ${\bf X}^n$ associated with some piecwise smooth $\RR^\ell$-valued paths $(x^n_t)_{0\leq t\leq T}$. We have noticed that the solution flow $\varphi^n$ to the rough differential equation 
$$
d\varphi^n = Vdt + \textrm{\emph{F}}^\otimes {\bfX}^n(dt)
$$
is the flow associated with the ordinary differential equation 
$$
\dot y_u = V(y_u)du+V_i(y_u)\,d(x^n_u)^i. 
$$
As $\|\varphi^n-\varphi\|_\infty=o_n(1)$, from the continuity of the solution flow \wrt the driving rough path, the flow $\varphi$ appears in that case as a uniform limit of the elementary flows $\varphi^n$. A H{\"o}lder weak geometric $p$-rough path with the above property is called a \textit{H{\"o}lder geometric $p$-rough path}; not all H{\"o}lder weak geometric $p$-rough path are H{\"o}lder geometric $p$-rough path \cite{FVGeometric}, although there is little difference. \vspace{0.2cm}
   \item \textbf{Time-inhomogeneous dynamics.} The above results have a straightforward generalization to dynamics driven by a time-dependent bounded drift $V(s;\cdot)$ which is Lipschitz continuous \wrt the time variable and $\mcC^2_b$ \wrt the space variable, uniformly \wrt time, and time-dependent  vector fields $V_i(s;\cdot)$ which are Lipschitz continuous \wrt time, and $\mcC^3_b$ \wrt the space variable, uniformly \wrt time. We define in that case a $\mcC^1$-approximate flow by defining $\mu_{ts}$ as the time 1 map associated with the ordinary differential equation
$$
\dot y_u = (t-s)V(s;y_u) + X^i_{ts}V_i(s;y_u) + \bbX^{jk}_{ts}\big[V_j,V_k\big](y_u),  \quad 0\leq u\leq 1.	
$$
\end{enumerate}
\end{Rems}

\subsection{Classical rough differential equations} 
\label{SubsubsectionClassicalRDEs}
\index{Classical rough differential equations}

In the classical setting of rough differential equations, one is primarily interested in a notion of \textit{solution path}, defined in terms of local Taylor-Euler expansion.

\begin{defn}
\label{DefnPathSolRDE}
A \textbf{path} $\big(z_s\big)_{0\leq s\leq T}$ is said to \emph{\textbf{solve the rough differential equation}}
\begin{equation}
\label{PathRDEGeneral}
dz = Vdt + \textrm{\emph{F}}\,\bfX(dt)
\end{equation}
with initial condition $x$, if $z_0=x$ and there exists a constant $a>1$ independent of $\bfX$, and two possibly $\bfX$-dependent positive constants $\delta$ and $c$, such that
\begin{equation}
\label{DefnSolRDEGeneral}
\Big|f(z_t) - \Big\{f(z_s)+(t-s)(Vf)(z_s)+X^i_{ts}\big(V_if\big)(z_s)+\bbX^{jk}_{ts}\big(V_jV_kf\big)(z_s)\Big\}\Big| \leq c\,\|f\|_{\mcC^3}\,|t-s|^a
\end{equation}
holds for all $0\leq s\leq t\leq T$, with $t-s\leq\delta$, for all $f\in\mcC^3_b$.
\end{defn}

\medskip

\begin{thm}[Lyons' universal limit theorem]
\label{ThmLyonsUniversalLimitThm}
The rough differential equation \eqref{PathRDEGeneral} has a unique solution path; it is a continuous function of $\bfX$ in the uniform norm topology.	\index{Lyons' universal limit theorem}
\end{thm}

\medskip

\begin{Dem}
\noindent \textbf{a) Existence.} It is clear that if $\big(\varphi_{ts}\big)_{0 \leq s\leq t\leq 1 }$ stands for the solution flow to the equation 
$$
d\varphi = Vdt + \textrm{F}^\otimes{\bfX}(dt),
$$
then the path $z_t := \varphi_{t0}(x)$ is a solution path to the rough differential equation \eqref{PathRDEGeneral} with initial condition $x$.

\medskip

\noindent \textbf{b) Uniqueness.} Let agree to denote by $O_c(m)$ a quantity whose norm is bounded above by $c\,m$. Let $\alpha$ stand for the minimum of $\frac{3}{p}$ and the constant $a$ in definition \ref{DefnPathSolRDE}, and let $y_\bullet$ be any other solution path. It satisfies by proposition \ref{PropEstimatesRDEWarmUp} the estimate
$$
\big|y_t - \varphi_{ts}(y_s)\big| \leq c|t-s|^\alpha.
$$
Using the fact that the maps $\varphi_{ts}$ are uniformly Lipschitz continuous, with a Lipschitz constant bounded above by $L$ say, one can write for any $\ep>0$ and any integer $k\leq \frac{T}{\ep}$
\begin{equation*}
\begin{split}
y_{k\ep} &= \varphi_{k\ep,(k-1)\ep}\big(y_{(k-1)\ep}\big) + O_c\big(\ep^\alpha\big) \\
             &= \varphi_{k\ep,(k-1)\ep}\Big(\varphi_{(k-1)\ep,(k-2)\ep}\big(y_{(k-2)\ep}\big) + O_c\big(\ep^\alpha\big) \Big) + O_c\big(\ep^\alpha\big) \\
             &= \varphi_{k\ep,(k-2)\ep}\big(y_{(k-2)\ep}\big) + O_{cL}\big(\ep^\alpha\big) + O_c\big(\ep^\alpha\big),
\end{split}
\end{equation*}
and see by induction that 
\begin{equation*}
\begin{split}
y_{k\ep} &= \varphi_{k\ep,(k-n)\ep}\big(y_{(k-n)\ep}\big) + O_{cL}\big((n-1)\ep^\alpha\big) + O_c\big(\ep^\alpha\big) \\
             &= \varphi_{k\ep,0}(x) + O_{cL}\big(k\ep^\alpha\big) + o_\ep(1) \\
             &= z_{k\ep} + O_{cL}\big(k\ep^\alpha\big) + o_\ep(1).
\end{split}
\end{equation*}
Taking $\ep$ and $k$ so that $k\ep$ converges to some $t\in [0,T]$, we see that $y_t=z_t$, since $\alpha>1$.

\medskip

The continuous dependence of the solution path $z_\bullet$ \wrt $\bfX$ is transfered from $\varphi$ to $z_\bullet$.
\end{Dem}

\medskip

The map that associates to the rough path $\bfX$ the solution to the rough differential equation \eqref{PathRDEGeneral} is called the {\bf Ito map}.

\bigskip

\section{The general case}
\label{SubsectionRDEGeneralCase}

We have defined in the previous section a solution to the rough differential equation
$$
d\varphi = Vdt+\textrm{F}^\otimes{\bfX}(dt),
$$
driven by a weak geometric H\"older $p$-rough path, for $2\leq p<3$, as a flow with $(s,t;x)-$uniform Taylor-Euler expansion of the form
$$
f\big(\varphi_{ts}(x)\big) = f(x) + (t-s)(Vf)(x) + X^i_{ts}\big(V_if\big)(x) + \bbX^{jk}_{ts}\big(V_jV_kf\big)(x) +  O\big(|t-s|^{>1}\big).
$$
The definition of a solution flow in the general case will require from $\varphi$ that it satisfies a similar expansion, of the form
\begin{equation}
\label{EqDefnSolFlowGeneral}
f\big(\varphi_{ts}(x)\big) = f(x) + (t-s)(Vf)(x) + \sum_{|I|\leq [p]} X^I_{ts}\big(V_If\big)(x) + O\big(|t-s|^{>1}\big).
\end{equation}
As in the previous section, we shall obtain $\varphi$ as the unique flow associated with some $\mcC^1$-approximate flow $\big(\mu_{ts}\big)_{0\leq s\leq t\leq 1}$, where $\mu_{ts}$ is the time 1 map associated with an ordinary differential equation constructed from the $V_i$ and their brackets, and $V$ and ${\bfX}_{ts}$. In order to avoid writing expressions with loads of indices (the ${\bfX}^I_{ts}$), I will first introduce in subsection \ref{SubsubsectionDifferentialOperators} a coordinate-free way of working with rough paths and vector fields. A $\mcC^1$-approximate flow with the awaited Euler expansion will be constructed in subsection \ref{SubsubsectionFromTaylorToRDEsGeneral}, leading to a general well-posedness result for rough differential equations on flows.

\bigskip

To make the crucial formula \eqref{EqExactFormulaMuTs} somewhat shorter we assume in this section that $V=0$. The reader is urged to workout by herself/himself the infinitesimal changes that have to be done in what follows in order to work with a non-null drift $V$. From hereon, the vector fields $V_i$ are assumed to be of class $\mcC^{[p]+1}_b$. We denote by $\mcC^{[p]+1}_b(\textrm{E},\textrm{E})$ the set of $\mcC^{[p]+1}_b$ vector fields on E. We denote for by $\pi_k : T^\infty_\ell\rightarrow (\RR^\ell)^k$ the natural projection operator and set $\pi_{\leq k} = \sum_{j\leq k} \pi_j$.

\medskip

\subsection{Differential operators}
\label{SubsubsectionDifferentialOperators}

Let F be a continuous linear map from $\RR^\ell$ to $\mcC^{[p]+1}_b(\textrm{E},\textrm{E})$ -- one usually calls such a map a \textit{vector field valued 1-form on $\RR^\ell$}. For any $v\in\RR^\ell$, we identify the $\mcC^{[p]+1}$ vector field $\textrm{F}(v)$ on E with the first order differential operator
$$
\textrm{F}^\otimes(v) \,:\,  g\in \mcC^1(\textrm{E})\mapsto (D_{\cdot} g)\big(\textrm{F}(v)(\cdot)\big)\in\mcC^0(\textrm{E});
$$
in those terms, we recover the vector field $\textrm{F}(v)$ as $\textrm{F}^\otimes (v)\textrm{Id}$. The map $\textrm{F}^\otimes$ is extended to $T^{[p]+1}_\ell$ by setting 
$$
\textrm{F}^\otimes(1) := \textrm{Id} : \mcC^0(\textrm{E})\mapsto \mcC^0(\textrm{E}),
$$
and defining $\textrm{F}^\otimes (v_1\otimes\cdots\otimes v_k)$, for all $1\leq k\leq [p]+1$ and $v_1\otimes\cdots\otimes v_k\in (\RR^\ell)^{\otimes k}$, as the $k^{\textrm{th}}$-order differential operator from $\mcC^k(\textrm{E})$ to $\mcC^0(\textrm{E})$, defined by the formula
$$
\textrm{F}^\otimes (v_1\otimes\cdots\otimes v_k) := \textrm{F}^\otimes\big(v_1\big)\cdots \textrm{F}^\otimes\big(v_k\big),
$$
and by requiring linearity. So, we have the morphism property
\begin{equation}
\label{EqMorphismF}
\textrm{F}^\otimes(\be)\,\textrm{F}^\otimes(\be') = \textrm{F}^\otimes(\be\be')
\end{equation}
for any $\be,\be'\in T^{[p]+1}_\ell$ with $\be\be'\in T^{[p]+1}_\ell$. This condition on $\be,\be'$ is required for if $\be'=v_1\otimes\cdots\otimes v_k$ with $v_i\in\RR^\ell$, the map $\textrm{F}^\otimes(\be')\textrm{Id}$ from E to itself is $\mcC^{[p]+1-k}_b$, so $\textrm{F}^\otimes(\be)\,\textrm{F}^\otimes(\be')$ only makes sense if $\be\be'\in T^{[p]+1}_\ell$. We also have 
$$
\Big[\textrm{F}^\otimes(\be),\textrm{F}^\otimes(\be')\Big] = \textrm{F}^\otimes\big([\be,\be']\big)
$$
for any $\be,\be'\in T^{[p]+1}_\ell$ with $\be\be'$ and $\be'\be$ in $T^{[p]+1}_\ell$. This implies in particular that $\textrm{F}^\otimes(\Lambda)$ is actually a first order differential operator for any $\Lambda\in\frak{g}^{[p]+1}_\ell$, that is a vector field. Note that for any $\Lambda\in\frak{g}^{[p]+1}_\ell$ and $1\leq k\leq [p]+1$, then $\Lambda^k:=\pi_k(\Lambda)$ is an element of $\frak{g}^{[p]}_\ell$, and the vector field $\textrm{F}^\otimes\big(\Lambda^k\big)\textrm{Id}$ is $\mcC^{[p]+1-k}_b$.

\ssk

We extend $\textrm{F}^\otimes$ to the unrestricted tensor space $T^\infty_\ell$ setting
\begin{equation}
\label{EqExtensionF}
\textrm{F}^\otimes(\be) = \textrm{F}^\otimes\big(\pi_{\leq [p]+1}\be\big)
\end{equation}
for any $\be\in T^\infty_\ell$.

\medskip

Consider as a particular case the map F defined for $u\in\RR^\ell$ by the formula
$$
\textrm{F}(u) = u^i\,V_i(\cdot).
$$
Using the formalism of this paragraph, an Euler expansion of the form 
$$
f\big(\varphi_{ts}(x)\big) = f(x) +\sum_{|I|\leq [p]} X^I_{ts}\big(V_If\big)(x) + O\big(|t-s|^{>1}\big),
$$
as in equation \eqref{EqDefnSolFlowGeneral}, becomes
$$
f\big(\varphi_{ts}(x)\big) = \big(\textrm{F}^\otimes\big({\bfX}_{ts}\big)f\big)(x) + O\big(|t-s|^{>1}\big).
$$

\medskip

\subsection{From Taylor expansions to flows driven by rough paths: bis}
\label{SubsubsectionFromTaylorToRDEsGeneral}

Let $2\leq p$ be given, together with a $\frak{G}^{[p]}_\ell$-valued weak-geometric H\"older $p$-rough path $\bfX$, defined on some time interval $[0,T]$, and some continuous linear map F from $\RR^\ell$ to the set $\mcC^{[p]+1}_b(\textrm{E},\textrm{E})$ of vector fields on E. For any $0\leq s\leq t\leq T$, denote by ${\bf \Lambda}_{ts}$ the logarithm of ${\bfX}_{ts}$,  and let $\mu_{ts}$ stand for the well-defined time $1$ map associated with the ordinary differential equation
\begin{equation}
\label{EqRPODE}
\dot y_u = \textrm{F}^\otimes\big({\bf\Lambda}_{ts}\big)(y_u),\quad 0\leq u\leq 1.
\end{equation}
This equation is indeed an ordinary differential equation since ${\bf\Lambda}_{ts}$ is an element of $\frak{g}^{[p]}_\ell$. For $2\leq p<3$, it reads 
\begin{equation*}
\dot y_u = X^i_{ts}V_i(y_u) + \frac{1}{2}\left(\bbX^{jk}_{ts}+\frac{1}{2}X^j_{ts}X^k_{ts}\right)\big[V_j,V_k\big](y_u),    \quad   0\leq u\leq 1.
\end{equation*}
As the tensor $X_{ts}\otimes X_{ts}$ is symmetric and the map $(j,k)\mapsto \big[V_j,V_k\big]$ is antisymmetric, this equation actually reads
$$
\dot y_u = X^i_{ts}V_i(y_u) + \frac{1}{2}\bbX^{jk}_{ts}\,\big[V_j,V_k\big](y_u),
$$
which is nothing else than equation \eqref{EqODEApproachRDE}, whose time 1 map defined the $\mcC^1$-approximate flow we studied in section \ref{SubsubsectionFromTaylorToRDEs}. 

\ssk

It is a consequence of classical results from ordinary differential equations, and the definition of the norm on the space of weak-geometric H\"older $p$-rough paths, that the solution map $(r,x)\mapsto y_r$, with $y_0=x$, depends continuously on $\big((s,t),\bfX\big)$ in $\mcC^0$-norm, and satisfies the following basic estimate. The next proposition shows that $\mu_{ts}$ has precisely the kind of Taylor-Euler expansion property that we expect from a solution to a rough differential equation.

\begin{equation}
\label{EqRoughODEsEstimate}
\big\|y_r-\textrm{Id}\big\|_{\mcC^1} \leq c\Big(1+\|{\bfX}\|^{[p]}\Big)|t-s|^\frac{1}{p}, \quad 0\leq r\leq 1
\end{equation}

\begin{prop}
\label{PropFundamentalEstimate}
There exists a positive constant $c$, depending only on the $V_i$, \st the inequality
\begin{equation}
\label{EqFundamentalEstimate}
 \Big\|f\circ\mu_{ts} - \textrm{\emph{F}}^\otimes\big({\bfX}_{ts}\big)f\Big\|_\infty \leq c\Big(1+\|{\bfX}\|^{[p]}\Big)\,\|f\|_{\mcC^{[p]+1}}\,|t-s|^{\frac{[p]+1}{p}}
\end{equation}
holds for any $f\in\mcC^{[p]+1}_b(\textrm{E})$. 
\end{prop}

\ssk

In the classical setting of an ordinary differential equation 
$$
\dot z_u = W(z_u), \quad z_0=x,
$$
driven by a $\mcC^1$ vector field $W$ on E, we would get a Taylor expansion formula for $f(z_1)$ from the elementary formula
$$
f(z_1) = \sum_{k=0}^{n-1}\frac{1}{k!}\,\big(W^{\circ k}f\big)(x) +\int_{\Delta_n}\big(W^{\circ n}f\big)\big(y_{s_n}\big)\,ds
$$
obtained by induction, where $W^{\circ 0}f=f$ and $W^{\circ n+1}f := W\big(W^{\circ n}f\big)$, and 
$$
\Delta_n := \big\{(s_1,\dots,s_n)\in [0,T]^n\,;\,s_n\leq \cdots\leq s_1\big\},
$$ 
with the notation $ds$ for $ds_n\dots ds_1$. We proceed along these lines to obtain a similar formula for the solution to the preceeding equation with 
$$
W=\textrm{F}^\otimes\big({\bf\Lambda}_{ts}\big) = \sum_{k_1=0}^{[p]+1} \textrm{F}^\otimes\big(\Lambda^{k_1}_{ts}\big).
$$ 
Some care however is needed to take into account the fact that the vector fields $\textrm{F}^\otimes\big(\Lambda^{k_1}_{ts}\big)$ have different regularity properties.

\ssk

Recall $\textrm{F}^\otimes (0)$ is the null map from $\mcC^0(E)$ to itself and $\pi_0\Lambda=0$ for any $\Lambda\in\frak{g}^{[p]}_\ell$. The proof of this proposition and the following one are based on the elementary identity \eqref{EqExactFormulaMuTs} below, obtained by applying repeatedly the identity
\begin{equation*}
\begin{split}
f\big(y_r\big) &= f(x) + \int_0^r\Big(\textrm{F}^\otimes\big({\bf\Lambda}_{ts}\big)f\Big)(y_u)\,du \\
                      &= f(x) + \sum_{k_1=0}^{[p]+1}\int_0^r\Big(\textrm{F}^\otimes\big({\bf\Lambda}^{k_1}_{ts}\big)f\Big)(y_u)\,du,      \quad 0\leq r\leq 1
\end{split}
\end{equation*}
together with the morphism property \eqref{EqMorphismF}. As emphasized above, the above sum over $k_1$ is needed to take care of the different regularity properties of the maps $\textrm{F}^\otimes\big({\bf\Lambda}^{k_1}_{ts}\big)f$. 
\begin{equation*}
\begin{split}
f\big(\mu_{ts}(x)\big) &= f(x) + \Big(\textrm{F}^\otimes({\bf\Lambda}_{ts})f\Big)(x) + \sum_{k_1+k_2\leq [p]+1} \int_0^1\int_0^{s_1}\Big(\textrm{F}^\otimes({\bf\Lambda}^{k_2}_{ts})\textrm{F}^\otimes({\bf\Lambda}^{k_1}_{ts})f\Big)\big(y_{s_2}\big)\,ds_2\,ds_1 \\
                                  &= f(x) + \Big(\textrm{F}^\otimes({\bf\Lambda}_{ts})f\Big)(x) + \int_0^1\int_0^{s_1}\Big(\textrm{F}^\otimes\big({\bf\Lambda}_{ts}^{\bullet 2}\big)f\big)\big(y_{s_2}\Big)\,ds_2\,ds_1 
\end{split}
\end{equation*}
We use here the notation $\bullet 2$ to denote the multiplication $\Lambda_{ts}^{\bullet 2} = \Lambda_{ts}\Lambda_{ts}$, not to be confused with the second level ${\bf\Lambda}_{ts}^2$ of ${\bf\Lambda}_{ts}$; the product is done here in $T^\infty_\ell$, and definition \eqref{EqExtensionF} used to make sense of $\textrm{F}^\otimes\big({\bf\Lambda}_{ts}^{\bullet 2}\big)f$. Repeating $(n-1)$ times the above procedure in an iterative way, we see that
\begin{equation*}
\begin{split}
f\big(\mu_{ts}(x)\big) &= f(x) + \sum_{k=1}^{n-1} \frac{1}{k!}\,\Big(\textrm{F}^\otimes\big({\bf\Lambda}_{ts}^{\bullet k}\Big)f\big)(x) +\int_{\Delta_n} \Big(\textrm{F}^\otimes\big({\bf\Lambda}_{ts}^{\bullet n}\big)f\Big)\big(y_{s_n}\big)\,ds \\
                                  &= f(x) + \sum_{k=1}^n \frac{1}{k!}\,\Big(\textrm{F}^\otimes\big({\bf\Lambda}_{ts}^{\bullet k}\big)f\Big)(x) +\int_{\Delta_n} \Big\{ \Big(\textrm{F}^\otimes\big({\bf\Lambda}_{ts}^{\bullet n}\big)f\Big)\big(y_{s_n}\big) - \Big(\textrm{F}^\otimes\big({\bf\Lambda}_{ts}^{\bullet n}\big)f\Big)(x)\Big\}\,ds. 
\end{split}
\end{equation*}
Note that $\pi_j{\bf\Lambda}_{ts}^{\bullet n} = 0$, for all $j\leq n-1$, and 
$$
\pi_{\leq [p]} \left(\sum_{k=1}^{[p]} \frac{1}{k!}{\bf\Lambda}_{ts}^{\bullet k}\right) = {\bfX}_{ts};
$$
also $\pi_{\leq [p]}\Big({\bf\Lambda}_{ts}^{\bullet [p]}\Big) = \big(X^1_{ts}\big)^{\otimes {[p]}}$ is of size $|t-s|^\frac{[p]}{p}$. We separate the different terms in the above identity according to their size in $|t-s|$; this leads to the following expression for $f\big(\mu_{ts}(x)\big)$.
{\small 
\begin{equation}
\label{EqExactFormulaMuTs}
\begin{split}
&f(x) + \left(\textrm{F}^\otimes\Big(\pi_{\leq [p]}\Big\{\sum_{k=1}^n \frac{1}{k!}\,\Lambda_{ts}^{\bullet k}\Big\}\Big)f\right)(x) + \int_{\Delta_n} \left\{ \Big(\textrm{F}^\otimes\big(\pi_{\leq [p]} {\bf\Lambda}_{ts}^{\bullet n}\big)f\Big)\big(y_{s_n}\big) - \Big(\textrm{F}^\otimes\big(\pi_{\leq [p]}{\bf\Lambda}_{ts}^{\bullet n}\big)f\Big)(x)\right\}\,ds \\
        &+ \left(\textrm{F}^\otimes\Big(\pi_{[p]+1}\Big\{\sum_{k=1}^n \frac{1}{k!}\,{\bf\Lambda}_{ts}^{\bullet k}\Big\}\Big)f\right)(x) + \int_{\Delta_n} \Big\{ \Big(\textrm{F}^\otimes\Big(\pi_{[p]+1}{\bf\Lambda}_{ts}^{\bullet n}\Big)f\Big)\big(y_{s_n}\big) - \Big(\textrm{F}^\otimes\Big(\pi_{[p]+1}{\bf\Lambda}_{ts}^{\bullet n}\Big)f\Big)(x)\Big\}\,ds
\end{split}
\end{equation}
}
\noindent We denote by $\ep^{f\,;\,n}_{ts}(x)$ the sum of the two terms involving $\pi_{[p]+1}$ in the above line, made up of terms of size at least $|t-s|^\frac{[p]+1}{p}$. Note that for $n=[p]$, the integral term in the first line involves $\pi_{\leq [p]}\Big({\bf\Lambda}_{ts}^{[p]}\Big) = \big(X^1_{ts}\big)^{\otimes {[p]}}$ and the increment $y_{s_n}-x$, of size $|t-s|^\frac{1}{p}$, by estimate \eqref{EqRoughODEsEstimate}, so this term is of size $|t-s|^\frac{[p]+1}{p}$; we include it in $\ep^{f\,;\,[p]}_{ts}(x)$.

\medskip

\begin{DemPropFundamentalEstimate}
Applying the above formula with $n=[p]$, we get the identity
\begin{equation*}
f\big(\mu_{ts}(x)\big) = \Big(\textrm{F}^\otimes\big({\bfX}_{ts}\big)f\Big)(x) + \ep^{f\,;\,[p]}_{ts}(x).
\end{equation*}
It is clear on the formula for $\ep^{f\,;\,[p]}_{ts}(x)$ that its absolute value is bounded above by a constant multiple of $\Big(1+\|{\bfX}\|^{[p]}\Big)|t-s|^\frac{[p]+1}{p}$, for a constant depending only on the data of the problem and $f$ as in \eqref{EqFundamentalEstimate}. 
\end{DemPropFundamentalEstimate}

\medskip

A further look at formula \eqref{EqExactFormulaMuTs} makes it clear that if $2\leq n\leq [p]$, and $f$ is $\mcC^{n+1}_b$, the estimate 
\begin{equation}
\label{EqEstimateC1SizeEpsilon}
\Big\|\ep^{f\,;\,n}_{ts}\Big\|_{\mcC^1} \leq c\Big(1+\|{\bfX}\|^{[p]}\Big)\,\|f\|_{\mcC^{n+1}}|t-s|^\frac{[p]+1}{p},
\end{equation}
holds as a consequence of formula \eqref{EqRoughODEsEstimate}, for a constant $c$ depending only on the $V_i$. 

\begin{prop}
\label{PropSummaryPropertiesGeneralCase}
The family of maps $\big(\mu_{ts}\big)_{0\leq s\leq t\leq T}$ is a $\mcC^1$-approximate flow. 
\end{prop}

\medskip

\begin{Dem}
As the vector fields $V_i$ are of class $\mcC^{[p]+1}_b$, with $[p]+1\geq 3$, the identity 
$$
\big\|\mu_{ts}-\textrm{Id}\big\|_{\mcC^2} = o_{t-s}(1)
$$
holds as a consequence of classical results on ordinary differential equations; we turn to proving the $\mcC^1$-approximate flow property \eqref{EqMuMu}. Recall $X^m_{ts}$ stands for $\pi_m{\bfX}_{ts}$. We first use for that purpose formula \eqref{EqExactFormulaMuTs} to write
\begin{equation}
\label{EqMuMuGeneral}
\begin{split}
\mu_{tu}\big(\mu_{us}(x)\big) &= \Big(\textrm{F}^\otimes\big({\bfX}_{tu}\big)\textrm{Id}\Big)\big(\mu_{us}(x)\big) +\ep^{\textrm{Id}\,;\,[p]}_{tu}\big(\mu_{us}(x)\big) \\
&= \mu_{us}(x) + \sum_{m=1}^{[p]} \Big(\textrm{F}^\otimes\big(X^m_{tu}\big)\textrm{Id}\Big)\big(\mu_{us}(x)\big) +\ep^{\textrm{Id}\,;\,[p]}_{tu}\big(\mu_{us}(x)\big).
\end{split}
\end{equation}
We splitted the function $\textrm{F}^\otimes\big({\bfX}_{tu}\big)\textrm{Id}$ into a sum of functions $\textrm{F}^\otimes\big(X^m_{tu}\big)\textrm{Id}$ with different regularity properties, so one needs to use different Taylor expansions for each of them. One uses \eqref{EqEstimateC1SizeEpsilon} and inequality \eqref{EqRoughODEsEstimate} to deal with the remainder 
$$
\Big\|\ep^{\textrm{Id}\,;\,[p]}_{tu}\big(\mu_{us}(x)\big)\Big\|_{\mcC^1} \leq c\Big(1+\|{\bfX}\|^{[p]}\Big)^2|t-u|^\frac{[p]+1}{p}.
$$

\ssk

To deal with the term $\Big(\textrm{F}^\otimes\big(X^m_{tu}\big)\textrm{Id}\Big)\big(\mu_{us}(x)\big)$, we use formula \eqref{EqExactFormulaMuTs} with $n=[p]-m$ and $f=\textrm{F}^\otimes\big(X^m_{tu}\big)\textrm{Id}$. Writing $ds$ for $ds_{[p]-m}\dots ds_1$, we have
{\small
\begin{equation}
\label{EqControlFOtimes}
\Big(\textrm{F}^\otimes\big(X^m_{tu}\big)\textrm{Id}\Big)\big(\mu_{us}(x)\big) = \Big(\textrm{F}^\otimes\big(X^m_{tu}\big)\textrm{Id}\Big)(x) + \left(\textrm{F}^\otimes\Bigg(\Big\{\pi_{\leq [p]}\sum_{k=1}^{[p]-m} \frac{1}{k!}\,{\bf\Lambda}_{us}^{\bullet k}\Big\}\,X^m_{tu}\Bigg)\textrm{Id}\right) (x) + \ep^{\star\,;\,{p}-m}_{us}(x).
\end{equation}}
The notation $\star$ in the above identity stands for the $\mcC^{[p]+2-m}_b$ function $\textrm{F}^\otimes\big(X^m_{tu}\big)\textrm{Id}$; it has $\mcC^1$-norm controlled by \eqref{EqEstimateC1SizeEpsilon}. The result follows directly from \eqref{EqMuMuGeneral} and \eqref{EqControlFOtimes} writing
$$
\mu_{us}(x) = \Big(\textrm{F}^\otimes\big({\bfX}_{us}\big)\textrm{Id}\Big)(x) + \ep^{\textrm{Id}\,;\,[p]}_{us}(x),
$$
and using the identities $\exp\big({\bf \Lambda}_{us}\big)=\bfX_{us}$ and ${\bfX}_{ts} = {\bfX}_{us}{\bfX}_{tu}$ in $T^{[p]}_\ell$.
\end{Dem}

\medskip

\begin{defn}
\label{DefnGeneralRDESolution}
A \textbf{\emph{flow}} $(\varphi_{ts}\,;\,0\leq s\leq t\leq T)$ is said to \textbf{\emph{solve the rough differential equation}}
\begin{equation}
\label{RDEGeneral}
d\varphi = \textrm{\emph{F}}^\otimes\,\bfX(dt)
\end{equation}
if there exists a constant $a>1$ independent of $\bfX$ and two possibly $\bfX$-dependent positive constants $\delta$ and $c$ such that
\begin{equation}
\label{DefnSolRDEGeneral}
\|\varphi_{ts}-\mu_{ts}\|_\infty \leq c\,|t-s|^a
\end{equation}
holds for all $0\leq s\leq t\leq T$ with $t-s\leq\delta$.
\end{defn}

This definition can be equivalently reformulated in terms of uniform Taylor-Euler expansion of the form
$$
f\big(\varphi_{ts}(x)\big) = f(x) + \sum_{|I|\leq [p]} X^I_{ts}\big(V_If\big)(x) + O\big(|t-s|^{>1}\big).
$$
The following well-posedness result follows directly from theorem \ref{ThmConstructingFlows} and proposition \ref{PropSummaryPropertiesGeneralCase}; its proof is identical to the proof of theorem \ref{ThmLyonsUniversalLimitThm}, without a single word to be changed, except for the power of $\|\bfX\|$ in estimate \eqref{EqApproxPhiMu}, which needs to be taken as $[p]+1$ instead of $3$.

\begin{thm}
\label{ThmMainResultGeneral}
The rough differential equation 
$$
d\varphi = \textrm{\emph{F}}^\otimes\,\bfX(dt)
$$ 
has a unique solution flow; it takes values in the space of uniformly Lipschitz continuous homeomorphisms of $E$ with uniformly Lipschitz continuous inverses, and depends continuously on $\bfX$.
\end{thm}

\medskip

Remarks \ref{RemsThm2p3} on Friz-Victoir's approach to rough differential equations and time-inhomogeneous dynamics also hold in the general setting of this section.  
Section \ref{SubsubsectionClassicalRDEs} on classical rough differential equations has a straightforward analogue in the general setting of this section. We leave the reader the pleasure to adapt it and check that Lyons' universal limit theorem holds, with exactly the same proof as given in section \ref{SubsubsectionClassicalRDEs}. Note only that we ask the vector fields $V_i$ to be $\mcC^{[p]+1}_b$ when working with a weak geometric H\"older $p$-rough path; the drift vector field $V$ is only required to be $\mcC^2_b$.

\bigskip
\bigskip

\section{Exercices on flows driven by rough paths}
\label{SubectionExoRoughFlows}

We first see in exercice 12 how some general result on flows and $\mcC^1$-approximate flows proved in exercice 6 can be used to strengthen the result on the continuous dependence of the solution to a rough differential equation \wrt the driving rough path into a local Lipschitz dependence. The next two exercices are variations on the notion of solution flow. Roughly speaking, they are defined in terms of uniform Taylor-Euler expansion property. What happens if the driving vector fields allow for a priori higher order expansion? How robust are these expansions \wrt perturbation of the driving rough path? Exercices 13 and 14 partly answer these questions. Exercise 15 makes a crucial link between the 'differential formulation' of a solution path to a rough differential equation introduced above and the 'integral formulation' that can be set in the setting of controlled paths. The equivalence between these two formulations will be fundamental for the applications to stochastic analysis exposed in the last part of the course.

\bigskip

{\small 

{\bf 12. Local Lipschitz continuity of $\varphi$ \wrt $\bfX$.} Use the result proved in exercice 6 to prove that the solution flow to a rough differential equation driven by $\bfX$ is a locally Lipschitz continuous function of $\bfX$, in the uniform norm topology. \vspace{0.3cm}

{\bf 13. Taylor expansion of solution flows.} Let $V_1,\dots,V_\ell$ be $\mcC^{[p]+1}_b$ vector fields on a Banach space E, and $\bfX$ be a weak geometric H\"older $p$-rough path over $\RR^\ell$, with $2\leq p$. Set $\textrm{F}=\big(V_1,\dots,V_\ell\big)$. The solution flow to the rough differential equation 
$$
d\varphi = \textrm{F}^\otimes{\bfX}(dt)
$$
enjoys, by definition, a uniform Taylor-Euler expansion property, expressed either by writing
$$
\big\|\varphi_{ts}-\mu_{ts}\big\|_\infty \leq c|t-s|^a
$$
for the $\mcC^1$-approximate flow $\big(\mu_{ts}\big)_{0\leq s\leq t\leq 1}$ contructed in section \ref{SubsubsectionFromTaylorToRDEsGeneral}, or by writing
$$
\left\|f\circ\varphi_{ts} - \sum_{|I|\leq [p]} X^I_{ts}V_If\right\|_\infty \leq c|t-s|^a.
$$
\textit{What can we say if the vector fields $V_i$ are actually more regular than $\mcC^{[p]+1}_b$?}

\ssk

Assume $N\geq [p]+2$ is given and the $V_i$ are $\mcC^N_b$. Let $\bfY$ be the canonical lift of $\bfX$ to a $\frak{G}^N_\ell$-valued weak geometric H\"older $N$-rough path, given by Lyons' extension theorem proved in exercice 7. Let $\Theta_{ts}\in\frak{g}^N_\ell$ stand for $\log {\bfY}_{ts}$. For any $0\leq s\leq t\leq 1$, let $\nu_{ts}$ be the time 1 map associated with the ordinary differential equation
$$
\dot z_u = \textrm{F}^\otimes\big(\Theta_{ts}\big)(z_u), \quad 0\leq u\leq 1. \vspace{0.1cm}
$$

{\bf a)} Prove that $\nu_{ts}$ enjoys the following Euler expansion property. For any $f\in\mcC^{N+1}_b$ we have
\begin{equation}
\label{EqEulerExpansionNuTs}
\big\|f\circ\nu_{ts} - \textrm{F}^\otimes\big({\bfY}_{ts}\big)f\big\|_\infty \leq c|t-s|^\frac{N+1}{p},
\end{equation}
where the contant $c$ depends only on the $V_i$ and $\bfX$. \vspace{0.1cm}

{\bf b)} Prove that $\big(\nu_{ts}\big)_{0\leq s\leq t\leq 1}$ is a $\mcC^1$-approximate flow. \vspace{0.1cm}

{\bf c)} Prove that $\varphi_{ts}$ satisfies the high order Euler expansion formula \eqref{EqEulerExpansionNuTs}. \vspace{0.3cm}

{\bf 14. Perturbing the signal or the dynamics?} Let $2\leq p$ be given and $V_1,\dots,V_\ell$ be $\mcC^{[p]+1}_b$ vector fields on E. Let $\bfX$ be a weak geometric H\"older $p$-rough path over $\RR^\ell$, and ${\bfa}\in\frak{g}^{[p]}_\ell$ be such that $\pi_j{\bfa} = 0$ for all $j\leq [p]-1$. Write it
$$
{\bfa} = \sum_{|I|=[p]} a^I {\be}_{[I]},
$$
where $\big(e_1,\dots,e_\ell\big)$ stand for the canonical basis of $\RR^\ell$, and for $I=\big(i_1,\dots,i_k\big)$,
$$
{\be}_{[I]} = \Big[e_{i_1},\big[\dots,\big[e_{i_{k-1}},e_{i_k}\big]\dots\Big]
$$
in $T^{[p]}_\ell$. The ${\be}_{[I]}$'s form a basis of $\frak{g}^{[p]}_\ell$ with $\pi_n{\be}_{[I]} = 0$ if $n\neq |I|$. Recall the definition of $\exp : T^{[p],0}_\ell\rightarrow T^{[p],1}_\ell$ and its reciprocal $\log$. \vspace{0.2cm}

{\bf a)} Show that one defines a weak geometric H\"older $p$-rough path $\overline\bfX$ over $\RR^\ell$ setting 
$$
\overline\bfX_{ts} = \exp\Big(\log {\bfX}_{ts} + (t-s)\bfa\Big). \vspace{0.1cm}
$$

{\bf b)} Show that the solution flow to the rough differential equation 
$$
d\psi = \textrm{F}^\otimes\,\overline{\bfX}(dt)
$$
coincides with the solution flow to the rough differential equation 
$$
d\varphi = Vdt + \textrm{F}^\otimes\,{\bfX}(dt),
$$
where the vector field $V$ is defined by the formula
$$
V = a^IV_{[I]}. \vspace{0.3cm}
$$

{\bf 15. Differential and integral formulations of a rough differential equation.} Recall the setting of controlled paths investigated in exercise 11, and let $\textrm{F}=\big(V_1,\dots,V_\ell\big)$ be a collection of $\mcC^3_b$ vector fields on $\RR^d$, seen as a linear map from $\RR^\ell$ to the set of vector fields on $\RR^d$. Let also $\bfX$ be a weak geometric H\"older $p$-rough path, with $2\leq p<3$. Prove that the path $(x_t)_{0\leq t\leq 1}$ is a solution path to the rough differential equation 
$$
dx_t = \textrm{F}(x_t)\,{\bfX}(dt),
$$
in the sense of definition \ref{DefnPathSolRDE}, if and only if it is a solution to the integral equation
$$
x_t = x_0 + \int_0^t \textrm{F}(x_s)\,d{\bfX}_s
$$
(in the set of $\RR^d$-valued paths controlled by $\bfX$), where the above integral is the rough integral defined in exercise 11. Note here that the setting of controlled paths offers the possibility to define the above integral for non-weak geometric H\"older $p$-rough paths, so one can also define fixed-point problems and solve rough differential equations in that setting. (See the excellent forthcoming lecture notes \cite{FrizHairer} for this point of view on rough differential equations.) Given what we have done in section \ref{SubsectionControlledPath} on controlled paths, we are bound however to working with H\"older $p$-rough paths, with $2<p<3$. Fortunatley, this will be sufficent to deal with stochastic differential equations driven by Brownian motion in the next part of the course.

\vfill
\pagebreak

\chapter{Applications to stochastic analysis}
\label{SectionFlowsRoughPaths}

\medskip

\todo[inline, backgroundcolor=white, bordercolor=black]{Guide for this chapter}

\medskip

So far, I have presented the theory of rough differential equations as a purely deterministic theory of differential equations driven by multi-scale time indexed signals. Lyons, however, constructed his theory first as a deterministic alternative to It\^o's integration theory, after some hints by F\"ollmer in the early 80's that It\^o's formula can be understood in a deterministic way, and other works (by Bichteler, Karandikhar...) on the pathwise construction of stochastic integrals. (Recall that stochastic integrals are obtained as limits \textit{in probability} of Riemann sums, with no hope for a stronger convergence to hold as a rule.) Lyons was not only looking for a deterministic way of constructing It\^o integrals, he was also looking for a way of obtaining them as \textit{continuous} functions of their integrator! This required a notion of integrator different from the classical one... Rough paths were born as such integrators, with the rough integral of controlled integrands, defined in exercise 13, in the role of It\^o integrals. What links these two notions of integrals is the following fudamental fact. \textit{Brownian motion has a natural lift into a H\"older $p$-rough path}, for any $2<p<3$, called the \textbf{Brownian rough path}. This object is constructed in section \ref{SubsectionBrownianRoughPath} using Kolmogorov's classical regularity criterion, and used in section \ref{SubsectionRoughStochIntegrals} to see that the stochastic and rough integrals coincide whenever they both make sense. This fundamental fact is used in section \ref{SubsectionRoughStochasticDE} to see that stochastic differential equations can be solved in a two step process.
\begin{enumerate}
   \item[{\bf (i)}] {\bf Purely probabilistic step.} Lift Brownian motion into the Brownian rough path.   
   \item[{\bf (ii)}] {\bf  Purely deterministic step.} Solve the rough differential equation associated with the stochastic differential equation.
\end{enumerate}
This requires from the driving vector fields to be $\mcC^3_b$, for the machinery of rough differential equations to make sense, which is more demanding than the Lipschitz regularity required in the It\^o setting. This constraint comes with an enormous gain yet: the solution path to the stochastic differential equation is now a \textit{continuous} function of the driving Brownian rough path, this is Lyons' universal limit theorem, in striking contrast with the measurable character of this solution, when seen as a function of Brownian motion itself. (The twist is that the second level of the Brownian rough path is itself just a measurable function of the Brownian path.) Together with the above solution scheme for solving stochastic differential equations, this provides a simple and deep understanding of some fundamental results on diffusion processes, as section \ref{SubsectionFreidlinWentzell} on Freidlin-Wentzell theory of large deviation will demonstrate.

\ssk

We follow the excellent forthcoming lecture notes \cite{FH14} in sections \ref{SubsectionRoughStochIntegrals} and \ref{SubsectionRoughStochasticDE}.

\medskip

\section{The Brownian rough path}
\label{SubsectionBrownianRoughPath}

\subsection{Definition and properties}
\label{SubsubsectionBrownianRPDefn}

Let $\big(B_t\big)_{0\leq t\leq 1}$ be an $\RR^\ell$-valued Brownian motion defined on some probability space $(\Omega,\mcF,\PP)$. There is no difficulty in using It\^o's theory of stochastic integrals to define the two-index continuous process
\begin{equation}
\label{EqDefnIteratedIto}
\bB_{ts}^\textrm{It\^o} := \int_s^t\int_s^u dB_r\otimes dB_u = \int_s^t B_{us}\otimes dB_u.
\end{equation}
This process satisfies Chen's relation 
$$
\bB_{ts}^\textrm{It\^o} = \bB_{tu}^\textrm{It\^o} + \bB_{us}^\textrm{It\^o} + B_{us}\otimes B_{tu}
$$
for any $0\leq s\leq u\leq t\leq 1$. As is $B$ well-known to have almost-surely $\frac{1}{p}$-H\"older continuous sample paths, for any $p>2$, the process 
$$
{\bfB}^\textrm{It\^o} = \big(B,\bB^\textrm{It\^o}\big)
$$
will appear as a H\"older $p$-rough path if one can show that $\bB^\textrm{It\^o}$ is almost-surely $\frac{2}{p}$-H\"older continuous. This can be done easily using \textit{Kolmogorov's regularity criterion}, which we recall and prove for completeness. Denote  for that purpose by $\DD$ the set of dyadic rationals in $[0,1]$ and write $\DD_n$ for $\big\{k2^{-n}\,;\,k=0..2^n\big\}$.

\begin{thm}[Kolmogorov's criterion]
\label{ThmKolmogorovCriterion}
\index{Kolmogorov's regularity criterion}
Let $(S,d)$ be a metric space, and $q\geq 1$ and $\beta>1/q>0$ be given. Let also $\big(X_t\big)_{t\in \DD}$be an $S$-valued process defined on some probability space, such that one has
\begin{equation}
\label{EqKolmogorovCriterion}
\big\|d(X_t,X_s)\big\|_{\LL^q} \leq C\,|t-s|^\beta,
\end{equation}
for some finite constant $C$, for all $s,t\in\DD$. Then, for all $\al\in \bigl[0,\beta-\frac{1}{q}\bigr)$, there exists a \emph{random variable} $\textrm{\emph{C}}_\al\in \LL^q$ such that one has almost-surely
$$
d\big(X_s,X_t\big) \leq \textrm{\emph{C}}_\alpha\,|t-s|^\alpha,
$$
for all $s,t\in\DD$; so the process $X$ has an $\al$-H\"older modification defined on $[0,1]$.
\end{thm}

\begin{Dem}
{\small 
Given $s,t\in\DD$ with $s<t$, let $m\geq 0$ be the only integer \st $2^{-(m+1)}\leq t-s < 2^{-m}$. The interval $[s,t)$ contains at most one interval $\bigl[r_{m+1},r_{m+1}+2^{-(m+1)}\bigr)$ with $r_{m+1}\in\DD_{m+1}$. If so, each of the intervals $[s,r_{m+1})$ and $\bigl[r_{m+1}+2^{-(m+1)},t\bigr)$ contains at most one interval $\bigl[r_{m+2},r_{m+2}+2^{-(m+2)}\bigr)$ with $r_{m+2}\in\DD_{m+2}$. Repeating this remark up to exhaustion of the dyadic interval $[s,t)$ by such dyadic sub-intervals, we see, using the triangle inequality, that 
$$
d\big(X_t,X_s\big) \leq 2\sum_{n\ge m+1}S_n,
$$
where $S_n=\sup_{u\in \DD_n} d\big(X_u,X_{u+2^{-n}}\big)$. So we have
$$
\frac{d\big(X_t,X_s\big) }{(t-s)^\al}\leq 2\sum_{n\geq m+1}S_n\,2^{(m+1)\al}\leq \textrm{\emph{C}}_\al
$$
where $\textrm{\emph{C}}_\al := 2\sum_{n\geq 0}2^{n\al}S_n$. But as the assumption \eqref{EqKolmogorovCriterion} implies
$$
\EE\big[S_n^q\big] \leq \EE\left[\sum_{t\in \DD_n} d\big(X_t,X_{t+2^{-n}}\big)^q\right]\leq 2^nC(2^{-n})^{q\,\beta},
$$
we have
$$
\big\|C_\alpha\big\|_{\LL^q} \leq 2\sum_{n\geq 0}2^{n\al}\|S_n\|_q \leq 2C\sum_{n\geq 0}2^{\bigl(\al-\beta+\frac{1}{p}\bigr)n}<\infty,
$$
so $C_{\al}$ is almost-surely finite. The conclusion follows in a straightforward way.
}
\end{Dem}

\medskip

Recall the definition of the homogeneous norm on $T^{2,1}_\ell$ 
$$
\|{\bfa}\| = \big\|1\oplus a^1\oplus a^2\big\| = \big|a^1\big|+ \sqrt{\big|a^2\big|}
$$
introduced in equation \eqref{EqDefnHomogeneousNormAmbiantSpace}, with its associated distance function $d({\bfa},{\bfb}) = \|{\bfa}^{-1}{\bfb}\|$. To see that $\bfB^\textrm{It\^o}$ is a H\"older $p$-rough path we need to see that it is almost-surely $\frac{1}{p}$-H\"older continuous as a $\big(T^{2,1}_\ell,\|\cdot\|\big)$-valued path. This can be obtained from Kolmologorv's criterion provided one has
$$
\big\|{\bfB}^\textrm{It\^o}_{ts}\big\|_{\LL^q} \leq C\,|t-s|^\frac{1}{2},
$$
for some constants $q$ with $0<\frac{1}{2} - \frac{1}{q}<\frac{1}{p}$, and $C$. Given the form of the norm on $T^{2,1}_\ell$, this is equivalent to requiring
\begin{equation}
\label{EqConditionBrownianRoughPath}
\big\|B_{ts}\big\|_{\LL^q} \leq C\,|t-s|^\frac{1}{2}, \quad \quad \big\|{\bB}^\textrm{It\^o}_{ts}\big\|_{\LL^\frac{q}{2}} \leq C\,|t-s|.
\end{equation}
These two inequalities holds as a straightforward consequence of the scaling property of Brownian motion. (The random variable ${\bB}^\textrm{It\^o}_{10}$ is in any $\LL^q$ as a consequence of the BDG inequality for instance.)

\begin{cor}
The process $\bfB^\textrm{\emph{It\^o}}$ is almost-surely a H\"older $p$-rough path, for any $p$ with $\frac{1}{3}<\frac{1}{p}<\frac{1}{2}$. It is called the \emph{\textbf{It\^o Brownian rough path}}. \index{Brownian rough path}
\end{cor}

Note that $\bfB^\textrm{It\^o}$ is not weak geometric as the symmetric part of $\bB^\textrm{It\^o}_{ts}$ is equal to $\frac{1}{2}\,B_{ts}\otimes B_{ts}-\frac{1}{2}(t-s)\,\textrm{Id}$. Note also that we may as well have used Stratonovich integral in the definition of the iterated integral
$$
\bB_{ts}^\textrm{Str} := \int_s^t\int_s^u \circ\,dB_r\otimes \circ dB_u = \int_s^t B_{us}\otimes \circ dB_u;
$$
this does not make a big difference a priori since
$$
\bB_{ts}^\textrm{Str} = \bB_{ts}^\textrm{It\^o} + \frac{1}{2}(t-s)\textrm{Id}.
$$
So one can define another H\"older $p$-rough path ${\bfB}_{ts}^\textrm{Str} = \big(B_{ts},\bB_{ts}^\textrm{Str}\big)$ above Brownian motion, called \textbf{Stratonovich Brownian rough path}. Unlike It\^o Brownian rough path, it is weak geometric. (Compute the symmetric part of $\bB_{ts}^\textrm{Str}$!) Whatever choice of Brownian rough path we do, its definition seems to involve It\^o's theory of stochastic integral. It will happen to be important for applications these two rough paths can actually be constructed in a pathwise way from the Brownian path itself.

\medskip

Given $n\geq 1$, define on the ambiant probability space the $\sigma$-algebra $\mcF_n := \sigma\big\{B_{k2^{-n}}\,;\,0\leq k\leq 2^n\big\}$, and let $B^{(n)}_\bullet$ stand for the continuous piecewise linear path that coincides with $B$ at dyadic times in $\DD_n$ and is linear in between. Denote by $B^{(n),i}$ the coordinates of $B^{(n)}$. There is no difficulty in defining 
$$
\bB^{(n)}_{ts} := \int_s^tB^{(n)}_{us}\otimes dB^{(n)}_u
$$
as a genuine integral as $B^{(n)}$ is piecewise linear, and one has acutally, for $j\neq k$,
\begin{equation}
\label{EqMartingaleApprox}
B^{(n)}_{ts} = \EE\big[B_{ts}\big|\mcF_n\big], \quad\quad \bB^{(n),jk}_{ts} = \EE\big[\bB^{\textrm{Str}, jk}_{ts} \big|\mcF_n\big]
\end{equation}
and $\bB^{(n),ii}_{ts} = \frac{1}{2}\,\Big(B^{(n),i}_{ts}\Big)^2$.

\begin{prop}
\label{PropApproximateBrownianRP}
The H\"older $p$-rough path ${\bfB}^{(n)} = \big(B^{(n)},\bB^{(n)}\big)$ converges almost-surely to ${\bfB}^\textrm{\emph{Str}}$ in the H\"older $p$-rough path topology.
\end{prop}

\begin{Dem}
We use the interpolation result stated in proposition \ref{PropInterpolation} to prove the above convergence result. The almost-sure pointwise convergence follows from the martingale convergence theorem applied to the martingales in \eqref{EqMartingaleApprox}. To get the almost-sure uniform bound
\begin{equation}
\label{EqUniformRPEstimate}
\sup_n\,\big\|{\bfB}^{(n)}\big\| <\infty
\end{equation}
it suffices to notice that the estimates
$$
\big|B_{ts}\big| \leq C_p |t-s|^\frac{1}{p}, \quad\quad \big|\bB^{\textrm{Str},jk}_{ts}\big| \leq C^2_p |t-s|^\frac{2}{p}
$$
obtained from Kolmogorov's regularity criterion with $C_p\in\LL^q$ for (any) $q>2$, give
$$
\big|B^{(n)}_{ts}\big| \leq \EE\big[C_p\big|\mcF_n\big]\,|t-s|^\frac{1}{p}, \quad\quad \big|\bB^{(n),jk}_{ts}\big| \leq \EE\big[C^2_p\big|\mcF_n\big] |t-s|^\frac{2}{p},
$$
so the uniform estimate \eqref{EqUniformRPEstimate} follows from Doob's maximal inequality, which implies that almost-sure finite character of the maximum of the martingales $\EE\big[C_p^{1 \textrm { or } 2}\big|\mcF_n\big]$, since this maximum is integrable.
\end{Dem}

\medskip

\subsection{How big is the Brownian rough path?}
\label{SubsubsectionHowBigBrownianRP}

The upper bound of $\big\|{\bfB}^{\textrm{It\^o}}\big\|_{\frac{1}{p}}$ provided by the constant $C_{\frac{1}{p}}$ of Kolmogorov's regularity result says us that $\big\|{\bfB}^{\textrm{It\^o}}\big\|_{\frac{1}{p}}$ is in all the $\LL^q$ spaces. The situation is actually much better! As a first hint, notice that since ${\bfB}^{\textrm{It\^o}}_{ts}$ has the same distribution as $\delta_{\sqrt{t-s}}{\bfB}^{\textrm{It\^o}}_{10}$, and the norm of ${\bfB}^{\textrm{It\^o}}_{10}$ has a Gaussian tail (this is elementary), we have
\begin{equation}
\label{EqGaussianTailB10}
\EE\left[\exp\left(\frac{\big\|{\bfB}^{\textrm{It\^o}}_{ts}\big\|^2}{t-s}\right)\right] = \EE\Big[\exp\Big(\big\|{\bfB}^{\textrm{It\^o}}_{10}\big\|^2\Big)\Big] < \infty
\end{equation}
The following Besov embedding is useful in estimating the H\"older norm of a path from its two-point moments. \index{Besov's embedding}

\begin{thm}[Besov]
Given $\alpha\in \big[0,\frac{1}{2}\big)$ there exists an integer $k_\alpha$ and a positive constant $C_\alpha$ with the following properties. For any metric space $(S,d)$ and any $S$-valued continuous path $(x_t)_{0\leq t\leq 1}$ we have
$$
\|x_\bullet\|_\alpha \leq C_\alpha\left(\int_0^1\int_0^1\left(\frac{d(x_t,x_s)}{\sqrt{t-s}}\right)^{2k}\,ds\,dt\right)^{\frac{1}{2k}}.
$$
\end{thm}
It can be proved as a direct consequence of the famous Garsia-Rodemich-Rumsey lemma. Applied to the Brownian rough path $\bfB^{\textrm{It\^o}}$, Besov's estimate gives
$$
\EE\Big[\big\|{\bfB}^{\textrm{It\^o}}\big\|^{2k}_{\frac{1}{p}}\Big] \leq C^{2k}_{\frac{1}{p}}\int_0^1\int_0^1 \EE\left[\left(\frac{\big\|{\bfB}_{ts}\big\|}{\sqrt{t-s}}\right)^{2k}\right]\,ds\,dt =  C^{2k}_{\frac{1}{p}}\,\EE\Big[\|{\bfB}_{10}\big\|^{2k}\Big].
$$
So it follows from \eqref{EqGaussianTailB10} that we have for any positive constant $c$ 
$$
\EE\left[\sum_{k\geq k_{\frac{1}{p}}} \frac{c^k\big\|{\bfB}^{\textrm{It\^o}}\big\|^{2k}_{\frac{1}{p}}}{k!} \right] \leq \EE\Big[\exp\big(cC^2_{\frac{1}{p}} \|{\bfB}^{\textrm{It\^o}}_{10}\big\|^2\big)\Big]
$$
so $\exp \|{\bfB}^{\textrm{It\^o}}\big\|^2_{\frac{1}{p}}$ will be integrable provided $c$ is small enough, by \eqref{EqGaussianTailB10}.

\begin{cor}
\label{CorGaussianTailNormBrownianRP}
The $p$-rough path norm of the Brownian rough path has a Gaussian tail.
\end{cor}

\section{Rough and stochastic integral}
\label{SubsectionRoughStochIntegrals}

Let $\bfX$ be any $\RR^\ell$-valued H\"older $p$-rough path, with $2<p<3$. Recall a linear map $A$ from $\RR^\ell$ to $\RR^d$ acts on $(\RR^\ell)^{\otimes 2}$ as follows: $A(a\otimes b) = (Aa)\otimes b$. Recall also that we defined in section \ref{SubsectionControlledPath} the integral of an $\textrm{L}(\RR^\ell,\RR^d)$-valued path $\big(\textrm{F}_s\big)_{0\leq s\leq 1}$ controlled by ${\bfX} = (X,\bbX)$ as the well-defined limit
$$
\int_0^1 \textrm{F}\,d{\bfX} = \lim\,\sum\,\textrm{F}_{t_i}X_{t_{i+1}t_i} + \textrm{F}'_{t_i}{\bbX}_{t_{i+1}t_i},
$$
where the sum is over the times $t_i$ of finite partitions $\pi$ of $[0,1]$  whose mesh tends to $0$. This makes sense in particular for $\bfX = \bfB^\textrm{It\^o}$. At the same time, if F is adapted to the Brownian filtration, the Riemann sums $\sum\,\textrm{F}_{t_i}B_{t_{i+1}t_i}$ converge in probability to the stochastic integral $\int_0^1 \textrm{F}_s\,dB_s$, as the mesh of the partition $\pi$ tends to $0$. Taking subsequences if necessary, one defines simultaneously the stochastic and the rough integral on an event of probability $1$. They actually coincide almost-surely if F' is adapted to the Brownian filtration! To see this, it suffices to see that $\sum \textrm{F}'_{t_i}{\bbX}_{t_{i+1}t_i}$ converges in $\LL^2$ to $0$ along the subsequence of partitions used to define the stochastic integral $\int_0^1 \textrm{F}_sdB_s$. Assume first that F' is bounded, by $M$ say. Then, since it is adapted and $\textrm{F}'_{t_i}$ is independent of $\bB_{t_{i+1}t_i}$, an elementary conditioning gives
$$
\Big\|\sum \textrm{F}'_{t_i}\bB^{\textrm{It\^o}}_{t_{i+1}t_i}\Big\|^2_{\LL^2} = \sum \Big\|\textrm{F}'_{t_i}\bB^{\textrm{It\^o}}_{t_{i+1}t_i}\Big\|^2_{\LL^2} \leq M^2\sum \Big\|\bB_{t_{i+1}t_i}\Big\|^2_{\LL^2} \leq M^2 |\pi|, 
$$
which proves the result in that case. If F' is not bounded, we use a localization argument and stop the process at the stopping time
$$
\tau_M := \inf \big\{u\in [0,1]\,;\,|\textrm{F}'_u|>M\big\}\wedge 1.
$$
The above reasoning shows in that case that we have the almost-sure equality
$$
\int_0^{\tau_M}\textrm{F}\,d{\bfB} = \int_0^1 \textrm{F}^{\tau_M}_s\,dB_s,
$$
from which the result follows since $\tau_M$ tends to $\infty$ as $M$ increases indefinitely.

\begin{prop}
\label{PropRoughItoIntegrals}
Let $\big(\textrm{\emph{F}}_s\big)_{0\leq s\leq 1}$ be an $\textrm{\emph{L}}(\RR^\ell,\RR^d)$-valued path controlled by ${\bfB}^{\textrm{\emph{It\^o}}} = (B,\bB)$, adapted to the Brownian filtration, with a derivative process $\textrm{\emph{F}}'$ also adapted to that filtration. Then we have almost-surely
$$
\int_0^1\textrm{\emph{F}}\,d{\bfB}^{\textrm{\emph{It\^o}}} = \int_0^1 \textrm{\emph{F}}_s\,dB_s.
$$
\end{prop}

\ssk

If one uses $\bfB^{\textrm{Str}}$ instead of ${\bfB}^{\textrm{It\^o}}$ in the above rough integral, an additional well-defined term
$$
(\star) := \underset{|\pi|\searrow 0}{\lim}\;\sum\textrm{F}'_{t_i}\,\frac{1}{2}\big(t_{i+1}-t_i\big)\textrm{Id}
$$
appears in the left hand side, and we have almost-surely
$$
\int_0^1\textrm{F}\,d{\bfB}^{\textrm{Str}}  = \int_0^1\textrm{F}\,d{\bfB}^{\textrm{It\^o}}  + (\star) = \int_0^1 \textrm{F}_s\,dB_s + (\star).
$$
To identify that additional term, denote by $\textrm{Sym}(A)$ the symmetric part of a matrix $A$ and recall that 
$$
\frac{1}{2}\big(t_{i+1}-t_i\big)\textrm{Id} = \textrm{Sym}\big({\bfB}^{\textrm{Str}}_{t_{i+1}t_i}\big) - \textrm{Sym}\big({\bfB}^{\textrm{It\^o}}_{t_{i+1}t_i}\big) = \frac{1}{2}\,\bB^{\otimes 2}_{t_{i+1}t_i} - \textrm{Sym}\big({\bfB}^{\textrm{It\^o}}_{t_{i+1}t_i}\big);
$$
note also that the above reasoning showing that $\sum \textrm{F}'_{t_i}\bB^{\textrm{It\^o}}_{t_{i+1}t_i}$ converges to $0$ in $\LL^2$ also shows that  $\sum \textrm{F}'_{t_i}\textrm{Sym}\Big(\bB^{\textrm{It\^o}}_{t_{i+1}t_i}\Big)$ converges to $0$ in $\LL^2$. So $(\star)$ is almost-surely equal to the limit as $|\pi|\searrow 0$ of the sums 
$$
\frac{1}{2}\sum \textrm{F}'_{t_i} \bB^{\otimes 2}_{t_{i+1}t_i}.
$$
Since
$$
\textrm{F}'_{t_i} \bB_{t_{i+1}t_i} = \textrm{F}_{t_{i+1}t_i} + \textrm{R}_{t_{i+1}t_i}
$$
for some $\frac{2}{p}$-H\"older remainder R, the above sum equals
$$
\frac{1}{2}\Big(\sum \textrm{F}_{t_{i+1}t_i} \bB_{t_{i+1}t_i}\Big) + o_{|\pi|}(1).
$$
We recognize in the right hand side sum a quantity which converges in probability to the bracket of F and $B$.

\begin{cor}
\label{CorRoughStochStrIntegral}
Under the assumptions of proposition \ref{PropRoughItoIntegrals}, we have almost-surely
$$
\int_0^1\textrm{\emph{F}}\,d{\bfB}^{\textrm{\emph{Str}}} = \int_0^1 \textrm{\emph{F}}_s\,\circ dB_s.
$$
\end{cor}

\medskip

\section{Rough and stochastic differential equations}
\label{SubsectionRoughStochasticDE}

Equipped with the preceeding two results, it is easy to see that the solution path to a rough differential equation driven by $\bfB^{\textrm{Str}}$ or $\bfB^{\textrm{It\^o}}$ coincides almost-surely with the solution of the corresponding Stratonovich or It\^o stochastic differential equation.

\begin{thm}
\label{ThmRoughStochDE}
Let $\textrm{\emph{F}} = \big(V_1,\dots,V_\ell\big)$ be $\mcC^3_b$ vector fields on $\RR^d$. The solution to the rough differential equation 
\begin{equation}
\label{EqRDESDE}
dx_t = \textrm{\emph{F}}(x_t)\,{\bfB}^{\textrm{Str}}(dt)
\end{equation}
coincides almost-surely with the solution to the Stratonovich differential equation
$$
dz_t = V_i(z_t)\,\circ dB^i_t.
$$
\end{thm}

A similar statement holds for the It\^o Brownian rough path and solution to It\^o equations.

\medskip

\begin{Dem}
Recall we have seen in exercise 16 that a path is a solution to the rough differential equation \eqref{EqRDESDE} if and only if it is a solution path to the integral equation 
$$
x_t = x_0 + \int_0^t \textrm{F}(x_s)\,d{\bfB}^{\textrm{Str}}_s.
$$
Given the result of corollary \ref{CorRoughStochStrIntegral}, the theorem will follow if we can see that $x_\bullet$ is adapted to the Brownian filtration; for if one sets $\textrm{F}_s := \textrm{F}(x_s)$ then its derivative $\textrm{F}'_s = D_{x_s}\textrm{F}'\textrm{F}(x_s)$, with $DF$ the differential of F with respect to $x$, will also be adapted. But the adaptedness of the solution $x_\bullet$ to equation \eqref{EqRDESDE} is clear from its construction in the proof of theorem \ref{ThmLyonsUniversalLimitThm}.
\end{Dem}

\ssk

We obtain as a corollary of theorem \ref{ThmRoughStochDE}, Lyons' universal limit theorem and the convergence result proved in proposition \ref{PropApproximateBrownianRP} for the rough path associated with the piecewise linear interpolation $B^{(n)}$ of $B$ the following fundamental result, first proved by Wong and Zakai in the mid 60'. \index{Wong-Zakai theorem}

\begin{cor}[Wong-Zakai theorem]
\label{CorWongZakai}
The solution path to the ordinary differential equation 
\begin{equation}
\label{EqODEWongZakai}
dx^{(n)}_t = \textrm{\emph{F}}\left(x^{(n)}_t\right)\,dB^{(n)}_t
\end{equation}
converges almost-surely to the solution path to the Stratonovich differential equation
$$
dx_t = \textrm{\emph{F}}(x_t)\,\circ dB_t.
$$
\end{cor}

\begin{Dem}
It suffices to notice that solving the rough differential equation
$$
dz^{(n)}_t = \textrm{F}\left(z^{(n)}_t\right)\,{\bfB}^{(n)}(dt)
$$
is equivalent to solving equation \eqref{EqODEWongZakai}.
\end{Dem}

\bigskip

\section{Freidlin-Wentzell large deviation theory}
\label{SubsectionFreidlinWentzell}

We shall close this course with a spectacular application of the continuity property of the solution map to a rough differential equation, by showing how one can recover the basics of Freidlin-Wentzell theory of large deviations for diffusion processes from a unique large deviation principle for the Stratonovich Brownian rough path. Exercise 18 also uses this continuity property to deduce Stroock-Varadhan's celebrated support theorem for diffusion laws from the corresponding statement for the Brownian rough path.

\subsection{A large deviation principle for the Stratonovich Brownian rough path}
\label{SubsubsectionLDPBrownianRP}

Let start this section by recalling Schilder's large deviation principle for Brownian motion. \medskip

{\bf a) Schilder's theorem.} L Define for that purpose the real-valued function I on $\mcC^0\big([0,1],\RR^d\big)$ equal to $\frac{1}{2}\|h\|^2_{H^1} = \frac{1}{2}\int_0^1 \big|\dot h_s\big|^2\,ds$ on $H^1$, and $\infty$ elsewhere. We agree to write $\textrm{I}(\mcA)$ for $\inf\{\textrm{I}\big(h_\bullet\big)\,;\,h\in\mcA\}$, for any Borel subset $\mcA$ of $\mcC^0\big([0,1],\RR^d\big)$, endowed with the $\mcC^0$ topology. \index{Schiler's theorem}

\begin{thm}
Let $\PP$ stand for Wiener measure on $\mcC^0\big([0,1],\RR^d\big)$ and $B$ stand for the coordinate process. Given any Borel subset $\mcA$ of $\mcC^0\big([0,1],\RR^d\big)$, we have
$$
-\textrm{\emph{I}}\Big(\overset{\circ}{\mcA}\Big) \leq \underline{\overline{\lim}}\;\varep^2\log\PP\big(\varep B_\bullet\in \mcA\big) \leq  -\textrm{\emph{I}}\big(\overline\mcA\big).
$$
\end{thm}

\begin{Dem}
The traditional proof of the lower bound is a simple application of the Cameron-Martin theorem. Indeed, if $\mcA$ is the ball of centre $h\in H^1$ with radius $\delta$, and if we define the probability $\QQ$ by its density
$$
\frac{d\QQ}{d\PP} = \exp\left(-\varep^{-1}\int_0^1h_sdB_s - \frac{\varep^{-2}}{2}\textrm{I}(h)\right)
$$
with respect to $\PP$, the process $\overline B_\bullet := B_\bullet - \varep^{-1}h$ is a Brownian motion under $\QQ$, and we have
\begin{equation*}
\begin{split}
\PP\big(|\varep B-h| \leq \delta\big) &= \PP\big(\big|\overline B\big| \leq \varep^{-1}\delta\big) = \EE_\QQ\left[{\bf 1}_{\big|\overline B\big| \leq \varep^{-1}\delta} \exp\left(-\varep^{-1}\int_0^1h_sdB_s - \frac{\varep^{-2}}{2}\textrm{I}(h)\right)\right] \\
&\geq e^{- \frac{\varep^{-2}}{2}\textrm{I}(h)}\,\QQ\Big(\big|\overline B\big| \leq \varep^{-1}\delta\Big) = e^{- \frac{\varep^{-2}}{2}\textrm{I}(h)}\big(1-o_\varep(1)\big).
\end{split}
\end{equation*}

\ssk

One classically uses three facts to prove the upper bound. \vspace{0.1cm}
\begin{enumerate}
   \item The piecewise linear approximation $B^{(n)}$ of $B$ introduced above obviously satisfies the upper bound, as $B^{(n)}$ lives (as a random variable) in a finite dimensional space where it defines a Gaussian random variable. \vspace{0.1cm}
   \item The sequence $\varep B^{(n)}_\bullet$ provides an exponentially good approximation of $\varep B_\bullet$, in the sense that 
   $$
   \limsup_{\varep\searrow 0}\,\varep^2\log\PP\Big(\big|\varep B^{(n)}-\varep B\big|_\infty \geq \delta\Big) \underset{m\rightarrow\infty}{\longrightarrow} -\infty.  \vspace{0.1cm}  
   $$
   \item The map I enjoys the following 'continuity' property. With $\mcA^\delta := \{x\,;/,d(x,\mcA)\leq \delta\}$, we have
   $$
   \textrm{I}\big(\overline\mcA\big) = \underset{\delta\searrow 0}{\lim}\,\textrm{I}\big(\mcA^\delta\big).   
   $$
\end{enumerate}

The result follows from the combination of these three facts. The first and third points are easy to see. As for the second, just note that $B^{(n)}-B$ is actually made up of $2^n$ independent copies of a scaled Brownian bridge $2^{-\frac{n+1}{2}}\,\overline B^k_\bullet$, with each $\overline B^k$ defined on the dyadic interval $\big[k2^{-n},(k+1)2^{-n}\big]$. As it suffices to look at what happens in each coordinate, classical and easy estimates on the real-valued Brownian bridge provide the result.
\end{Dem}

\bigskip

{\bf b) Schilder's theorem for Stratonovich Brownian rough path.} The extension of Schilder's theorem to the Brownian rough path requires the introduction of the function J, defined on the set of $\frak{G}^2_\ell$-valued continuous paths ${\be}_\bullet = \big(e^1_\bullet,e^2_\bullet\big)$ by the formula
$$
\overline{\textrm{J}}\big({\be}_\bullet\big) = \textrm{I}\big(e^1_\bullet\big).
$$
Recall the definition of the dilation $\delta_\lambda$ on $T^2_\ell$, given in \eqref{EqDefnDeltaLambda}. Given any $0\leq \frac{1}{p}<\frac{1}{2}$, one can see the distribution ${\bfP}_\varep$ of $\delta_\varep\bfB^{\textrm{Str}}$ as a probability measure on the space of $\frac{1}{p}$-H\"older $\frak{G}^2_\ell$-valued functions, with the corresponding norm.

\begin{thm}
\label{ThmSchilderBrownianRP}
The family ${\bfP}_\varep$ of probability measures on $\mcC^{\frac{1}{p}}\big([0,1],\frak{G}^2_\ell\big)$ satisfies a large deviation principle with good rate function $\overline{\textrm{\emph{J}}}$.
\end{thm}

It should be clear to the reader that it is sufficent to prove the claim for the Brownian rough path above a 2-dimensional Brownian motion $B = \big(B^1,B^2\big)$, defined on $\mcC^0\big([0,1],\RR^2\big)$ as the coordinate process. We shall prove this theorem as a consequence of Schilder's theorem; this would be straightforward if the second level process $\bB$ -- or rather just its anti-symmetric part -- were a continuous function of the Brownian path, in uniform topology, which does not hold true of course. However, proposition \ref{PropApproximateBrownianRP} on the approximation of the Brownian rough path by its 'piecewise linear' counterpart makes it clear that it is almost-surely equal to a limit of continuous functional of the Brownian path. So it is tempting to try and use the following general contraction principle for large deviations. (See the book \cite{Kallenberg} by Kallenberg for an account of the basics of the theory, and a proof of this theorem.) We state it here in our setting to avoid unnecessary generality, and define the approximated L\'evy area $\bbA^m_{ts}$ as a real-valued function on the $\mcC^0\big([0,1],\RR^2\big)$ setting
$$
\bbA^m_t := \frac{1}{2}\,\int_0^t \Big(B^1_{\frac{[ms]}{m}}\,dB^2_s - dB^1_s\,B^2_{\frac{[ms]}{m}}\Big).
$$
It is a continuous function of $B$ in the uniform toplogy. The maps $\bbA^m$ converge almost-surely uniformly to the L\'evy area process $\bbA_\bullet$ of $B$. We see the process $\bbA_\bullet$ as a map defined on the space $\mcC^0\big([0,1],\RR^2\big)$, equal to L\'evy's area process on a set of probability 1 and defined in a genuine way on $H^1$ using Young integrals. (Note that elements of $H^1$ are $\frac{1}{2}$-H\"older continuous.) \index{Extended contraction principle}

\begin{thm}(Extended contraction principle)
\label{ThmExtendedContractionPrinciple}
If
\begin{enumerate}
   \item {\bf (Exponentially good approximation property)} 
   $$
   \limsup_\varep\,\varep^2 \log {\bfP}_\varep\big(\big\|\bbA^m-\bbA\big\|_\infty>\delta\big)\,\underset{m\rightarrow \infty}{\longrightarrow} -\infty,
   $$
   \item {\bf (Uniform convergence on I-level sets)}  for each $r>0$ we have 
   $$
   \big\|(\bbA^m-\bbA)_{\big| \{\textrm{\emph{I}} \leq r\}}\big\|_{\infty} \underset{m\rightarrow \infty}{\longrightarrow}\,0,
   $$
\end{enumerate}
then the distribution of $\bbA_\bullet$ under ${\bfP}_\varep$ satisfies a large deviation principle $\mcC^0\big([0,1],\frak{g}^2_2\big)$ with good rate function $\inf\{\textrm{\emph{I}}(\omega)\,;\,\bfa = \bbA(\omega)\}$.
\end{thm}

\medskip

\begin{DemThmSchilderBrownianRP}
The proof amounts to proving points (1) and (2) in theorem \ref{ThmExtendedContractionPrinciple}. The second point is elementary if one notes that for $h\in H^1\big([0,1],\RR^2\big)$, we have
\begin{equation*}
\begin{split}
\left| \int_0^t \Big(h_{{\frac{[ms]}{m}}}-h_s\Big)\otimes dh_s\right| &\leq \|h\|_{\frac{1}{2}}\Big(\frac{1}{m}\Big)^\frac{1}{2}\|h\|_1 \\
&\leq \|h\|_{H^1}^2 m^{-\frac{1}{2}}.
\end{split}
\end{equation*}
As for the first point, it suffices to prove that
$$
\limsup_\varep\,\varep^2 \log{\bfP}_1\left(\sup_{t\in [0,1]}\, \int_0^t\Big(B^1_s - B^1_{{\frac{[ms]}{m}}}\Big)\,dB^2_s\geq \varep^{-2}\delta\right)\,\underset{m\rightarrow \infty}{\longrightarrow} -\infty,
$$
which we can do using elementary martingale inequalities. Indeed, denoting by $M_t$ the martingale defined by the above stochastic integral, with bracket $\int_0^t \Big| B^1_s - B^1_{{\frac{[ms]}{m}}}\Big|^2ds$, the classical exponential inequality gives
$$
{\bfP}_1\Big(M^*_1 \geq \delta\varep^{-2},\; \langle M\rangle_1\geq\varep^{-2}\,m^{-\frac{1}{p}}\Big) \leq  \exp\left(-\frac{\delta^2\varep^{-2}\,m^{\frac{1}{p}}}{2}\right),
$$
while we also have
$$
{\bfP}_1\Big(\langle M\rangle_1\geq\varep^{-2}\,m^{-\frac{1}{p}}\Big) \leq {\bfP}_1\Big(\big\|B^1\big\|^2_{\frac{1}{p}}\,m^{-\frac{2}{p}} \geq\varep^{-2}\,m^{-\frac{1}{p}}\Big)
$$
So the conclusion follows from the fact that the $\frac{1}{p}$-H\"older norm of $B^1$ has a Gaussian tail.
\end{DemThmSchilderBrownianRP}

\bigskip

\subsection{Freidlin-Wentzell large deviation theory for diffusion processes}
\label{SubsubsectionFWLDTheory}

All together, theorem on the rough path interpretation of Stratonovich differential equations, Lyons' universal limit theorem and the large deviation principle satisfied by the Brownian rough path prove the following basic result of Freidlin-Wentzell theory of large deviation for diffusion processes. Given some $\mcC^3_b$ vector fields $V_1,\dots,V_\ell$ on $\RR^d$, and $h \in H^1$, denote by $y^h$ the solution to the well-defined controlled ordinary differential equation
$$
dy^h_t = \varep V_i\big(y^h_t\big)\,\circ dh^i_t.
$$

\begin{thm}[Freidlin-Wentzell]
Denote by ${\bfP}_\varep$ the distribution of the solution to the Stratonovich differential equation \index{Freidlin-Wentzell large deviation theorem}
$$
dx_t = \varep V_i(x_t)\,\circ dB^i_t,
$$
started from some initial condition $x_0$. Given any $\frac{1}{p}<\frac{1}{2}$, one can consider ${\bfP}_\varep$ as a probability measure on $\mcC^\frac{1}{p}\big([0,1],\RR^d\big)$. Then the family ${\bfP}_\varep$ satisfies a large deviation principle with good rate function 
$$
\textrm{\emph{J}}(z_\bullet) = \inf\big\{\textrm{\emph{I}}(h)\,;\, y^h_\bullet = z_\bullet\big\}.
$$
\end{thm}

\bigskip

\section{Exercises on rough and stochastic analysis}
\label{SubsectionexercisesStochastic}

Exercise 16 provides another illustration of the power of Lyons universal limit theorem and the continuity of the solution map to a rough differential equation, called the \textit{It\^o map}. It shows how to obtain a groundbreaking result of Stroock and Varadhan on the support of diffusion laws by identifying the support of the distribution of the Brownian rough path. Exercise 17 gives an interesting example of a rough path obtained as the limit of a 2-dimensional signal made up of a Brownian path and a delayed version of it. While the first level concentrates on a degenerate signal with identical coordinates and null area process as a consequence, the second level converges to a non-trivial function. Last, exercise 18 is a continuation of exercise 11 on the pairing of two rough paths. \vspace{0.2cm}

{\bf 16. Support theorem for the Brownian rough path and diffusion laws.} \index{Support theorem of Stroock and Varadhan}  We show in this exercise how the continuity of the It\^o map leads to a deep result of Stroock and Varadhan on the support of diffusion laws. The reader unacquainted with this result may have a look at the poloshed proof given in the book by Ikeda and Watanabe \cite{IkedaWatanabe} to see the benefits of the rough path approach.

\medskip

{\bf a) Translating a rough path.} Given a Lipschitz continuous path $h$ and a $p$-rough path ${\bfa} = 1\oplus a^1\oplus a^2$, with $2<p<3$, check that we define another $p$-rough path  setting
$$
\tau_h({\bfa})_{ts} := 1 \oplus \big(a^1_{ts}+h_{ts}\big) \oplus \Big(a^2_{ts}+\int_s^t a^1_{us}\otimes dh_u+ \int_s^t h_{us}\otimes da^1_u +\int_s^t h_{us}\otimes dh_u\Big),
$$
where the integral $\int_s^t h_{us}\otimes da^1_u$ is defined as a Young integral by the integration be parts formula 
$$
\int_s^t h_{us}\otimes da^1_u := h_{ts}\otimes a^1_{ts} - \int_s^t dh_u\otimes a^1_{us}. \vspace{0.2cm}
$$

{\bf b)} Given any $\RR^\ell$-valued coninuous path $x_\bullet$, denote as in section \ref{SubsectionBrownianRoughPath} by $x^{(n)}$ the piecewise linear coninuous interpolationof $x_\bullet$ on dyadic times of order $n$, and let ${\bfX}^{(n)}$ stand for its associated rough path, for $2<p<3$. We define a map ${\bfX} : \mcC^0\big([0,1],\RR^\ell\big)\rightarrow \big({\frak{G}^{2,1}_\ell}\big)^{[0,1]}$ setting
$$
\pi_1\big({\bfX}(x_\bullet)\big) := x_\bullet, \quad \quad \pi_2\big({\bfX}(x_\bullet)\big)^{jk}_t := \limsup_n \int_0^tx^{(n),j}_u\otimes dx^{(n),k}_u.
$$
So the random variable ${\bfX}(x_\bullet)$ is almost-surely equal to Stratonovich Browian rough path under Wiener measure $\PP$. \vspace{0.1cm}

$\quad\quad${\bf (i)} Show that one has $\PP$-almost-surely 
$$
{\bfX}(x+h) = \tau_h\big({\bfX}(x)\big)
$$
for any Lipschitz continuous path $h$. \vspace{0.1cm}

$\quad\quad${\bf (ii)} Prove that th law of the random variable $\tau_h\circ{\bfX}$ is equivalent to the law of ${\bfX}$ under $\PP$. \vspace{0.1cm}

\noindent Recall that the support of a probability measure on a topological space if the smallest closed setof full measure. We consider $\bfX $, under $\PP$, as a $\mcC^{}\frac{1}{p}\big([0,1],\frak{G}^2_\ell\big)$-valued random variable. \vspace{0.1cm}

$\quad\quad${\bf (iii)} Prove that if ${\bfa}_\bullet$ is an element of the support of the law of $\bfX$ under $\PP$, then $\tau_h{\bfa}$ as well. \vspace{0.2cm}

{\bf c) (i)} Use the same kind of arguments as in proposition \ref{PropApproximateBrownianRP} to show that one can ind an element ${\bfa}_\bullet$ in the support of the law of $\bfX$ under $\PP$, and some Lipschitz coninuous paths $x^{(n)}_\bullet$ such that $\big\|\tau_{x^{(n)}}{\bfa}\big\|$ tends to 0 as $n\rightarrow \infty$. \vspace{0.1cm}

$\quad\quad$ {\bf (ii)} Prove that the support of the law of ${\bfB}_{\bullet 0}^{\textrm{Str}}$ in $\mcC^\frac{1}{p}\big([0,1],\frak{G}^{2}_\ell\big)$ is the closure in $\frac{1}{p}$ H\"older topology of the set of of Lipschitz continuous paths. \vspace{0.2cm}

{\bf d) Stroock-Varadhan support theorem.} Let $\bfP$ stand for the distribution of the solution to the rough differential equation in $\RR^d$
$$
dx_t = V_i(x_t)\,\circ dB^i_t,
$$
driven by Brownian motion and some $\mcC^3_b$ vector fields $V_i$. Justify that one can see $\bfP$ as a probability on $\mcC^\frac{1}{p}\big([0,1],\RR^d\big)$. Let also write $y^h$ for the solution to the ordinary differential equation 
$$
dy^h_t = V_i\big(y^h_t\big)\,dh^i_t
$$
 driven by a Lipschitz $\RR^\ell$-valued path $h$. Prove that the support of $\bfP$ is the closure in $\mcC^\frac{1}{p}\big([0,1],\RR^d\big)$ of the set of all $y^h$, for $h$ ranging in the set of Lipschitz $\RR^\ell$-valued paths. \vspace{0.3cm}

 {\bf 17. Delayed Brownian motion.} Let $(B_t)_{0\leq t\leq 1}$ be a real-valued Brownian motion. Given $\ep>0$, we define a 2-dimensional process setting
 $$
 x_t = \big(B_{t-\ep},B_t\big);
 $$
its area process 
 $$
 {\bbA}^\ep_{ts} :=\frac{1}{2}\,\int_s^t\Big(B_{u-\ep,s-\ep}dB_u - B_{us}dB_{u-\ep}\Big) 
 $$
 is well-defined for $0\leq t-s<\ep$, as $B_{\bullet-\ep}$ and $B_\bullet$ are independent  on $[s,t]$ in that case. \vspace{0.2cm}
 
 {\bf 1)} Show that we define a rough path ${\bfX}^\ep$ setting
 $$
 {\bfX}^\ep_t := \exp\big(x_t + {\bbA}^\ep_t\big)\in T^2_2. \vspace{0.1cm} 
 $$
 
 {\bf 2)} Recall that $d$ stands for the ambiant metric in $T^2_2$. Prove that one can find a positive constant $a$ such that the nequality
 $$
 \EE\left[\exp\left(a\,\frac{d\big({\bfX}^\ep_t,{\bfX}^\ep_s\big)^2}{t-s}\right)\right]\leq C<\infty 
 $$
 holds for a positive constant $C$ independent of $0<\ep\leq 1$ and $0\leq s\leq t\leq 1$. As in section \ref{SubsubsectionHowBigBrownianRP}, it follows from Besov's embedding theorem that, for any $2<p<3$, the weak geometric H\"older $p$-rough path ${\bfX}^\ep$ has a Gaussian tail, with 
 $$
 \underset{0<\ep\leq 1}{\sup}\,\EE\left[\exp\Big(a\big\|{\bfX}^\ep\big\|^2\Big)\right] < \infty 
 $$
 for some positive constant $a$. \vspace{0.1cm}
 
 {\bf 3)} Define ${\bf 1}$ as the vector of $\RR^2$ with coordinates $1$ and $1$ in the canonical basis, and set 
 $$
 {\bfY}_t := \exp\Big(B_t{\bf 1}-\frac{t}{2}\textrm{Id}\Big).
 $$
Write $\textrm{d}_p$ for the distance on the set of H\"older $p$-rough paths given in definition \ref{DefnRoughPath}. Prove that $\textrm{d}_p\big({\bfX}^\ep,{\bfY}\big)$ converges to 0 in $\LL^q$, for any $1\leq q<\infty$. \vspace{0.3cm}

{\bf 18. Joint lift of a random and a deterministic rough path.} Let $2<p<3$ and ${\bfX}=(X,\bbX)$ be an $\RR^d$-valued H\"older $p$-rough path. Denote by $\bfB$ the It\^o Brownian rough path over $\RR^\ell$. Given $j\in\llbracket 1,d\rrbracket$ and $k\in\llbracket 1,\ell\rrbracket$, tdefine he integral 
$$
{\bbZ}^{jk}_{ts} := \int_s^t X^k_{us}dB^j_u
$$ 
as a genuine It\^o integral, and define the integral $\int_s^t B^j_{us}dX^k_u$ by integration by part, setting
$$
{\bbZ}^{kj}_{ts} := \int_s^t B^j_{us}dX^k_u := B^j_{ts}X^k_{ts} - \int_s^t X^k_{us}dB^j_u.
$$
Prove that one defines a H\"older $p$-rough path $\bfZ$ over $(X,B)\in\RR^{d+\ell}$ defining the $(jk)$-component of its second order level, as equal to ${\bbX}^{jk}$ if $1\leq j,k\leq d$, equal to ${\bbB}^{jk}$ if $d+1\leq j,k\leq d+\ell$, and by the above formulas otherwise.

\vfill
\pagebreak

\chapter{Looking backward}
\label{SectionLookingBackward}

\section{Summary}

It is now time to forget the details and summarize the main ideas.

\medskip

Chapter 2 provides a toolbox for constructing flows on Banach spaces from approximate flows. The interest of working with this notion comes from the fact that approximate flows pop up naturally in a number of situations, more or less under the form of "numerical schemes", as the step-1 Euler scheme for ordinary differential equations encountered in exercise 1, or the higher order Milstein-type schemes $\mu_{ts}$ used in chapter 4 to define solutions to rough differential equations. Rough paths appear in that setting as coefficients in the numerical schemes, and as natural generalizations of multiple integrals in some H\"older scales space. The miracle that takes place here is essentially algebraic and rests on the fact that solving an ordinary differential equation is algebraically very close to an exponentiation operation. This echoes the fact that the tensor space $T^{[p]}_\ell$ in which rough paths live also has natural notions of exponential and logarithm. As a matter of fact, this "pairing" ODE-exponential-rough paths works in exactly the same way with the branched rough paths introduced in \cite{GubinelliBranched} by Gubinelli in order to deal with rough differential equations driven by non-weak geometric rough paths.

\medskip

We have concentrated in this course on one approach to rough differential equations and rough paths. There are other approaches, with their own benefits, to start with Lyons' original formulation of his theory, as exposed in Lyons' seminal article \cite{Lyons97} or his book \cite{LyonsQianBook} with Qian. Its core concept is a notion of rough integral which associates a rough path to another rough path. As you may guess, it can be shown to be a continuous functions of both its integrand and integrator. In that setting, as solution path to a rough differential equation is a fixed point to an integral equation \textit{in the space of rough paths}; it was first solved using a Picard iteration process. Two crucial features of Lyons' original formulation were spotted by Gubinelli and Davie. The first level of a solution path $x_\bullet$ to a rough differential equation locally look like the first level of $\bfX$, and it suffices to know $x_\bullet$ and $\bfX$  to get back the entire rough path solution to the rough differential equation in Lyons' sense. This led Gubinelli to the introduction of the notion of controlled paths, which have far reaching applications to difficult problems on stochastic partial differential equations (SPDEs), as illustrated in the recent and brilliant works of Gubinelli and his co-authors. It also was one of the seeds of Hairer's groundbreaking theory of regularity structures \cite{Hairer}, which enabled him to construct a robust solution theory for some important up to now ill-posed SPDEs coming from physics deep problems. The forthcoming lecture notes \cite{FH14} by Friz and Hairer provides a very nice introduction to Gubinelli's point of view on rough paths theory.

\ssk

On the other hand, Davie uncovered in \cite{Davie} the fact that one can characterize the first level of a solution path to a rough differential equation in Lyons' sense in terms of numerical schemes of Milstein type. Has was able in that setting to prove sharp well-posedness and existence results, under essentially optimal regularity conditions on the driving vector fields, for $p$-rough paths with $2\leq p<3$. No notion of integral is needed in this approach, and we only work with $\RR^d$-valued paths, as opposed to Lyons' formulation. 

\ssk

His ideas were reworked and generalized by Friz and Victoir \cite{FrizVictoirEuler}, who defined solution paths to rough differential equations driven by some rough path $\bfX$ as limits of solution paths to controlled ordinary differential equations, in which the rough path canonically associated with the (smooth or absolutely continuous) control converges in a rough paths sense to $\bfX$. The book \cite{FVBook} by Friz and Victoir provides a thorough account of their approach. See also the short 2009 lecture notes \cite{LNPeter2009} of Friz. (My presentation in section \ref{SubsectionFreidlinWentzell} of the material on Freidlin-Wentzell large deviation theory follows his approach.) The point of view presented in these notes builds on Davie's approach and on the inspiring work \cite{FdlP} of Feyel and de la Pradelle; it is mainly taken from the article \cite{BailleulFlows}.

\bigskip
\bigskip

I hope you enjoyed the tour.\footnote{Please do not hesitate to send me any comments or suggestions to the email address {\sf ismaelbailleul@univ-rennes1.fr}.}

\section{A guided tour of the litterature}

The litterature on rough paths theory is increasing rapidly. To help you find your way in this bush, I comment below on a few references whose reading may be helpful to get a better view of the domain. Reference numbers refer to the bibliography following this paragraph, while stared references refer to the above bibliography. \vspace{0.3cm}

\ssk

\noindent {\bf Lecture notes.} Here are collected a few references that aim at giving a pedagogical presentation of rough paths theory, from different point of views.

\begin{itemize}
   \item[-] Lejay gave in \cite{LejayYoung} a self-contained and easily accessible account of the theory of Young differential equations, which correspond to rough differential equations driven by $p$-rough paths, with $1<p<2$.  \vspace{0.2cm} 

   \item[-] The lecture notes \cite{IntroLejay2} by Lejay provides a well-motivated and detailled study of the algebraic setting in which rough paths theory \textit{needs} to be formulated. It is easily readable.\vspace{0.2cm}

   \item[-] The approach of Friz-Victoir to rough differential equations, as described in remark \ref{RemsThm2p3}, was put forward in \cite{FrizVictoirEuler}*; it is developped thouroughly in their monograph \cite{FVBook}*. The lecture notes \cite{LNPeter2009}* by Friz, and \cite{BaudoinLectureNotes}* by Baudoin, provide an easy access to that approach. \vspace{0.2cm}

   \item[-] We warmly recommend the reading of the forthcoming lecture notes \cite{FH14}* of Friz and Hairer on the theory of rough paths and rough differential equations seen from the point of view of controlled paths. Although it does not lend itself to an easy access when the roughness index $p$ is greater than 3, this approach is the seed of the very exciting development of a new framework for handling SPDEs which were previously untractable. Have a look for instance at the (hard) work \cite{HairerRegularity,HairerRegularityIntroduction} of Hairer on regularity structures, or the somewhat more "down-to-earth" work \cite{GubinelliImkellerPerkowski} of Gubinelli, Imkeller and Perkowski to see how ideas from controlled paths can enable you to do some forbidden operation: multiplying two distributions!  \vspace{0.2cm}
\end{itemize}

\ssk

\noindent {\bf Historical works.} I have chosen to put forward here a few references that illustrate the development of the theory.

\begin{itemize}
   \item[-] Lyons' amazing seminal work \cite{Lyons98} is a must.    \vspace{0.2cm}

   \item[-] One owes to Strichartz \cite{Strichartz} a far reaching generaiztion of the well-known \textit{Baker-Campbell-Dynkin-Hausdorff} formula expressing the multiplication in a Lie group as an operation in the Lie algebra. This fantastic paper was the basis of basis of groudbreaking works by Castell \cite{Castell}, Ben Arous and others on \textit{Taylor expansions} for stochastic differential equations, and can somehow be seen as a precursor to the approach to rough differential equations put forward in this course. \vspace{0.2cm}

   \item[-] The other paper that inspired our flow-based approach is the \textit{sewing lemma} proved in \cite{FdlP,FdlPM} by Feyel, de la Pradelle and Mokobodzki.   \vspace{0.2cm}

   \item[-] The work \cite{Davie}* by Davie showing that one could understand rough differential equations in terms of numerical schemes -- or Taylor expansions -- was also instrumental in the development of Friz-Victoir's approach to the subject, as well as to the present approach.  \vspace{0.2cm}
\end{itemize}

\medskip

\noindent {\bf Today.} The theory of rough paths is presently experiencing a fantastic development in all sorts of directions. Just a few of them: Malliavin related business, differential geometry and machine learning, to testify of the diversity of directions that are actively being investigated. \vspace{0.2cm}

\begin{itemize}
   \item[-] There has been much industry in proving that one can use Malliavin calculus methods and rough paths theory to show that solutions of rough differential equations driven by Gaussian rough paths have a density at any fixed time under some bracket-type conditions on the driving vector fields and some non-degeneracy conditions on the Gaussian noise. For two landmark results in this direction, see the works \cite{CassFriz} of Cass and Friz, and what may be a temporarily final point \cite{CassHairerLittererTindel} by Cass, Hairer, Litterer and Tindel. \vspace{0.2cm}
   
   \item[-] Rough path theory has an inherent geometric content built in, to start with the nilpotent Lie group on which rough paths live. Given that most all of physics takes place on finite or infinite dimensional manifolds (configuration spaces), it is natural to try and give an intrinsic notion of rough path on a manifold. Starting with the seminal work \cite{CassLittererLyons} of Cass, Litterer and Lyons that makes a first step in this direction, this important question is being investigated. See the work \cite{CassDriverLitterer} of Cass, Driver and Litterer, for a nice reworking of the ideas of \cite{CassLittererLyons}, giving an intrinsic notion of rough path on a compact finite dimensional manifold, and my own work \cite{RoughIntegrators} for a general framework for dealing with rough integrators on Banach manifolds. \vspace{0.2cm}

   \item[-] Lyons' group in Oxford is presently exploring the potential application of the use of signature (the set of all iterated integrals of a rough path) to investigate learning questions! See for instance the works \cite{LevinLyonsNi} by Levin, Lyons, Ni and \cite{GyurkoLyonsKontkowskiField}, by Gyurko, Lyons, Kontkowski and Field.
   
\end{itemize}


\end{document}